\documentclass[10pt, a4paper]{article}

\usepackage{amsmath, amsthm, amssymb, float}
\usepackage{graphicx,latexsym, amsmath,amsfonts}
\usepackage{amsthm}
\usepackage{hyperref}
\usepackage[nameinlink]{cleveref}
\usepackage{subcaption}
\usepackage{graphicx}
\usepackage{setspace}
\usepackage{makecell}
\usepackage{verbatim}
\usepackage{epstopdf}
\usepackage{color}
\usepackage{enumitem}
\usepackage{a4}
\usepackage[numbers]{natbib}
\allowdisplaybreaks[4]
\hypersetup{colorlinks=true}
\hypersetup{linkcolor=blue}

\numberwithin{equation}{section}
\theoremstyle{plain}

\newtheorem{definition}{Definition}[section]
\newtheorem{Theorem}{Theorem}[section]
 
\newtheorem{Lemma}{Lemma}[section]
\newtheorem{remark}{Remark}[section] 

\newtheorem{Example}{Example}[section]
\newtheorem{Assumption}{Assumption}[section]

\makeatletter \@addtoreset{equation}{section}

\begin{document}

\linespread{1.3}
\pagestyle{plain}

\title
{Strong convergence and Mittag-Leffler stability of stochastic theta method for time-changed stochastic differential equations}

\author{Jingwei Chen\footnotemark[2] \footnotemark[1], Jun Ye\footnotemark[3] , Jinwen Chen\footnotemark[3],  Zhidong Wang\footnotemark[4]}

\renewcommand{\thefootnote}{\fnsymbol{footnote}}
\footnotetext[1]{Corresponding author. (\texttt{chenj22@mails.tsinghua.edu.cn})}
\footnotetext[2]{Yau Mathematical Sciences Center, 
Tsinghua University,
Beijing, 100084, China.}
\footnotetext[3]{Department of Mathematical Sciences,
Tsinghua University,
Beijing, 100084, China.}
\footnotetext[4]{School of Economics and Management,
Tsinghua University,
Beijing, 100084, China.}

\date{}

\maketitle

\begin{abstract}
We propose the first $\alpha$-parameterized framework for solving time-changed stochastic differential equations (TCSDEs), explicitly linking convergence rates to the driving parameter of the underlying stochastic processes.
Theoretically, we derive exact moment estimates and exponential moment estimates of inverse $\alpha$-stable subordinator $E$ using Mittag-Leffler functions.
The stochastic theta (ST) method is investigated for a class of SDEs driven by a time-changed Brownian motion, whose coefficients are time-space-dependent and satisfy the local Lipschitz condition.
We prove that the convergence order dynamically responds to the stability index $\alpha$ of stable subordinator $D$, filling a gap in traditional methods that treat these factors independently.
We also introduce the notion of Mittag-Leffler stability for TCSDEs, and investigate the criterion of Mittag-Leffler stability for both the exact and numerical solutions.
Finally, some numerical simulations are presented to illustrate the theoretical results.


\medskip
\noindent \textbf{AMS Subject Clasification:} 65C30,\;60H10
\end{abstract}

\section{Introduction}
Stochastic differential equations (SDEs) driven by time-changed Brownian motion, also known as time-changed SDEs (TCSDEs), have emerged as a fundamental tool for modeling complex systems with trapping, waiting, or delay phenomena. 
These phenomena are ubiquitous in fields such as mathematical finance \cite{magdziarz2011option, clark1973subordinated}, physics \cite{scalas2006five}, biology \cite{kou2004generalized,saxton1997single}, hydrology \cite{benson2000application}.
For example, time-changed processes can describe constant stagnation of periods in financial markets. Magdziarz \cite{magdziarz2011option} introduced subdiffusion into pricing of options, and showed that time-changed Black-Scholes model performed better than the original model.
For the general theory and the triangular relationship between time-changed stochastic differential equations, time-fractional Fokker-Planck equations (FFPE) and subdiffusion, a thorough investigation can be found in \cite{umarov2018beyond}.

However, most TCSDEs encountered in practice do not admit closed-form solutions, necessitating the development of numerical approximation schemes.
Existing works have not applied the stochastic theta (ST) method to TCSDEs, and the case of local Lipschitz condition was not analyzed either.

Two challenging questions also arise:

\textbf{(Q1)} Does the stability index $\alpha$ intrinsically limit the strong convergence order of the numerical method?

\textbf{(Q2)} Could the traditional stability analysis be improved to explicitly incorporate the effect of $\alpha$? Specifically, can we develop a stability criterion that not only ensures boundedness but also captures the characteristic Mittag-Leffler decay governed by the stability index $\alpha$?

In this paper, we investigate the $\alpha$-scaled strong convergence of the ST method, the newly-defined Mittag-Leffler stability of both the exact and numerical solutions for the following TCSDEs with time-space-dependent coefficients
\begin{equation}
    X_{t}=x_{0}+\int_{0}^{t}F\left(s,X_{s}\right)dE_{s}+\int_{0}^{t}G\left(s,X_{s}\right)dB_{E_{s}} \label {SDE},
\end{equation}
where the coefficients $F$ and $G$ satisfy some non-globally Lipschitz conditions, $B=\left(B_{t}\right)_{t\geq0}$ is a standard Brownian motion, and $E=\left(E_{t}\right)_{t\geq0}$ is the inverse of a $\alpha$-stable subordinator $D=\left(D_{t}\right)_{t\geq0}$ with $\alpha\in(0,1)$, independent of $B$.
The establishment of such $\alpha$-sensitive framework would provide both theretical foundation and extended applicability for numerical analysis for such SDEs.

%Diffuculties$
A time-changed process refers to a new process $Z_{T_t}$ obtained by transforming the time of a random process $Z_t$ through another nondecreasing process $T_t$, called a random clock or dependent process.
Two processes can be either independent or dependent.
In this way, the classical SDE 
\begin{equation}
  X_{t}=x_{0}+\int_{0}^{t}F\left(s,X_{s}\right)ds+\int_{0}^{t}G\left(s,X_{s}\right)dZ_s \notag
\end{equation}
turns into a TCSDEs, i.e.
\begin{equation}
    X_{t}=x_{0}+\int_{0}^{t}F\left(s,X_{s}\right)dT_{s}+\int_{0}^{t}G\left(s,X_{s}\right)dZ_{T_{s}}.  \notag
\end{equation}
One common subordinated process is $B\circ E=\left(B_{E_t}\right)_{t\geq 0}$, called a time-changed Brownian motion, where $B=\left(B_{t}\right)_{t\geq0}$ is a standard Brownian motion, $E=\left(E_{t}\right)_{t\geq0}$ is the inverse of a subordinator $D=\left(D_{t}\right)_{\geq0}$ with infinite Lévy measure, independent of $B$.
Generalizations of time-changed Brownian motion mainly include time-changed fractional Brownian motion \cite{hahn2011time,hahn2011fokker,meerschaert2009correlated} and SDEs with time-changed semimartingales \cite{2011Kobayashi}.
Furthermore, TCSDEs can also be generalized to be driven by Lévy noise, see \cite{hahn2012sdes,nane2016stochastic,magdziarz2016stochastic}.

When comparing the generalized class of TCSDEs with spatiotemporal coefficients to \eqref{SDE} 
\begin{equation}
  X_{t}=x_{0}+\int_{0}^{t}F\left(X_{s}\right)dE_{s}+\int_{0}^{t}G\left(X_{s}\right)dB_{E_{s}},	\label{1.4}
\end{equation}
the primary challenge stems from the asynchronization between the internal time change $E_{s}$ governing process revolution and the external time variable $s$ in coefficient functions.
If they are synchronized, i.e. 
\begin{equation}
  dX_{t}=F\left(E_{t},X_{t}\right)dE_{t}+G\left(E_{t},X_{t}\right)dB_{E_{t}},X_{0}=x_{0}, \label{SDE2}
\end{equation}
then we can utilize duality principle between TCSDEs and classical Itô SDEs, and much of the work is simplified. 
To be more precise, the duality principle connects TCSDEs \eqref{SDE2} with the dual classical Itô SDE
\begin{equation}
  dY_{t}=F\left(t,Y_{t}\right)dt+G\left(t,Y_{t}\right)dB_{t},Y_{0}=x_{0}. \label{Ito SDE}
\end{equation}
If $Y_{t}$ solves \eqref{Ito SDE},
then $X_{t}:=Y_{E_{t}}$ solves \eqref{SDE2};
while if $X_{t}$ solves \eqref{SDE2},
then $Y_{t}:=X_{D_{t}}$ solves \eqref{Ito SDE}.

Therefore, there are mainly two approaches in analyzing the strong convergence of TCSDEs.
By utilizing the duality principle, 
Jum and Kobayashi \cite{jum2016strong} discussed both strong and weak form of Euler-Maruyama scheme for SDE under standard Lipschitz assumption on the coefficients.
Liu, et al. \cite{liu2020truncated} adopted truncated Euler-Maruyama method to approximate \eqref{SDE} with the Hölder continuity in the temporal variable and the super-linear growth in the state variable. 
Deng \cite{deng2020semi} established backward Euler-Maruyama scheme for \eqref{SDE2} with superlinearly growing coefficients.
For a large class of SDEs driven by a time-changed Brownian motion \eqref{SDE}, since the form of coefficients of SDEs are more complicated, the duality principle cannot be applied. 
Jin and Kobayashi \cite{2019JinKobayashi} investigated Euler-Maruyama scheme for \eqref{SDE}.
Jin and Kobayashi \cite{jin2021strong} studied convergence of Euler-Maruyama scheme and Milstein scheme for TCSDEs involving drifts with random and non-random integrators, under time-varying Lipschitz bound.

Stability bahaviours of solutions of SDEs include the stability of the exact solutions and the numerical stability of the numerical methods. For SDEs driven by time-changed semimartingale, there are already abundant literatures regarding the analytical solutions of the TCSDEs, but there are rarely literature regarding the numerical stability. 
For SDEs driven by time-changed Brownian motion, 
Kobayashi\cite{2011Kobayashi} proved the existence and uniqueness theorem for a strong solution for TCSDEs and showed the stability results on the linear TCSDEs by their explicit solution.
Wu \cite{wu2016stability} investigated the sample-path stability and the $p$th moment asymptotic stability of SDEs time-changed by inverse of a stable subordinator, but without time-space-dependent coefficient.
Deng \cite{deng2020semi} showed that the numerical mean-square polynomial stability of the Backward Euler-Maruyama (BEM) method for TCSDEs without a time-space-dependent coefficient is preserved.
Zhu, Li, Liu \cite{zhu2021exponential} considered the almost sure exponential stability and the $p$th moment exponential stability for the TCSDEs with time-space-dependent coefficient. 
Li, Xu and Ma \cite{li2022global} provided the sufficient condition ensuring the existence of the global attracting sets and proved the exponential stability of stochastic functional differential equations (SFDEs) driven by the time-changed Brownian motion.
He et al. \cite{he2024eta} proved sufficient condtions of the $p$th moment $\eta$-stability and the mean-square $\eta$-stability for a class of SFDEs driven by time-changed Brownian motion. 
Zhang and Yuan \cite{zhang2019razumikhin,zhang2021asymptotic} presented sufficient conditions on the exponential stability for the time-changed SFDEs with Markovian switching, and the asymptotic stability of solutions to the time-changed stochastic delay differential equations with Markovian switching, respectively.
Shen, Zhang, Song and Wu \cite{shen2023class} proved the existence and uniqueness theorem of strong solutions for distribution dependent SDEs driven by time-changed Brownian motion, and obtained sufficients condition of stability in different senses. 
For SDEs driven by time-changed Lévy noise, 
Nane and Ni \cite{nane2016stochastic,nane2017stability} generalized the result in \cite{wu2016stability} to SDEs time-changed by the inverse of a general Lévy subordinator, and studied the exponential stability, the stability in probability, and the asymptotical stability.
Yin \cite{yin2021stability} studied the stability of solutions to nonlinear SDEs driven by the time-changed Lévy noise with impulsive effects, and also the stochastic perturbation for some unstable time-changed differential equations with impulses. 
Xu and Li \cite{xu2023h} investigated $h$-stability of a class of SFDEs driven by time-changed Lévy noise. 

%Contributions
This paper provides affirmative answers to questions (Q1) and (Q2), with its primary contributions summarized as follows:
\begin{itemize}
  \item Generalized numerical scheme and coefficient conditions: while the ST method is a very important implicit method and has been well studied for Itô SDEs and other variants \cite{mao2013strong, zhang2022convergence}, its application to TCSDEs remains unexplored. Building on the motivation from \cite{jin2021strong}, the novelty of our work includes: (i) apply the ST method to TCSDEs whose coefficients are time-space-dependent, without using duality principle, (ii) relax the coefficients assumptions to local Lipschitz condition.
  \item Generalized moment estimates for inverse subordinator $E$: we derive exact formulas using the properties of Mittag-Leffler functions (Theorem 2.1 and 2.2).
  \item Generalized moment inequalities for the exact solution: we establish refined,$\alpha$-sensitive bounds for the exact solution (Theorem 3.1 and Lemma 3.1). 
  \item New dynamic convergence order: our analysis not only determines the convergence rate of the ST scheme but also demonstrates that its strong convergence order varies explicitly with the stability index $\alpha$ of the stable subordinator $D$. 
  This finding, presented in \autoref{thm: convergence}, addresses a gap in traditional analyses that treat the stochastic clock and discretization scheme as independent factors.
  \item Introduction of Mittag-Leffler stability for TCSDEs: we propose a novel stability concept tailored to TCSDEs (Definition 5.1), which accurately describes the slower, algebraic-like decay behavior governed by the Mittag-Leffler function $E_\alpha$, which is a characteristic in subdiffusive systems. Furthermore, we establish a Lyapunov-based criterion for this stability (Theorem 5.1). The theorem explicitly shows that the stability index $\alpha$ of the time change is a fundamental parameter determining the long-term decay rate of the system's moments, providing a unified framework that generalizes classical exponential stability (the special case of $\alpha$=1).
  \item Rigorous analysis of numerical Mittag-Leffler stability: we introduce the concept of numerical Mittag-Leffler stability for discretized schemes (Definition 6.1) and provide a sharp stability criterion for the ST method (\autoref{numerical stability}). Our results reveal that for $\theta \in [1/2, 1]$, the method preserves stability for any step size, while for $\theta \in [0, 1/2)$, stability is maintained under a step size restriction reflective of the system's intrinsic dissipation and noise structure.

\end{itemize}

The paper is organized as follows: 
In Section 2, we survey the properties of the underlying subordinator $D$, derive the exact formulas for moment estimates and exponential estimates of the inverse stable subordinator $E$.
In Section 3, we derive the exact moment bound of the exact solution to \eqref{SDE} under the monotone condition, and the probability of the exact solution will not explode on a finite time interval is estimated. We also establish moment bounds for the difference of two arbitrary exact solutions. All of the above estimates are $\alpha$-sensitive.
In Section 4, we investigate the strong convergence of the ST method. Before proving the main theorem, we define a function $f$ depending on coefficient $F$ and $\theta$, and prove the existence and uniqueness of solution to $f(x)=b$. The mean-square boundedness of the STM scheme \eqref{continuous STM2} is proved. Next we show that the FBEM scheme \eqref{continuous FBEM} and the STM scheme \eqref{continuous STM2} stay close to each other on a compact domain, and estimate the probability that both schemes will not explode on a finite time interval. Finally we utilize the previous results to prove the $\alpha$-scaled strong convergence order of the ST method.
In Section 5, we introduce the notion of the Mittag-Leffler stability for TCSDES and prove the stability theorem. Sections 6 is dedicated to the analysis of mean-square Mittag-Leffler stability for the numerical solutions.
Finally, in Section 6, several examples are studied to illustrate the interest and usefulness of the main results.

\section{Inverse subordinators and associated moments}
Throughout this paper, 
let $(\Omega,\mathcal{F},\mathbb{P})$ be a complete probability space with filtration $\{\mathcal{F}_{t}\}_{t\geq0}$ satisfying the usual conditions (i.e. right continuous and increasing while $\mathcal{F}_{0}$ contains all $\mathbb{P}$-null sets).

Let $D=(D_{t})_{t\geq0}$ be a subordinator starting at 0 with Laplace exponent $\psi$ with killing rate 0, drift 0, and Lévy measure $\mathit{v}$;
i.e. $D$ is a one-dimensional nondecreasing Lévy process with càdlàg paths starting at 0 with Laplace transform
\begin{equation}
  \mathbb{E}[e^{-\xi D_t}] = e^{-t \psi(\xi)}, \notag
\end{equation}
where the Laplace exponent $\psi:(0,\infty)\rightarrow(0,\infty)$ is 
\begin{equation}
  \psi(\xi) = \int_0^\infty (1 - e^{-\xi y}) v(dy), \xi >0, \notag
\end{equation}
satisfying $\int_{0}^{\infty}(y\wedge1)v(dy)<\infty$, where $y\wedge 1$ stands for the minimum of $y$ and $1$.
In this paper, we focus on the infinite Lévy measure case, i.e. $v(0,\infty)=\infty$. This assumption implies that $D$ has strictly increasing paths with infinitely many jumps. 
%E
Let $E=(E_{t})_{t\geq0}$ be the inverse of $D$, i.e.
\begin{equation}
    E_{t}:=\inf\{u>0:D_{u}>t\},t\geq0. \notag
\end{equation}
It is called an inverse subordinator or time change as it "reverses" the time scale of the subordinator.

If the subordinator $D$ is stable with index $\alpha\in(0,1)$, then $\psi(\xi)=\xi ^\alpha$,
and the corresponding time change $E$ is called an inverse $\alpha$-stable subordinator.
If $D$ is a stable subordinator, then $E$ has Mittag-Leffler distributions, see \cite{meerschaert2004limit}.

% subordinated process
The composition $B\circ E=\left(B_{E_t}\right)_{t\geq 0}$ denotes the time-changed process. If $B$ is a standard Brownian motion independent of $D$, we can think of particles represented by the time-changed Brownian motion $B\circ E$ as being trapped and immobile during the constant periods.

Next we derive the following moment property of the inverse subordinator $E$.
\begin{Theorem}\label{Theorem:moment of inverse subordinator E}
  Let $E$ be the inverse of a stable subordinator $D$ with stability index $\alpha\in(0,1)$. Then for any integer $p\geq 1$ and $0\leq s<t$,
  \begin{equation}
    \mathbb{E}\left[E_{t}^{p}\right]=\frac{\Gamma(p+1)}{\Gamma(\alpha p+1)}t^{\alpha p}.
  \end{equation}
  Furthermore,
  \begin{equation}
    \mathbb{E}\left[\left|E_{t}-E_{s}\right|^{p}\right]=\frac{\Gamma(p+1)}{\Gamma(\alpha p+1)}|t-s|^{p\alpha}.
  \end{equation}
\end{Theorem}
\begin{proof}
  Recall that the Laplace transform of $E$ is 
  \begin{equation}
    \mathbb{E}\left[e^{-\xi E_{t}}\right]=E_{\alpha}(-\xi t^{\alpha}), \notag
  \end{equation}
  where $E_{\alpha}(z)=\sum_{k=0}^{\infty}\frac{z^{k}}{\Gamma(\alpha k+1)}$ is a one-parameter Mittag-Leffler function, see \cite{bingham1971limit}.
  The series of the Mittag-Leffler function converges absolutely for $\xi \geq 0$, allowing term-by-term differentiation.
  For $k\geq p$, we have
  \begin{equation}
    \frac{d^{p}}{d\xi^{p}}(-\xi)^{k}=(-1)^{k}\frac{k!}{(k-p)!}\xi{}^{k-p}. \notag
  \end{equation}
  Note that 
  \begin{align}
    \left.(-\xi)^{k-p}\right|_{\xi=0}=\begin{cases}
      1, & k=p\\
      0, & k\neq p
      \end{cases}, \notag
  \end{align}
  then
  \begin{equation}
    \left.\frac{d^{p}}{d\xi^{p}}\mathbb{E}\left[e^{-\xi E_{t}}\right]\right|_{\xi=0}
    =\left.\sum_{k=p}^{\infty}\frac{(-1)^{k}t^{\alpha k}k!}{\Gamma(\alpha k+1)(k-p)!}\xi{}^{k-p}\right|_{\xi=0}
    =\frac{(-1)^{p}p!t^{\alpha p}}{\Gamma(\alpha p+1)}. \notag
  \end{equation}
  Therefore
  \begin{equation}
    \mathbb{E}\left[E_{t}^{p}\right]
    =\left.(-1)^{p}\frac{d^{p}}{d\xi^{p}}\mathbb{E}\left[e^{-\xi E_{t}}\right]\right|_{\xi=0}
    =\frac{\Gamma(p+1)}{\Gamma(\alpha p+1)}t^{\alpha p}. \notag
  \end{equation}
  For the second part of the result, by definition $E_{t}	:=\inf\left\{ u>0:D_{u}>t\right\} $, we have
  \begin{equation}
    E_{t}-E_{s}=\inf\left\{ v>0:D_{E_{s}+v}>t\right\} 
    \overset{d}{=}\inf\left\{ v>0:D_{v}>t-s\right\}. \notag
  \end{equation}
  Since $D_{E_{s}+v}-D_{E_{s}}\overset{d}{=}D_{v}$, we have
  \begin{equation}
    E_{t}-E_{s}\overset{d}{=}E_{t-s}. \notag
  \end{equation}
  Therefore
  \begin{equation}
    \mathbb{E}\left[\left|E_{t}-E_{s}\right|^{p}\right]
    =\mathbb{E}\left[E_{t-s}^{p}\right]
    =\frac{\Gamma(p+1)}{\Gamma(\alpha p+1)}\left|t-s\right|^{\alpha p}. \notag
  \end{equation}
\end{proof}

We now proceed to the exponential moment property of the inverse subordinator $E$. 
Jum and Kobayashi \cite{jum2016strong} proved that an inverse subordinator $E$ with infinite Lévy measure has finite exponential moment; i.e. $\mathbb{E}\left[exp\left(\xi E_{t}\right)\right]<\infty$ for all $\xi,t >0$.
Kobayashi \cite{2019JinKobayashi} established a criterion for the existence of $\mathbb{E}\left[exp\left(\xi E_{t}^{r}\right)\right]$ for all $\xi,t$, and $r>0$ by combining tail probability estimates with properties of regularly varying functions.
Here, we prove the result using the properties of Mittag-Leffler functions. Not only have we proved its finiteness, but we have also provided the explicit expression depending on the stability index $\alpha$.

\begin{Theorem}\label{Theorem: exponential moments of powers of inverse subordinator E}
    Let $E$ be the inverse of a stable subordinator $D$ with stability index $\alpha\in(0,1). $
    Then for $\xi >0$ and $t>0$,
    \begin{equation}
      \mathbb{E}\left[exp\left(\xi E_{t}\right)\right]=\sum_{k=0}^{\infty}\frac{(\xi t^{\alpha})^{k}}{\Gamma(\alpha k+1)}<\infty. \label{2.6}
    \end{equation}
    For $r>0$, 
    \begin{equation}
      \mathbb{E}\left[exp\left(\xi E_{t}^{r}\right)\right]
      =\sum_{k=0}^{\infty}\frac{\xi^{k}}{k!}\frac{\Gamma(rk+1)}{\Gamma(\alpha rk+1)}t^{\alpha rk}. \label{2.8}
    \end{equation}
    1. If $r<1/(1-\alpha)$, then $\mathbb{E}\left[exp\left(\xi E_{t}^{r}\right)\right]<\infty$.

    2. If $r>1/(1-\alpha)$, then $\mathbb{E}\left[exp\left(\xi E_{t}^{r}\right)\right]=\infty$.
\end{Theorem}
\begin{proof}
   The expression \eqref{2.6} comes diretly from the Laplace transform of $E_t$, and is clearly finite.
   To prove the convergence of \eqref{2.8}, we employ ratio test.
   Set $a_{k}:=\frac{\xi^{k}}{k!}\frac{\Gamma(rk+1)}{\Gamma(\alpha rk+1)}t^{\alpha rk}$ as the general term of the series. 
   In order to analyze its asymptotic behaviour,
   apply the Stirling formula for gamma function to the two gamma functions in $a_k$.
   Recall the Stirling formula for gamma function is  $\Gamma(z)\sim\sqrt{\frac{2\pi}{z}}\left(\frac{z}{e}\right)^{z}$ as $z\rightarrow\infty$.
   Therefore,
   \begin{align}
    \frac{\Gamma(rk+1)}{\Gamma(\alpha rk+1)}
    &\sim \frac{\sqrt{2\pi rk}\left(\frac{rk}{e}\right)^{rk}}{\sqrt{2\pi\alpha rk}\left(\frac{\alpha rk}{e}\right)^{\alpha rk}} \notag \\
    &=\frac{1}{\sqrt{\alpha}}\frac{(rk)^{rk}}{(\alpha rk)^{\alpha rk}}\frac{e^{\alpha rk}}{e^{rk}} \notag \\
    &=\frac{1}{\sqrt{\alpha}}\frac{(rk)^{rk}}{(\alpha rk)^{\alpha rk}}e^{-(1-\alpha)rk}\notag \\
    &=\frac{1}{\sqrt{\alpha}}(rk)^{(1-\alpha)rk}\alpha^{-\alpha rk}e^{-(1-\alpha)rk}.\label{ratio}
   \end{align}
   Plugging \eqref{ratio} back into $\alpha_k$ yields
   \begin{align}
    a_{k}
    &\sim \frac{\xi^{k}}{k!}\frac{1}{\sqrt{\alpha}}\left[r^{(1-\alpha)rk}\alpha^{-\alpha rk}e^{-(1-\alpha)rk}t^{\alpha rk}\right]k^{(1-\alpha)rk} \notag \\
    &=\frac{\xi^{k}}{k!}\frac{1}{\sqrt{\alpha}}\left[\left(\frac{r^{(1-\alpha)}}{\alpha^{\alpha}e^{1-\alpha}}\right)^{rk}t^{\alpha rk}\right]k^{(1-\alpha)rk} \notag \\
    &=\frac{\xi^{k}}{k!}\frac{1}{\sqrt{\alpha}}\left[\left(\frac{r}{\alpha^{\alpha/(1-\alpha)}e}\right)^{(1-\alpha)rk}t^{\alpha rk}\right]k^{(1-\alpha)rk}. \notag
   \end{align}
   Then
   \begin{equation}
    \frac{a_{k+1}}{a_{k}}\sim\frac{\xi\left(\frac{r}{\alpha^{\alpha/(1-\alpha)}e}\right)^{(1-\alpha)r}t^{\alpha r}(k+1)^{(1-\alpha)r}}{k+1}.
   \end{equation}
   Taking $k\rightarrow\infty$, we have
   \begin{equation}
    \underset{k\rightarrow\infty}{lim}\frac{a_{k+1}}{a_{k}}\sim\xi\left(\frac{r}{\alpha^{\alpha/(1-\alpha)}e}\right)^{(1-\alpha)r}t^{\alpha r}\underset{k\rightarrow\infty}{lim}k{}^{(1-\alpha)r-1}. \notag
   \end{equation}
   Finally, choosing $r<\frac{1}{1-\alpha}$, the ratio test guarantees that the series converges. 
\end{proof}

\section{$\alpha$-dependent moment estimates}

Throughout the rest of the paper, let $(\Omega,\mathcal{F},\mathbb{P})$ be a complete probability space with filtration $\{\mathcal{F}_{t}\}_{t\geq0}$ satisfying the usual conditions. 
Let $E$ be the inverse of a stable subordinator $D$ with stability index $0<\alpha<1$ and infinite Lévy measure. 
Let $B$ be a standard Brownian motion independent of $D$. 
Consider the TCSDE
\begin{equation}
    X_{t}=x_{0}+\int_{0}^{t}F(s,X_s)dE_{s}+\int_{0}^{t}G(s,X_s)dB_{E_s}, \forall  t\in[0,T],	
\end{equation}
where $X(t_0)=x_0$,
$x_{0}\in\mathbb{R}^{d}$ is a non-random constant,
$T>0$ is a finite time horizon,
$F(t,x):[0,T]\times\mathbb{R}^{d}\rightarrow\mathbb{R}^{d}$ and 
$G(t,x):[0,T]\times\mathbb{R}^{d}\rightarrow\mathbb{R}^{d\times m}$
are measurable functions, and
$F(t,0)=G(t,0)=0$.
In the content going forward, we assume $m=d=1$ for simplicity of discussions and expressions when necessary; an extension to multidimensional case is straightforward. 
% product space
Since the time-changed process consists of two independent processes $B$ and $E$, we can put $B$ and $D$ on a product space with product measure $\mathbb{P}=\mathbb{P}_{B}\times\mathbb{P}_{D}$. 
% expectation
Let $\mathbb{E}_{B},\mathbb{E}_{D}$, and $\mathbb{E}$ denote the expectation under the probability measures $\mathbb{P}_{B},\mathbb{P}_{D}$ and $\mathbb{P}$, respectively. 

We need the following assumptions on the drift and diffusion coefficients.
%varying consant
Throughout this work, we use $C$ to denote generic positive constants that may change from line to line. 

\begin{Assumption}(Local Lipschitz Condition).\label{LL}
  For every $R>0$, there exists a positive constant $C(R)$, depending only on $R$, such that for all $x,y\in\mathbb{R}^{d}$ and $t\geq 0$,
  \begin{equation}
    \left|F(t,x)-F(t,y)\right|+|G(t,x)-G(t,y)|\leq C(R)|x-y|, \notag
  \end{equation}
  whenever $|x|\vee|y|\leq R$.
\end{Assumption}

The following polynomial growth condition plays a crucial role following non-global Lipschitz condition, as it ensures that the higher-order moments of numerical solutions do not explode.
\begin{Assumption}\label{polynomial growth}
  There exists a pair of constants $h\geq 1$ and $C(h)>0$ such that, for all $t\in[0,T]$ and $x\in\mathbb{R}^{d}$,
  \begin{equation}
    \left|F(t,x)\right|\vee|G(t,x)|\leq C(h)\left(1+|x|^{h}\right). \notag
  \end{equation}
\end{Assumption}

\begin{Assumption}\label{monotone}
    There exists a constant $K_1>0$ such that, for all $t\in[0,T]$ and $x\in\mathbb{R}^{d}$,
    \begin{equation}
      \left\langle x,F(t,x)\right\rangle +\frac{2h-1}{2}|G(t,x)|^{2}\leq K_1(1+|x|^{2}).  \notag
    \end{equation}
\end{Assumption}

We also impose the following assumption to deal with the time-dependent coefficients.
\begin{Assumption}\label{continuity}
  There exist constants $K_2, K_3>0,\eta_F\in(0,1]$ and $\eta_G\in(0,1]$ such that, for all $s,t\in[0,T]$ and $x,y\in\mathbb{R}^{d}$,
  \begin{equation}
    |F(s,x)-F(t,x)|\leq K_2(|1+|x|)|s-t|^{\eta_{F}}, \notag
  \end{equation}
  and
  \begin{equation}
    |G(s,x)-G(t,x)|\leq K_3(|1+|x|)|s-t|^{\eta_{G}}. \notag
  \end{equation}
\end{Assumption}

A semimartingale is composed of a local martingale and a finite variation process. 
Since $E$ is a finite variation process and $B$ is independent of $E$, $B\circ E$ is a continuous local martingale. Thus $E$ and $B\circ E$ are both semimartingales. The proof of the existence and uniqueness theorem of the strong solution to the TCSDE \eqref{SDE} under the local Lipschitz assumption is similar to the proof of the existence and uniquess theorem of the strong solution to SDEs driven by semimartingale, see \cite{Protter2004} and \cite{mao1991semimartingales}. 

Next we derive the upper bound for the exact solution of \eqref{SDE} and the probability that the solution stays within a compact domain for finite time $T>0$.
We show that both estimates quantitatively depend on the value of the stability index $\alpha$.

For each positive number $R$, we define the stopping time
\begin{equation}
  \kappa_{R}:=\inf\{t\geq0:|X_{t}|>R\}. \label{kappa}
\end{equation}
\begin{Theorem}\label{thm: exact pth moment, probability}
  Let Assumptions 3.1-3.3 hold.
  Then, for $h\geq 1$ and any $0\leq t<T$,
  \begin{equation}
    \mathbb{E}\left[\left|X_{t}\right|^{2h}\right]
    \leq 2^{h-1}\sum_{k=0}^{\infty}\frac{(2hK_{1}t^{\alpha})^{k}}{\Gamma(\alpha k+1)}\left(1+\mathbb{E}\left|x_{0}\right|^{2h}\right)
    <\infty, \notag
  \end{equation}
  and
  \begin{equation}
    \mathbb{P}\left(\kappa_{R}\leq T\right)
    \leq \frac{2^{h-1}}{R^{2h}}\sum_{k=0}^{\infty}\frac{(2hK_{1}T^{\alpha})^{k}}{\Gamma(\alpha k+1)}\left(1+\mathbb{E}\left|x_{0}\right|^{2h}\right). \notag
  \end{equation}
  Furthermore, let $Y_{t}^{(h)}:=1+\underset{0\leq s\leq t}{\sup}|X_{s}|^{h}$.
  Then $\mathbb{E}[Y_{t}^{(h)}]<\infty.$
\end{Theorem}

\begin{proof}
  Similar to the proof of Theorem 4.1 in Section 2.4 of \cite{maobook}, using the time-changed Itô's formula in \cite{2011Kobayashi} and \autoref{monotone}, we have
  \begin{align}
    \left(1+\left|X_{t}\right|^{2}\right)^{h}
    &\leq 2^{h-1}\left(1+\left|x_{0}\right|^{2h}\right)+2h\int_{0}^{t}\left[1+\left|X_{s}\right|^{2}\right]^{h-1}\left(\left\langle X_{s},F\left(s,X_{s}\right)\right\rangle +\frac{2h-1}{2}\left|G\left(s,X_{s}\right)\right|^{2}\right)dE_{s} \notag \\
    &\quad +2h\int_{0}^{t}\left(1+\left|X_{s}\right|^{2}\right)^{h-1}\left\langle X_{s},G\left(s,X_{s}\right)\right\rangle dB_{E_{s}} \notag \\
    &\leq 2^{h-1}\left(1+\left|x_{0}\right|^{2h}\right)+2hK_{1}\int_{0}^{t}\left(1+\left|X_{s}\right|^{2}\right)^{h}dE_{s} \notag \\
    &\quad +2h\int_{0}^{t}\left(1+\left|X_{s}\right|^{2}\right)^{h-1}\left\langle X_{s},G\left(s,X_{s}\right)\right\rangle dB_{E_{s}}.
  \end{align}
  Taking $\mathbb{E}_{B}$ on both sides, since the independence of $B$ and $E$ ensures that the time change $E$ does not affect the martingale property of $B$ under $\mathbb{P}_{B}$, we have
  \begin{align}
    \mathbb{E}_{B}\left[\left(1+\left|X_{t\wedge\kappa_{R}}\right|^{2}\right)^{h}\right]
    &\leq 2^{h-1}\left(1+\mathbb{E}_{B}\left[\left|x_{0}\right|^{2h}\right]\right)+2hK_{1}\mathbb{E}_{B}\left[\int_{0}^{t\wedge\kappa_{R}}\left(1+\left|X_{s}\right|^{2}\right)^{h}dE_{s}\right] \notag  \\
    &\leq 2^{h-1}\left(1+\mathbb{E}_{B}\left[\left|x_{0}\right|^{2h}\right]\right)+2hK_{1}\mathbb{E}_{B}\left[\int_{0}^{t}\left(1+\left|X\left(s\wedge\kappa_{R}\right)\right|^{2}\right)^{h}dE_{s}\right]. \notag
  \end{align}
  Applying the Grönwall-type inequality (Lemma 6.3 in Chapter IX.6a of \cite{Jacod2003}), taking $\mathbb{E}_{D}$ on both sides, noting that $\mathbb{E}_{D}\left[\mathbb{E}_{B}\left[\left|X_{t\wedge\kappa_{R}}\right|^{2}\right]\right]=\mathbb{E}\left[\left|X_{t\wedge\kappa_{R}}\right|^{2}\right]$, and letting $R\rightarrow\infty$ yields
  \begin{equation}
    \mathbb{E}\left[\left(1+\left|X_{t}\right|^{2}\right)^{h}\right]\leq2^{h-1}\left(1+\mathbb{E}\left|x_{0}\right|^{2h}\right)\mathbb{E}\left[e^{2hK_{1}E_{t}}\right]. \notag
  \end{equation}
  Then
  \begin{equation}
    \mathbb{E}\left[\left|X_{t}\right|^{2h}\right]
    \leq \mathbb{E}\left[\left(1+\left|X_{t}\right|^{2}\right)^{h}\right]
    \leq 2^{h-1}\sum_{k=0}^{\infty}\frac{(2hK_{1}t^{\alpha})^{k}}{\Gamma(\alpha k+1)}\left(1+\mathbb{E}\left|x_{0}\right|^{2h}\right)
    <\infty. \notag
  \end{equation}
  Therefore
  \begin{equation}
    \mathbb{P}\left(\kappa_{R}\leq T\right)
    \leq \mathbb{P}\left(\left|X_{T\wedge\kappa_{R}}\right|\geq R,\kappa_{R}\leq T\right)
    \leq \frac{2^{h-1}}{R^{2h}}\sum_{k=0}^{\infty}\frac{(2hK_{1}T^{\alpha})^{k}}{\Gamma(\alpha k+1)}\left(1+\mathbb{E}\left|x_{0}\right|^{2h}\right). \notag
  \end{equation}

\end{proof}

We also need the following lemma for the proof of the main convergence theorem.
\begin{Lemma}\label{exact any two difference}
  Let \autoref{polynomial growth} hold. 
  Then for $h\geq 1$ and any $0\leq s\leq t$, 
  \begin{equation}
    \mathbb{E}_{B}\left[\left|X_{t}-X_{s}\right|\right]
    \leq C(h)\mathbb{E}_{B}\left[Y_{t}^{(2h)}\right]^{1/2}\left\{ \left|t-s\right|^{\alpha}+\left|t-s\right|^{\alpha/2}\right\}, \notag
  \end{equation}
  where $Y_{t}^{(2h)}$ is defined in \autoref{thm: exact pth moment, probability}.
\end{Lemma}
\begin{proof}
  Applying the Cauchy-Schwarz inequality for $\int_{s}^{t}G(r,X_{r})dB_{E_{r}}$ and \autoref{polynomial growth}, we have
  \begin{align}
    \mathbb{E}_{B}\left[\left|X_{t}-X_{s}\right|\right]
    &\leq \mathbb{E}_{B}\left|\int_{s}^{t}F(r,X_{r})dE_{r}\right|+\mathbb{E}_{B}\left|\int_{s}^{t}G(r,X_{r})dB_{E_{r}}\right| \notag \\
    &\leq \mathbb{E}_{B}\left[\int_{s}^{t}\left|F(r,X_{r})\right|dE_{r}\right]+\mathbb{E}_{B}\left[\left|\int_{s}^{t}G(r,X_{r})dB_{E_{r}}\right|^{2}\right]^{1/2} \notag \\
    &\leq \mathbb{E}_{B}\left[\int_{s}^{t}C(h)\left(1+\left|X_{r}\right|^{h}\right)dE_{r}\right]+\mathbb{E}_{B}\left[\left|\int_{s}^{t}C(h)\left(1+\left|X_{r}\right|^{h}\right)dB_{E_{r}}\right|^{2}\right]^{1/2} \notag \\
    &\leq C(h)\left(E_{t}-E_{s}\right)\mathbb{E}_{B}\left[Y_{t}^{(h)}\right]+C(h)\mathbb{E}_{B}\left[\int_{s}^{t}\left(Y_{t}^{(h)}\right)^{2}dE_{r}\right]^{1/2} \notag \\
    &\leq C(h)\left(E_{t}-E_{s}\right)\mathbb{E}_{B}\left[Y_{t}^{(h)}\right]+\sqrt{2}C(h)(E_{t}-E_{s})^{1/2}\mathbb{E}_{B}\left[Y_{t}^{(2h)}\right]^{1/2}. \label{difference} 
  \end{align}
  According to \autoref{Theorem:moment of inverse subordinator E} we have
  \begin{equation}
    \mathbb{E}\left[\left|E_{t}-E_{s}\right|\right]
    =\frac{\Gamma(2)}{\Gamma(\alpha+1)}|t-s|^{\alpha}
    =\frac{1}{\Gamma(\alpha+1)}|t-s|^{\alpha}. \label{E1}
  \end{equation}
  And by the Jensen's inequality we have
  \begin{equation}
    \mathbb{E}\left[\left|E_{t}-E_{s}\right|^{1/2}\right]
    \leq \left(\mathbb{E}\left|E_{t}-E_{s}\right|\right)^{1/2}
    =\frac{1}{\sqrt{\Gamma(\alpha+1)}}|t-s|^{\alpha/2}. \label{E2}
  \end{equation}
  Taking $\mathbb{E}_{D}$ on both sides of \eqref{difference} and using the identity \eqref{E1} and \eqref{E2}, we have
  \begin{align}
    \mathbb{E}\left[\left|X_{t}-X_{s}\right|\right]
    &\leq C(h)\frac{1}{\Gamma(\alpha+1)}\left|t-s\right|^{\alpha}\mathbb{E}\left[Y_{t}^{(h)}\right]+C(h)\frac{1}{\sqrt{\Gamma(\alpha+1)}}|t-s|^{\alpha/2}\mathbb{E}\left[Y_{t}^{(2h)}\right]^{1/2} \notag \\
    &\leq C(h)\mathbb{E}\left[Y_{t}^{(2h)}\right]^{1/2}\left\{ \left|t-s\right|^{\alpha}+\left|t-s\right|^{\alpha/2}\right\} .  \notag
  \end{align}
\end{proof}

\section{$\alpha$-scaled strong convergence of the ST scheme} 

In order to approximate the solution $X_t$ of SDE \eqref{SDE} we first approximate the inverse subordinator $E$ and then the time-changed Brownian motion $B\circ E$ on the interval $[0,T]$.
%approximation of E
Approximation of $E$ relies upon the subordinator $D$.
Fix an equidistant step size $\delta\in(0,1)$ and a time horizon $T>0$.
We first simulate a sample path of the subordinator $D$, by setting $D_{0}=0$ and $D_{i\delta}:=D_{(i-1)\delta}+Z_{i},i=1,2,3,...,$ with an i.i.d. sequence $\left\{Z_{i}\right\}_{i\in\mathbb{N}}$ distributed as $Z_{i}\overset{d}{=}D_{\delta}.$
One may generate the random variables $Z_{i}$ via the Chambers-Mallows-Stuck (CMS) algorithm \cite{chambers1976method}.
We stop this procedure upon finding the integer $N$ satisfying $T\in[D_{N\delta},D_{(N+1)\delta}).$ 
Next, let 
\begin{equation}
  \tilde{E}_{t}:=\left(\min\left\{ n\in\mathbb{N};D_{n\delta}>t\right\} -1\right)\delta,\forall t\in[0,T]. \label{Et}
\end{equation}
The sample paths of $\tilde{E}=(\tilde{E}_{t})_{t\geq0}$ are nondecreasing step functions with constant jump size $\delta$.
$Z_{i}=D_{i\delta}-D_{(i-1)\delta}$ represents the $i$th waiting time.
We can see that for $n=0,1,2,...,N$, we have $\tilde{E}_{t}=n\delta$ whenver $t\in[D_{n\delta},D_{(n+1)\delta})$.
It is also proved that the process $\tilde{E}_{t}$ efficiently approximates $E$, i.e. 
\begin{equation}
  E_{t}-\delta\leq\tilde{E}_{t}\leq E_{t},\forall t\in[0,T]. \label{approximates}
\end{equation}
For proofs, see \cite{jum2016strong} and \cite{magdziarz2009stochastic}.
%approximate B
Since we assume that $B$ and $D$ are independent, we can simply approximate $B\circ E$ over the time steps $\{0,\delta,2\delta,...,N\delta\}$.

%approximate ST
To approximate the solution we consider the ST scheme. Let $\tau_{n}=D_{n\delta},n=0,1,2,...,N$. Given any random step size $\tau_{n+1}-\tau_n,$ define a discrete-time process $(\tilde{X}_{\tau_{n}})$ by setting
\begin{align}
    & \tilde{X}_{0}=x_{0}, \\ 
    & \tilde{X}_{\tau_{n+1}}	=\tilde{X}_{\tau_{n}}+(1-\theta)F\left(\tau_{n},\tilde{X}_{\tau_{n}}\right)\delta+\theta F\left(\tau_{n+1},\tilde{X}_{\tau_{n+1}}\right)\delta+G\left(\tau_{n},\tilde{X}_{\tau_{n}}\right)\left(B_{(n+1)\delta}-B_{n\delta}\right), 	\label{discrete STM}
\end{align} 
for $n=0,1,2,...,N-1$,
where $\theta \in[0,1]$. 
The parameter $\theta$ controls the implicitness of the scheme. 
It is well-known that the ST method includes the Euler-Maruyama (EM) method ($\theta=0$), the trapezoidal method ($\theta=1/2$) and the Backward Euler-Maruyama (BEM) method ($\theta=1$).

Now, we define a continuous-time process $\tilde{X}=\left(\tilde{X}_{t}\right)_{t\in[0,T]}$ by continuous interpolation; i.e. whenver $s\in[\tau_{n},\tau_{n+1}),$
\begin{equation}
  \tilde{X}_{s}:=\tilde{X}_{\tau_{n}}+(1-\theta)\int_{\tau_{n}}^{s}F\left(\tau_{n},\tilde{X}_{\tau_{n}}\right)dE_{r}+\theta\int_{\tau_{n}}^{s}F\left(\tau_{n+1},\tilde{X}_{\tau_{n+1}}\right)dE_{r}+\int_{\tau_{n}}^{s}G\left(\tau_{n},\tilde{X}_{\tau_{n}}\right)dB_{E_{r}}. \label{continuous STM1}
\end{equation}
Let
\begin{equation}
  n_{t}=\max\left\{ n\in\mathbb{N}\cup\{0\};\tau_{n}\leq t\right\} ,t\geq0. \notag
\end{equation}
Clearly, $\tau_{n_{t}}\leq t\leq\tau_{n_{t}+1}$ for any $t>0$.
Note by \eqref{continuous STM1} and the identity $\tilde{X}_{s}-x_{0}=\sum_{i=0}^{n_{s}-1}\left(\tilde{X}_{\tau_{i+1}}-\tilde{X}_{\tau_{i}}\right)+\left(\tilde{X}_{s}-\tilde{X}_{\tau_{n_{s}}}\right)$, we have
\begin{equation}
  \tilde{X}_{s}-x_{0}=\sum_{i=0}^{n_{s}-1}\left[F\left(\tau_{i},\tilde{X}_{\tau_{i}}\right)\delta+G\left(\tau_{i},\tilde{X}_{\tau_{i}}\right)B_{(i+1)\delta}-B_{i\delta}\right]+\left(\tilde{X}_{s}-\tilde{X}_{\tau_{n_{s}}}\right).\notag
\end{equation}
Furthermore, since $\tau_{i}=\tau_{n_{r}}$ for any $r\in[\tau_{n_{r}},\tau_{n_{r}+1}), $ \eqref{continuous STM1} can be rewritten as
\begin{equation}
  \tilde{X}_{s}:=x_0+(1-\theta)\int_{0}^{s}F\left(\tau_{n_{r}},\tilde{X}_{\tau_{n_{r}}}\right)dE_{r}+\theta\int_{0}^{s}F\left(\tau_{n_{r}+1},\tilde{X}_{\tau_{n_{r}+1}}\right)dE_{r}+\int_{0}^{s}G\left(\tau_{n_{r}},\tilde{X}_{\tau_{n_{r}}}\right)dB_{E_{r}}.    \label{continuous STM2}
\end{equation}

In order to implement numerical scheme \eqref{discrete STM} more easily, we define a function $f:\mathbb{R}^{d}\rightarrow\mathbb{R}^{d}$ by
\begin{equation}
  f(x)=x-\theta F(t,x)\delta,\forall x\in\mathbb{R}^{d}. \label{f}
\end{equation}

We further impose the following assumption.
\begin{Assumption}\label{one-sided Lipschitz}
  There exists a constant $K_4$ such that, for all $t\in[0,T]$ and $x,y \in\mathbb{R}^{d}$,
  \begin{equation}
    \left\langle x-y,F(t,x)-F(t,y)\right\rangle 	\leq K_{4}|x-y|^{2}. \notag
  \end{equation}
\end{Assumption}

The lemma below gives the existence and uniqueness of the solution to the equation $f(x)=b$. It is a direct application of the property of the monotone operator (See Theorem 26.A in \cite{zeidler2013nonlinear} and \cite{mao2013strongBEM}).
\begin{Lemma}\label{f(x)=b existence and uniqueness}
  Let Assumptions 3.2 and 4.1 hold.
  Let $f$ be defined by \eqref{f}. 
  Then,  for any $b\in\mathbb{R}^{d}$ and $\delta<\delta^{*}=\min\left\{1,\frac{1}{K_1\theta}, \frac{1}{K_4\theta}\right\}$, there exists a unique $x\in\mathbb{R}^{d}$ such that
  \begin{equation}
    f(x)=b. \notag
  \end{equation}
\end{Lemma}
\begin{proof}
  For $\delta<\min\left\{1,\frac{1}{ K_{4}\theta}\right\}$, applying \autoref{one-sided Lipschitz} to \eqref{f}, we have
  \begin{equation}
    \left\langle x-y,f(x)-f(y)\right\rangle 
    \geq \left(1-\theta\delta K_{4}\right)|x-y|^{2}
    >0. \notag
  \end{equation}
  Thus $f(x)$ is monotone. 
  Now by \autoref{monotone},
  \begin{equation}
    \left\langle x,f(x)\right\rangle =\left\langle x,x-\theta F(t,x)\delta\right\rangle 
    \geq|x|^{2}-\theta K_1(1+|x|^{2})\delta. \notag
  \end{equation}
  If  $\delta<\frac{1}{K_{1}\theta }$, then 
  \begin{equation}
    \underset{|x|\rightarrow\infty}{lim}\frac{\left\langle x,f(x)\right\rangle }{|x|}=\infty \notag
  \end{equation}
  holds on $\mathbb{R}^{d}$.
  Since $f(x)$ is monotone and coercive,  $f(x)$ is surjective. 
  The lemma is then proved.
\end{proof}
We can see that an inverse operator $f^{-1}$ exists and the solution to \eqref{discrete STM} can be represented in the following form
\begin{equation}
  \tilde{X}_{\tau_{n+1}}	=f^{-1}\left(\tilde{X}_{\tau_{n}}+(1-\theta)F\left(\tau_{n},\tilde{X}_{\tau_{n}}\right)\delta+G\left(\tau_{n},\tilde{X}_{\tau_{n}}\right)\left(B_{(n+1)\delta}-B_{n\delta}\right)\right). \notag
\end{equation} 

In this section we list moment properties of the numerical solution of the ST scheme.

%%f 2nd bound
\begin{Lemma}\label{solution, 2nd moment bound, using f}
  Let Assumptions 3.3 and 4.1 hold.
  Let $f$ be defined by \eqref{f}. 
  Then, for any $x\in\mathbb{R}^{d}$ and $\delta<\delta^{*}=\min\left\{1,\frac{1}{2K_1\theta}, \frac{1}{K_4\theta}\right\}$, 
  \begin{equation}
    |x|^{2} \leq \frac{|f(x)|^{2}+2K_{1}\theta\delta}{(1-2K_{1}\theta\delta)}. \notag
  \end{equation}
\end{Lemma}
\begin{proof}
  According to \autoref{f(x)=b existence and uniqueness}, applying \autoref{monotone} to \eqref{f}, we have
  \begin{align}
    |f(x)|^{2}
    &=\left\langle f(x),f(x)\right\rangle  \notag \\
    &=\left\langle x-\theta F(t,x)\delta,x-\theta F(t,x)\delta\right\rangle \notag \\
    & \geq |x|^{2}-2K_{1}\theta\delta-2K_{1}\theta|x|^{2}\delta+\theta^{2}|F(t,x)|^{2}\delta^{2}.  \notag
  \end{align}
  Ignoring the higher order term $\theta^{2}|F(t,x)|^{2}\delta^{2}$ gives 
  \begin{equation}
    |f(x)|^{2}
    \geq (1-2K_{1}\theta\delta)|x|^{2}-2K_{1}\theta\delta. \notag
  \end{equation}
  Choosing $\delta$ sufficiently small for $\delta<\frac{1}{2K_1\theta}$, we obtain the desired result.
\end{proof}

The following discrete Grönwall inequality also plays an important role in dealing with discrete solutions. 
\begin{Lemma}(The Discrete Grönwall Inequality).\label{lemma: Discrete Gronwall Inequality} Let $N$ be a positive integer. 
  Let $u_n,w_n$ be nonnegative number for $n=1,2,...,N.$ 
  If 
  \begin{equation} u_{n}\leq u_0 +\sum_{k=0}^{n-1}u_{k}w_{k},\forall n=1,2,...N, \notag 
  \end{equation} 
  then 
  \begin{equation} u_{n}\leq u_0 exp\left(\sum_{k=0}^{n-1}w_{k}\right),\forall n=1,2,...N. \notag 
  \end{equation} 
\end{Lemma} 
The proof can be found in \cite{mao1991semimartingales}.

The following Burkholder-Davis-Gundy inequality is a remarkable result relating the maximum of a local martingale with its quadratic variation, see Proposition 3.26 and Theorem 3.28 of \cite{Karatzas1998}. It gives an upper bound for the expectation of the maximum of the process in terms of the expected cumulative size of the fluctuations.
\begin{Lemma}\label{BDG}
    For any $p>0$, there exists postive constants $b_{p}$ such that
    \begin{equation}
        \mathbb{E}\left[\underset{0\leq t\leq S}{\sup}|M_{t}|^{p}\right]	\leq b_{p}\mathbb{E}\left[[M,M]_{S}^{p/2}\right], \notag
    \end{equation}
where 
$S$ is any stopping time, and 
$M$ is any continuous local martingale with quadratic variation $[M,M]$.
\end{Lemma}

\subsection{The ST scheme and associated moments}
We define the stopping time 
\begin{equation}
  \vartheta_{R}:=\inf\left\{ n\geq0:\left|\tilde{X}_{\tau_{n}}\right|>R\right\}.\notag
\end{equation}

\begin{Theorem}\label{thm: discrete STM, 2nd moment bound, sup}
  Let Assumptions 3.1-3.3 and 4.1 hold. 
  Then, for any $T>0$ and $\delta<\delta^{*}=\min\left\{ 1,\frac{1}{2K_{1}\theta},\frac{1}{K_{4}\theta}\right\}$, there exists a constant $C(T)>0$ such that 
  \begin{equation}
    \underset{\delta<\delta^{*}}{\sup}\underset{0\leq\tau_{n}\leq T}{\sup}\mathbb{E}_{B}\left[\left|\tilde{X}_{\tau_{n}}\right|^{2}\right]<C(T). \notag
  \end{equation}
\end{Theorem}
\begin{proof}
  By definition \eqref{f} of function $f$, we have
  \begin{equation}
    f\left(\tilde{X}_{\tau_{n+1}}\right)
    =f\left(\tilde{X}_{\tau_{n}}\right)+F\left(\tau_{n},\tilde{X}_{\tau_{n}}\right)\delta+G\left(\tau_{n},\tilde{X}_{\tau_{n}}\right)\left(B_{(n+1)\delta}-B_{n\delta}\right). \notag
  \end{equation}
  Squaring boths sides we have  
  \begin{align}
    \left|f\left(\tilde{X}_{\tau_{n+1}}\right)\right|^{2}
    &=\left|f\left(\tilde{X}_{\tau_{n}}\right)\right|^{2}+(1-2\theta)\left|F\left(\tau_{n},\tilde{X}_{\tau_{n}}\right)\right|^{2}\delta^{2} \notag \\
    &\quad +\left(2\left\langle \tilde{X}_{\tau_{n}},F\left(\tau_{n},\tilde{X}_{\tau_{n}}\right)\right\rangle +\left|G\left(\tau_{n},\tilde{X}_{\tau_{n}}\right)\right|^{2}\right)\delta+\delta M_{\tau_{n+1}}
  \end{align}
  where
  \begin{align}
    \delta M_{\tau_{n+1}}
    &=\left|G\left(\tau_{n},\tilde{X}_{\tau_{n}}\right)\left(B_{(n+1)\delta}-B_{n\delta}\right)\right|^{2}-\left|G\left(\tau_{n},\tilde{X}_{\tau_{n}}\right)\right|^{2}\delta \notag \\
    &\quad +2\left\langle f\left(\tilde{X}_{\tau_{n}}\right),G\left(\tau_{n},\tilde{X}_{\tau_{n}}\right)\left(B_{(n+1)\delta}-B_{n\delta}\right)\right\rangle  \notag \\
    &\quad +2\left\langle F\left(\tau_{n},\tilde{X}_{\tau_{n}}\right)\delta,G\left(\tau_{n},\tilde{X}_{\tau_{n}}\right)\left(B_{(n+1)\delta}-B_{n\delta}\right)\right\rangle.
  \end{align}
  By \autoref{LL}, \autoref{monotone} and the fact that $\theta\in[1/2,1]$, we have
  \begin{align}
    \left|f\left(\tilde{X}_{\tau_{n+1}}\right)\right|^{2}
    &\leq \left|f\left(\tilde{X}_{\tau_{n}}\right)\right|^{2}+(1-2\theta)\left|F\left(\tau_{n},\tilde{X}_{\tau_{n}}\right)\right|^{2}\delta^{2}+2K_{1}\left(1+\left|\tilde{X}_{\tau_{n}}\right|^{2}\right)\delta+\delta M_{\tau_{n+1}} \notag \\
    &\leq \left|f\left(\tilde{X}_{\tau_{n}}\right)\right|^{2}+2K_{1}\left(1+\left|\tilde{X}_{\tau_{n}}\right|^{2}\right)\delta+\delta M_{\tau_{n+1}} \notag \\
    &\leq \left|f\left(\tilde{X}_{\tau_{n}}\right)\right|^{2}+2K_{1}\delta+2K_{1}\left|\tilde{X}_{\tau_{n}}\right|^{2}\delta+\delta M_{\tau_{n+1}} \label{3.12}
  \end{align}
  Let $N$ be any nonnegative integer such that $N\delta\leq T$. Summing up both sides of the inequality \eqref{3.12} from $n=0$ to $N\wedge\vartheta_{R}$, we have
  \begin{align}
    \left|f\left(\tilde{X}_{\tau_{N\wedge\vartheta_{R}+1}}\right)\right|^{2}
    &\leq \left|f\left(x_0\right)\right|^{2}+2K_{1}T+2K_{1}\sum_{n=0}^{N\wedge\vartheta_{R}}\left|\tilde{X}_{\tau_{n}}\right|^{2}\delta+\sum_{n=0}^{N\wedge\vartheta_{R}}\delta M_{\tau_{n+1}} \notag \\
    &\leq \left|f\left(x_0\right)\right|^{2}+2K_{1}T+2K_{1}\sum_{n=0}^{N}\left|\tilde{X}_{\tau_{n\wedge\vartheta_{R}}}\right|^{2}\delta+\sum_{n=0}^{N}\delta M_{\tau_{n+1}}\boldsymbol{1}_{[0,\vartheta_{R}]}(n).
  \end{align}
  Applying \autoref{thm: exact pth moment, probability}, \autoref{polynomial growth} and noting that $\tilde{X}_{\tau_{n}}$ and $\boldsymbol{1}_{[0,\theta_{R}]}(n)$ are $\mathcal{F}_{\tau_{n}}$-measurable while $\left(B_{(n+1)\delta}-B_{n\delta}\right)$ is independent of $\mathcal{F}_{\tau_{n}}$, taking $\mathbb{E}_{B}$ on both sides gives
  \begin{equation}
    \mathbb{E}_{B}\left[\left|f\left(\tilde{X}_{\tau_{N\wedge\vartheta_{R}+1}}\right)\right|^{2}\right]
    \leq \left|f\left(\tilde{X}_{\tau_{0}}\right)\right|^{2}+2K_{1}T+2K_{1}\mathbb{E}_{B}\left[\sum_{n=0}^{N}\left|\tilde{X}_{\tau_{n\wedge\vartheta_{R}}}\right|^{2}\delta\right]. \notag
  \end{equation}
  For $\delta<\delta^*$, by \autoref{solution, 2nd moment bound, using f} we have
  \begin{align}
    \mathbb{E}_{B}\left[\left|f\left(\tilde{X}_{\tau_{N\wedge\vartheta_{R}+1}}\right)\right|^{2}\right]
    &\leq \left|f\left(\tilde{X}_{\tau_{0}}\right)\right|^{2}+2K_{1}T+2K_{1}\mathbb{E}_{B}\left[\sum_{n=0}^{N}(1-2K_{1}\theta\delta)^{-1}\left(\left|f\left(\tilde{X}_{\tau_{n}\wedge\vartheta_{R}}\right)\right|^{2}+2K_{1}\theta\delta\right)\delta\right] \notag \\
    &\leq \left|f\left(\tilde{X}_{\tau_{0}}\right)\right|^{2}+2K_{1}T+\frac{4K_{1}^{2}\theta\delta^{2}}{1-2K_{1}\theta\delta}+\frac{2K_{1}}{1-2K_{1}\theta\delta}\mathbb{E}_{B}\left[\sum_{n=0}^{N}\left|f\left(\tilde{X}_{\tau_{n\wedge\vartheta_{R}}}\right)\right|^{2}\delta\right]. \notag 
  \end{align}
  From this and the discrete Grönwall inequality, we have
  \begin{align}
    \mathbb{E}_{B}\left[\left|f\left(\tilde{X}_{\tau_{N\wedge\vartheta_{R}+1}}\right)\right|^{2}\right]
    &\leq \left[\left|f\left(\tilde{X}_{\tau_{0}}\right)\right|^{2}+\left(2K_{1}+\frac{4K_{1}^{2}\theta\delta}{1-2K_{1}\theta\delta}\right)(T+\delta)\right]exp\left(\frac{2K_{1}(T+\delta)}{1-2K_{1}\theta\delta}\right). \notag
  \end{align}
  Thus, letting $R\rightarrow\infty$ and applying the Fatou's lemma, we get
  \begin{equation}
    \mathbb{E}_{B}\left[\left|f\left(\tilde{X}_{\tau_{N+1}}\right)\right|^{2}\right]
    \leq \left[\left|f\left(\tilde{X}_{\tau_{0}}\right)\right|^{2}+\left(2K_{1}+\frac{4K_{1}^{2}\theta\delta}{1-2K_{1}\theta\delta}\right)(T+\delta)\right]exp\left(\frac{2K_{1}(T+\delta)}{1-2K_{1}\theta\delta}\right). \notag
  \end{equation}
  Applying \autoref{solution, 2nd moment bound, using f} again gives
  \begin{equation}
    \underset{\delta<\delta^{*}}{\sup}\underset{0\leq\tau_{n}\leq T}{\sup}\mathbb{E}_{B}\left[\left|\tilde{X}_{\tau_{n}}\right|^{2}\right]<C(T). \notag
  \end{equation}
\end{proof}

\subsection{FBEM scheme and associated estimates} 
In order to prove the main theorem, we first introduce the discrete Forward-Backward Euler-Maruyama (FBEM) scheme
\begin{equation}
  \hat{X}_{\tau_{n+1}}=\hat{X}_{\tau_{n}}+F\left(\tau_{n},\tilde{X}_{\tau_{n}}\right)\delta+G\left(\tau_{n},\tilde{X}_{\tau_{n}}\right)\left(B_{(n+1)\delta}-B_{n\delta}\right), \label{discrete FBEM}
\end{equation}
where $\hat{X}_{\tau_{n}}=x_{0}$.
Then the continuous FBEM scheme is 
\begin{equation}
  \hat{X}_{s}	=\hat{x}_{s}+\int_{0}^{s}F\left(\tau_{n_{r}},\tilde{X}_{\tau_{n_{r}}}\right)dE_{r}+\int_{0}^{s}G\left(\tau_{n_{r}},\tilde{X}_{\tau_{n_{r}}}\right)dB_{E_{r}}. \label{continuous FBEM}
\end{equation}

Define the stopping time $v_{R}:=\inf\{t\geq0:|\hat{X}_{t}|>R\ or\ |\tilde{X}_{t}|>R\}$, we can get the following lemma.
\begin{Lemma}\label{cor: continuous FBEM-STM, pth moment bound, compact}
  Let Assumptions 3.1-3.4 and 4.1 hold, and $\theta\in[1/2,1]$.
  Then, there exists a constant $C(R,\alpha)$ such that for any $T>0$ and $\delta<\delta^{*}=\min\left\{ 1,\frac{1}{2K_{1}\theta},\frac{1}{K_{4}\theta}\right\} $,
  \begin{equation}
    \mathbb{E}\left[\left|\hat{X}_{t}-\tilde{X}_{t}\right|\boldsymbol{1}_{[0,\nu_{R}]}(t)\right]
    \leq C(R,\alpha)\delta^{\min\{\eta_{F},1/2\}}. \notag
  \end{equation}
\end{Lemma}
\begin{proof}
  By \eqref{continuous STM2}, \eqref{continuous FBEM} and Fubini's theorem, we have
  \begin{equation}
    \mathbb{E}_{B}\left[\left|\hat{X}_{t}-\tilde{X}_{t}\right|\right]
    =\theta\int_{0}^{t}\mathbb{E}_{B}\left[\left|F\left(\tau_{n_{r}+1},\tilde{X}_{\tau_{n_{r}+1}}\right)-F\left(\tau_{n_{r}},\tilde{X}_{\tau_{n_{r}+1}}\right)\right|\right]dE_{r}. \label{3.19}
  \end{equation}
  Since for $\tau_{n_{r}}\leq t\leq\nu_{R}$, we have $\left|\tilde{X}_{\tau_{n_{r}}}\right|\leq R$. Then by Assumptions 3.4 and 3.1, we have
  \begin{align}
    \mathbb{E}_{B}\left[\left|F\left(\tau_{n_{r}},\tilde{X}_{\tau_{n_{r}}}\right)-F\left(\tau_{n_{r}+1},\tilde{X}_{\tau_{n_{r}+1}}\right)\right|\right]
    &\leq\mathbb{E}_{B}\left[\left|F\left(\tau_{n_{r}},\tilde{X}_{\tau_{n_{r}}}\right)-F\left(\tau_{n_{r}+1},\tilde{X}_{\tau_{n_{r}}}\right)\right|\right] \notag \\
    &\quad +\mathbb{E}_{B}\left[\left|F\left(\tau_{n_{r}+1},\tilde{X}_{\tau_{n_{r}}}\right)-F\left(\tau_{n_{r}+1},\tilde{X}_{\tau_{n_{r}+1}}\right)\right|\right] \notag \\
    &\leq K_{2}\left(1+\left|\tilde{X}_{\tau_{n_{r}}}\right|\right)\left|\tau_{n_{r}}-\tau_{n_{r}+1}\right|^{\eta_{F}}+C(R)\mathbb{E}_{B}\left[\left|\tilde{X}_{\tau_{n_{r}}}-\tilde{X}_{\tau_{n_{r}+1}}\right|\right] \notag \\
    &\leq K_{2}(1+R)\delta^{\eta_{F}}+C(R)\mathbb{E}_{B}\left[\left|\tilde{X}_{\tau_{n_{r}}}-\tilde{X}_{\tau_{n_{r}+1}}\right|\right]. \label{3.20}
  \end{align}
  By the assumption $F(t,0)=0$, \autoref{LL}, \autoref{thm: discrete STM, 2nd moment bound, sup} and the fact that $\left|\left(B_{(n_{r}+1)\delta}-B_{n_{r}\delta}\right)\right|\leq\delta^{1/2}$, we have
  \begin{align}
    \mathbb{E}_{B}\left[\left|\tilde{X}_{\tau_{n_{r}}}-\tilde{X}_{\tau_{n_{r}+1}}\right|\right]
    &\leq \mathbb{E}_{B}\left[\left|F\left(\tau_{n_{r}},\tilde{X}_{\tau_{n_{r}}}\right)\right|\right]\delta+\mathbb{E}_{B}\left[\left|F\left(\tau_{n_{r}+1},\tilde{X}_{\tau_{n_{r}+1}}\right)\right|\right]\delta \notag \\
    &\quad +\mathbb{E}_{B}\left[\left|G\left(\tau_{n_{r}},\tilde{X}_{\tau_{n_{r}}}\right)\right|\left|\left(B_{(n_{r}+1)\delta}-B_{n_{r}\delta}\right)\right|\right] \notag \\
    &=\mathbb{E}_{B}\left[\left|F\left(\tau_{n_{r}},\tilde{X}_{\tau_{n_{r}}}\right)-F\left(\tau_{n_{r}},0\right)\right|\right]\delta+\mathbb{E}_{B}\left[\left|F\left(\tau_{n_{r}+1},\tilde{X}_{\tau_{n_{r}+1}}\right)-F\left(\tau_{n_{r}},0\right)\right|\right]\delta \notag \\
    &\quad +\mathbb{E}_{B}\left[\left|G\left(\tau_{n_{r}},\tilde{X}_{\tau_{n_{r}}}\right)\right|\left|\left(B_{(n_{r}+1)\delta}-B_{n_{r}\delta}\right)\right|\right] \notag \\
    &\leq C(R)\mathbb{E}_{B}\left[\left|\tilde{X}_{\tau_{n_{r}}}\right|\right]\delta+C(R)\mathbb{E}_{B}\left[\left|\tilde{X}_{\tau_{n_{r}+1}}\right|\right]\delta+C(R)\mathbb{E}_{B}\left[\left|\tilde{X}_{\tau_{n_{r}}}\right|\right]\delta^{1/2} \notag \\
    &\leq C(R)\delta+C(R)\delta+C(R)\delta^{1/2} \notag \\
    &\leq C(R)\delta^{1/2}. \label{3.21}
  \end{align}
  Subsitute \eqref{3.21} into \eqref{3.20} gives 
  \begin{equation}
    \mathbb{E}_{B}\left[\left|F\left(\tau_{n_{r}},\tilde{X}_{\tau_{n_{r}}}\right)-F\left(\tau_{n_{r}+1},\tilde{X}_{\tau_{n_{r}+1}}\right)\right|\right]
    \leq K_{2}(1+R)\delta^{\eta_{F}}+C(R)\delta^{1/2}
    \leq C(R)\delta^{\min\{\eta_{F},1/2\}}. \notag 
  \end{equation}
  Then \eqref{3.19} becomes
  \begin{equation}
    \mathbb{E}_{B}\left[\left|\hat{X}_{t}-\tilde{X}_{t}\right|\right]
    \leq C(R)\delta^{\min\{\eta_{F},1/2\}}E_{t}. \notag
  \end{equation}
  Taking $\mathbb{E}_{D}$ on both sides and applying \autoref{Theorem:moment of inverse subordinator E},we have
  \begin{equation}
    \mathbb{E}\left[\left|\hat{X}_{t}-\tilde{X}_{t}\right|\boldsymbol{1}_{[0,\nu_{R}]}(t)\right]
    \leq C(R)\delta^{\min\{\eta_{F},1/2\}}\frac{1}{\Gamma(\alpha+1)}t^{\alpha}
    \leq C(R,\alpha)\delta^{\min\{\eta_{F},1/2\}}. \notag 
  \end{equation}
\end{proof}

\begin{Theorem} \label{thm: FBEM-STM probability}
  Let Assumptions 3.1-3.4 and 4.1 hold, and $\theta\in[1/2,1]$. 
  Then, for any given $\varepsilon>0$, there exists a positive constant $R_0$ such that for any $R\geq R_0$,
  we can find a positive number $\delta_{0}=\delta_{0}(R)$ such that whenever $\delta <\delta_0$, 
  \begin{equation}
    \mathbb{P}(\nu_{R}\leq T)\leq \frac{\varepsilon}{R^2}, \quad \text{for} \; T>0. \notag
  \end{equation}
\end{Theorem}
\begin{proof}
  For any $t\in[0,T]$, by the Itô's formula and \autoref{monotone}, we have
  \begin{align}
    \left|\hat{X}_{t\wedge\nu_{R}}\right|^{2}
    &=\left|x_{0}\right|^{2}+\int_{0}^{\tau_{n_{t}}\wedge\nu_{R}}\left(2\left\langle \hat{X}_{\tau_{n_{s}}},F\left(\tau_{n_{s}},\tilde{X}_{\tau_{n_{s}}}\right)\right\rangle +\left|G\left(\tau_{n_{s}},\tilde{X}_{\tau_{n_{s}}}\right)\right|^{2}\right)dE_{s} \notag \\
    &\quad +2\int_{0}^{\tau_{n_{s}}\wedge\nu_{R}}\left\langle \hat{X}_{\tau_{n_{s}}},G\left(\tau_{n_{s}},\tilde{X}_{\tau_{n_{s}}}\right)\right\rangle dB_{E_{s}} \notag \\
    &\leq \left|x_{0}\right|^{2}+2K_{1}\int_{0}^{\tau_{n_{t}}\wedge\nu_{R}}\left(1+\left|\tilde{X}_{\tau_{n_{s}}}\right|^{2}\right)dE_{s}+2\int_{0}^{\tau_{n_{s}}\wedge\nu_{R}}\left\langle \hat{X}_{\tau_{n_{s}}}-\tilde{X}_{\tau_{n_{s}}},F\left(\tau_{n_{s}},\tilde{X}_{\tau_{n_{s}}}\right)\right\rangle dE_{s} \notag \\
    &\quad +2\int_{0}^{\tau_{n_{s}}\wedge\nu_{R}}\left\langle \hat{X}_{\tau_{n_{s}}},G\left(\tau_{n_{s}},\tilde{X}_{\tau_{n_{s}}}\right)\right\rangle dB_{E_{s}} \notag 
  \end{align}
  Taking $\mathbb{E}_{B}$ on both sides, since the independence of $B$ and $E$ ensures that the time change $E$ does not affect the martingale property of $B$ under $\mathbb{P}_{B}$, we have
  \begin{align}
    \mathbb{E}_{B}\left[\left|\hat{X}_{t\wedge\nu_{R}}\right|^{2}\right]
    &\leq \left|x_{0}\right|^{2}+2K_{1}\mathbb{E}_{B}\left[\int_{0}^{\tau_{n_{t}}\wedge\nu_{R}}\left(1+\left|\tilde{X}_{\tau_{n_{s}}}\right|^{2}\right)dE_{s}\right] \notag \\
    &\quad +2\mathbb{E}_{B}\left[\int_{0}^{\tau_{n_{s}}\wedge\nu_{R}}\left|\hat{X}_{\tau_{n_{s}}}-\tilde{X}_{\tau_{n_{s}}}\right|\cdot\left|F\left(\tau_{n_{s}},\tilde{X}_{\tau_{n_{s}}}\right)\right|dE_{s}\right]. \notag
  \end{align}
  Using \autoref{LL}, \autoref{thm: discrete STM, 2nd moment bound, sup}, and $F(t,0)=0$ for $t\in[0,\nu_{R}]$, we have
  \begin{equation}
    \left|F\left(\tau_{n_{s}},\tilde{X}_{\tau_{n_{s}}}\right)\right|
    \leq \left|F\left(\tau_{n_{s}},\tilde{X}_{\tau_{n_{s}}}\right)-F\left(\tau_{n_{s}},\tilde{X}_{0}\right)\right|+\left|F\left(\tau_{n_{s}},\tilde{X}_{0}\right)\right|
    \leq C(R). \label{4.2}
  \end{equation}
  Taking $\mathbb{E}_{D}$ on both sides, using \eqref{4.2} and $(a-b+c)^{2}\leq2(|a-b|^{2}+|c|^{2})$, we obtain
  \begin{align}
    \mathbb{E}\left[\left|\hat{X}_{t\wedge\nu_{R}}\right|^{2}\right]
    &\leq \left|x_{0}\right|^{2}+2K_{1}\mathbb{E}\left[\int_{0}^{\tau_{n_{t}}\wedge\nu_{R}}\left(1+\left|\tilde{X}_{\tau_{n_{s}}}\right|^{2}\right)dE_{s}\right]+C(R)\mathbb{E}\left[\int_{0}^{\tau_{n_{s}}\wedge\nu_{R}}\left|\hat{X}_{\tau_{n_{s}}}-\tilde{X}_{\tau_{n_{s}}}\right|dE_{s}\right] \notag \\
    &\leq \left|x_{0}\right|^{2}+2K_{1}\mathbb{E}\left[E_{T}\right]+4K_{1}\mathbb{E}\left[\int_{0}^{\tau_{n_{t}}\wedge\nu_{R}}\left|\tilde{X}_{\tau_{n_{s}}}-\hat{X}_{\tau_{n_{s}}}\right|^{2}dE_{s}\right] \notag \\
    &\quad +4K_{1}\mathbb{E}_{B}\left[\int_{0}^{\tau_{n_{t}}\wedge\nu_{R}}\left|\hat{X}_{\tau_{n_{s}}}\right|^{2}dE_{s}\right]+C(R)\mathbb{E}\left[\int_{0}^{\tau_{n_{s}}\wedge\nu_{R}}\left|\hat{X}_{\tau_{n_{s}}}-\tilde{X}_{\tau_{n_{s}}}\right|dE_{s}\right]. \notag
  \end{align}
  Then applying $\left|\hat{X}_{\tau_{n_{s}}}-\tilde{X}_{\tau_{n_{s}}}\right|^{2}\leq2R\left|\hat{X}_{\tau_{n_{s}}}-\tilde{X}_{\tau_{n_{s}}}\right|$ and \autoref{cor: continuous FBEM-STM, pth moment bound, compact}, we obtain
  \begin{align}
  \mathbb{E}\left[\left|\hat{X}_{t\wedge\nu_{R}}\right|^{2}\right]
  &\leq \left|x_{0}\right|^{2}+2K_{1}\mathbb{E}\left[E_{T}\right]+C(R,\alpha)\delta^{\min\{\eta_{F},1/2\}}E_{T}+4K_{1}\mathbb{E}\int_{0}^{\tau_{n_{t}}}\left|\hat{X}_{\tau_{n_{s}}\wedge\nu_{R}}\right|^{2}dE_{s}. \notag 
  \end{align}
  Then
  \begin{equation}
  \underset{0\leq s\leq t}{\sup}\mathbb{E}\left[\left|\hat{X}_{s\wedge\nu_{R}}\right|^{2}\right]
  \leq \left|x_{0}\right|^{2}+2K_{1}\mathbb{E}\left[E_{T}\right]+C(R,\alpha)\delta^{\min\{\eta_{F},1/2\}}E_{T}+4K_{1}\mathbb{E}\int_{0}^{\tau_{n_{t}}}\underset{0\leq s\leq t}{\sup}\left|\hat{X}_{\tau_{n_{s}}\wedge\nu_{R}}\right|^{2}dE_{s}. \notag
  \end{equation}
  Setting $t=T$, using the Grönwall-type inequality (Lemma 6.3 in Chapter IX.6a of \cite{Jacod2003}), and applying \autoref{Theorem:moment of inverse subordinator E} and \autoref{Theorem: exponential moments of powers of inverse subordinator E}, we can easily get
  \begin{align}
    \mathbb{E}\left[\left|\hat{X}_{T\wedge\nu_{R}}\right|^{2}\right]
    &\leq \left[\left|x_{0}\right|^{2}+2K_{1}\mathbb{E}\left[E_{T}\right]+C(R,\alpha)\delta^{\min\{\eta_{F},1/2\}}\mathbb{E}\left[E_{T}\right]\right]\cdot\mathbb{E}\left[exp\left(4K_{1}E_{T}\right)\right] \notag \\
    &\leq \sum_{k=0}^{\infty}\frac{(4K_{1}T^{\alpha})^{k}}{\Gamma(\alpha k+1)}\left[\left|x_{0}\right|^{2}+\frac{2K_{1}}{\Gamma(\alpha+1)}T^{\alpha}+C(R,\alpha)\frac{1}{\Gamma(\alpha+1)}T^{\alpha}\delta^{\min\{\eta_{F},1/2\}}\right]. \notag 
  \end{align}
  Hence we have
  \begin{align}
      \mathbb{P}(\nu_{R}\leq T)
      &\leq \mathbb{P}(\hat{X}_{T\wedge\nu_{R}}\geq R,\nu_{R}\leq T)+\mathbb{P}(\tilde{X}_{T\wedge\nu_{R}}\geq R,\nu_{R}\leq T) \notag \\
      &\leq \frac{1}{R^{2}}\sum_{k=0}^{\infty}\frac{(4K_{1}T^{\alpha})^{k}}{\Gamma(\alpha k+1)}\left[\left|x_{0}\right|^{2}+\frac{2K_{1}}{\Gamma(\alpha+1)}T^{\alpha}+C(R,\alpha)\frac{1}{\Gamma(\alpha+1)}T^{\alpha}\delta^{\min\{\eta_{F},1/2\}}\right]. \notag
  \end{align}
  Now, for any given $\varepsilon>0$, we choose $R_{0}=R_{0}(\varepsilon)$ such that for $R\geq R_{0}$,
  \begin{equation}
    \sum_{k=0}^{\infty}\frac{(4K_{1}T^{\alpha})^{k}}{\Gamma(\alpha k+1)}\left[\left|x_{0}\right|^{2}+\frac{2K_{1}}{\Gamma(\alpha+1)}T^{\alpha}\right]\leq\frac{\varepsilon}{2}. \notag
  \end{equation}
  Then we choose $\delta_{0}=\delta_{0}(R)$ such that for $\delta<\delta_{0}$,
  \begin{equation}
    \sum_{k=0}^{\infty}\frac{(4K_{1}T^{\alpha})^{k}}{\Gamma(\alpha k+1)}C(R,\alpha)\frac{1}{\Gamma(\alpha+1)}T^{\alpha}\delta^{\min\{\eta_{F},1/2\}}\leq\frac{\varepsilon}{2}, \notag
  \end{equation}
  Whence $\mathbb{P}(\nu_{R}\leq T)\leq\frac{\varepsilon}{R^{2}}$.
\end{proof}

Now we are ready to prove the main theorem.
\begin{Theorem}\label{thm: convergence}
  Let Assumptions 3.1-3.4 and 4.1 hold.
  Assume $\theta\in[1/2,1]$, $R=\delta^{-1/2}$ and $\alpha\in(1/2,1)$.
  Then, there exist a constant $C(R,\alpha)\in (0,\infty)$ such that
  \begin{equation}
    \mathbb{E}\left[\underset{0\leq t\leq T}{\sup}\left|\tilde{X}_{t}-X_{t}\right|\right]
    \leq C(R,\alpha)\delta^{\min\{\eta_{F},\eta_{G},\alpha/2\}}. \notag
  \end{equation}
\end{Theorem}

\begin{proof}
  Define the stopping times 
  \begin{equation}
    \rho_{R}=\kappa_{R}\wedge\nu_{R}. \label{rho} 
  \end{equation}
  We turn the problem into seperate parts:
  \begin{align}
    \mathbb{E}\left[\underset{0\leq s\leq t}{\sup}\left|\tilde{X}_{s}-X_{s}\right|\right]
    &=\mathbb{E}\left[\underset{0\leq s\leq t}{\sup}\left|\tilde{X}_{s}-X_{s}\right|\boldsymbol{1}_{\{\kappa_{R}>t,\nu_{R}>t\}}\right]+\mathbb{E}\left[\underset{0\leq s\leq t}{\sup}\left|\tilde{X}_{s}-X_{s}\right|\boldsymbol{1}_{\{\kappa_{R}<t\ or\ \nu_{R}<t\}}\right] \notag \\
    &=\mathbb{E}\left[\underset{0\leq s\leq t}{\sup}\left|\tilde{X}_{s}-X_{s}\right|\boldsymbol{1}_{\{\rho_{R}>t\}}\right]+\mathbb{E}\left[\underset{0\leq s\leq t}{\sup}\left|\tilde{X}_{s}-X_{s}\right|\boldsymbol{1}_{\{\kappa_{R}<t\ or\ \nu_{R}<t\}}\right]  \notag \\
    &\leq \mathbb{E}\left[\underset{0\leq s\leq t}{\sup}\left|\tilde{X}_{s\wedge\rho_{R}}-X_{s\wedge\rho_{R}}\right|\right]+\mathbb{E}\left[\underset{0\leq s\leq t}{\sup}\left|\tilde{X}_{s}-X_{s}\right|\boldsymbol{1}_{\{\kappa_{R}<t\ or\ \nu_{R}<t\}}\right] \notag 
  \end{align}
  Applying the Young's inequality, we have
  \begin{align}
      \mathbb{E}\left[\underset{0\leq s\leq t}{\sup}\left|\tilde{X}_{s}-X_{s}\right|\right]
      &\leq \mathbb{E}\left[\underset{0\leq s\leq t}{\sup}\left|\tilde{X}_{s\wedge\rho_{R}}-\hat{X}_{s\wedge\rho_{R}}\right|\right]+\mathbb{E}\left[\underset{0\leq s\leq t}{\sup}\left|\hat{X}_{s\wedge\rho_{R}}-X_{s\wedge\rho_{R}}\right|\right] \notag \\
      &\quad +\frac{q}{2}\mathbb{E}\left[\underset{0\leq s\leq t}{\sup}\left|\tilde{X}_{s}-X_{s}\right|^{2}\right]+\frac{1}{2q}\mathbb{P}\left(\kappa_{R}<t\right)+\frac{1}{2q}\mathbb{P}\left(\nu_{R}<t\right) \notag \\
      &:=I_1+I_2+I_3+I_4+I_5, \label{error}
  \end{align}
  where
  \begin{align}
    I_{1}:=\mathbb{E}\left[\underset{0\leq s\leq t}{\sup}\left|\tilde{X}_{s\wedge\rho_{R}}-\hat{X}_{s\wedge\rho_{R}}\right|\right], 
    && I_{2}:=\mathbb{E}\left[\underset{0\leq s\leq t}{\sup}\left|\hat{X}_{s\wedge\rho_{R}}-X_{s\wedge\rho_{R}}\right|\right], \notag
  \end{align}
  and
  \begin{align}
    I_{3}:=\frac{q}{2}\mathbb{E}\left[\underset{0\leq s\leq t}{\sup}\left|\tilde{X}_{s}-X_{s}\right|^{2}\right], 
    && I_{4}:=\frac{1}{2q}\mathbb{P}\left(\kappa_{R}<t\right), 
    && I_{5}:=\frac{1}{2q}\mathbb{P}\left(\nu_{R}<t\right). \notag
  \end{align}
  Using stopping times \eqref{rho}, for any $0\leq s\leq t\leq T$, we have the truncated version of TCSDE \eqref{SDE}
    \begin{equation}
      X_{s\wedge\rho_{R}}:=x_{0}+\int_{0}^{s\wedge\rho_{R}}F(r,X_{r})dE_{r}+\int_{0}^{s\wedge\rho_{R}}G(r,X_{r})dB_{E_{r}}. \notag
    \end{equation}
  For $I_1$, by \autoref{cor: continuous FBEM-STM, pth moment bound, compact}, we have
    \begin{equation}
      \mathbb{E}\left[\underset{0\leq s\leq t}{\sup}\left|\tilde{X}_{s\wedge\rho_{R}}-\hat{X}_{s\wedge\rho_{R}}\right|\right]
      \leq C(R,\alpha)\delta^{\min\{\eta_{F},1/2\}}. \label{I1}
    \end{equation} 
    For $I_2$, we know that
    \begin{align}
      \left|\hat{X}_{s\wedge\rho_{R}}-X_{s\wedge\rho_{R}}\right|
      &\leq \left|\int_{0}^{s\wedge\rho_{R}}\left(F\left(\tau_{n_{r}},\tilde{X}_{\tau_{n_{r}}}\right)-F\left(r,X_{r}\right)\right)dE_{r}\right|+\left|\int_{0}^{s\wedge\rho_{R}}\left(G\left(\tau_{n_{r}},\tilde{X}_{\tau_{n_{r}}}\right)-G\left(r,X_{r}\right)\right)dB_{E_{r}}\right| \notag \\
      &\leq I_{21}+I_{22}, \notag
    \end{align}
    where
    \begin{equation}
      I_{21}:=
      \underset{0\leq s\leq t}{\sup}\left|\int_{0}^{s\wedge\rho_{R}}\left(F\left(\tau_{n_{r}},\tilde{X}_{\tau_{n_{r}}}\right)-F\left(r,X_{r}\right)\right)dE_{r}\right|, \notag \\
    \end{equation}
    and
    \begin{equation}
      I_{22}:=
      \underset{0\leq s\leq t}{\sup}\left|\int_{0}^{s\wedge\rho_{R}}\left(G\left(\tau_{n_{r}},\tilde{X}_{\tau_{n_{r}}}\right)-G\left(r,X_{r}\right)\right)dB_{E_{r}}\right|. \notag \\
    \end{equation}
    For any $0\leq s\leq t\leq T$, set
    \begin{equation}
      Z_{t}:=\underset{0\leq s\leq t}{\sup}|\tilde{X}_{s}-X_{s}|. \notag
    \end{equation}
    
    For $I_{21}$, applying the Cauchy-Schwarz inequality, the elementary inequality $(a+b+c)^{2}\leq3(a^{2}+b^{2}+c^{2})$, \autoref{LL} and \autoref{continuity}, we obtain
    \begin{align}
      \mathbb{E}_{B}\left[I_{21}^{2}\right]
      &\leq \mathbb{E}_{B}\left[\int_{0}^{t\wedge\rho_{R}}\left|\left(F\left(\tau_{n_{r}},\tilde{X}_{\tau_{n_{r}}}\right)-F\left(r,X_{r}\right)\right)^{2}\right|dE_{r}\right] \notag \\
      &\leq E_{T}\int_{0}^{t\wedge\rho_{R}}\mathbb{E}_{B}\left[\left|F\left(\tau_{n_{r}},\tilde{X}_{\tau_{n_{r}}}\right)-\left(F\left(r,X_{r}\right)\right)^{2}\right|\right]dE_{r} \notag \\
      &\leq 3E_{T}\int_{0}^{t\wedge\rho_{R}}\mathbb{E}_{B}\left[\left|F\left(\tau_{n_{r}},X_{r}\right)-F\left(r,X_{r}\right)\right|^{2}\right]dE_{r}\notag \\
      &\quad +3E_{T}\int_{0}^{t\wedge\rho_{R}}\mathbb{E}_{B}\left[\left|F\left(\tau_{n_{r}},X_{\tau_{n_{r}}}\right)-F\left(\tau_{n_{r}},X_{r}\right)\right|^{2}\right]dE_{r} \notag \\
      &\quad +3E_{T}\int_{0}^{t\wedge\rho_{R}}\mathbb{E}_{B}\left[\left|F\left(\tau_{n_{r}},\tilde{X}_{\tau_{n_{r}}}\right)-F\left(\tau_{n_{r}},X_{\tau_{n_{r}}}\right)\right|^{2}\right]dE_{r} \notag \\
      &\leq 3E_{T}K_{2}^{2}\int_{0}^{t\wedge\rho_{R}}\mathbb{E}_{B}\left[\left(1+\left|X_{r}\right|\right)^{2}\left|\tau_{n_{r}}-r\right|^{2\eta_{F}}\right]dE_{r}+E_{T}C(R)\int_{0}^{t\wedge\rho_{R}}\mathbb{E}_{B}\left[\left|X_{\tau_{n_{r}}}-X_{r}\right|^{2}\right]dE_{r} \notag \\
      &\quad +E_{T}C(R)\int_{0}^{t\wedge\rho_{R}}\mathbb{E}_{B}\left[\left|\tilde{X}_{\tau_{n_{r}}}-X_{\tau_{n_{r}}}\right|^{2}\right]dE_{r}. \label{I2_1}
    \end{align}
  Using the inequality $0\leq \tau_{n_{r}}-r\leq \delta$, \eqref{I2_1} becomes
  \begin{align}
      \mathbb{E}_{B}\left[I_{21}^{2}\right]
      &\leq 3E_{T}K_{2}^{2}\delta^{2\eta_{F}}\int_{0}^{t\wedge\rho_{R}}\mathbb{E}_{B}\left(1+\left|X_{r}\right|\right)^{2}dE_{r}+E_{T}C(R)\int_{0}^{t\wedge\rho_{R}}\mathbb{E}_{B}\left[\left|X_{\tau_{n_{r}}}-X_{r}\right|^{2}\right]dE_{r} \notag \\
      &\quad +E_{T}C(R)\int_{0}^{t\wedge\rho_{R}}\mathbb{E}_{B}\left[\left|\tilde{X}_{\tau_{n_{r}}}-X_{\tau_{n_{r}}}\right|^{2}\right]dE_{r} \notag \\
      &\leq 6E_{T}^{2}K_{2}^{2}\delta^{2\eta_{F}}\mathbb{E}_{B}\left[Y_{T}^{(2)}\right]+E_{T}C(R)\int_{0}^{t\wedge\rho_{R}}\mathbb{E}_{B}\left[\left|X_{\tau_{n_{r}}}-X_{r}\right|^{2}\right]dE_{r} \notag \\
      &\quad +E_{T}C(R)\int_{0}^{t\wedge\rho_{R}}\mathbb{E}_{B}\left[\left|\tilde{X}_{\tau_{n_{r}}}-X_{\tau_{n_{r}}}\right|^{2}\right]dE_{r}. \label{I2_2}
  \end{align} 
  By \autoref{exact any two difference}, we further obtain
  \begin{align}
    \mathbb{E}_B\left[\left|X_{\tau_{n_{r}}}-X_{r}\right|^{2}\right]
    &\leq C(h)\mathbb{E}_B\left[Y_{T}^{(2h)}\right]\left\{ \left(\tau_{n_{r}}-r\right)^{2\alpha}+\left(\tau_{n_{r}}-r\right)^{\alpha}\right\}       \notag \\
    &\leq C(h)\mathbb{E}_B\left[Y_{T}^{(2h)}\right]\left\{ \delta^{2\alpha}+\delta^{\alpha}\right\}. \label{I2_3}
  \end{align}
  Subsituting \eqref{I2_3} into \eqref{I2_2} gives
  \begin{align}
    \mathbb{E}_{B}\left[I_{21}^{2}\right]
    &\leq 6E_{T}^{2}K_{2}^{2}\delta^{2\eta_{F}}\mathbb{E}_{B}\left[Y_{T}^{(2)}\right]+E_{T}C(R)\int_{0}^{t\wedge\rho_{R}}\mathbb{E}_{B}\left[Y_{T}^{(2h)}\right]\left\{ \delta^{2\alpha}+\delta^{\alpha}\right\} dE_{r} \notag \\
    &\quad +E_{T}C(R)\int_{0}^{t\wedge\rho_{R}}\mathbb{E}_{B}\left[Z_{r}^{2}\right]dE_{r} \notag \\
    &\leq 6E_{T}^{2}K_{2}^{2}\delta^{2\eta_{F}}\mathbb{E}_{B}\left[Y_{T}^{(2)}\right]+E_{T}^{2}C(R)\left(2\delta^{\alpha}\right)\mathbb{E}_{B}\left[Y_{T}^{(2h)}\right]+E_{T}C(R)\int_{0}^{t\wedge\rho_{R}}\mathbb{E}_{B}\left[Z_{r}^{2}\right]dE_{r} \notag \\
    &\leq C(R)E_{T}^{2}\mathbb{E}_{B}\left[Y_{T}^{(2h)}\right]\delta^{\min\{2\eta_{F},\alpha\}}+C(R)E_{T}\int_{0}^{t\wedge\rho_{R}}\mathbb{E}_{B}\left[Z_{r}^{2}\right]dE_{r}. \label{I21}
  \end{align}

  For $I_{22}$, applying the Burkholder-Davis-Gundy inequality \eqref{BDG}, \autoref{LL} and \autoref{continuity}, similar to the proof for $I_{21}$, we obtain
  \begin{align}
    \mathbb{E}_{B}\left[I_{22}^{2}\right]
    &\leq b_{2}\mathbb{E}_{B}\int_{0}^{t\wedge\rho_{R}}\left(G\left(\tau_{n_{r}},\tilde{X}_{\tau_{n_{r}}}\right)-G\left(r,X_{r}\right)\right)^{2}dE_{r} \notag \\
    &\leq 3b_{2}K_{3}^{2}\delta^{2\eta_{G}}\int_{0}^{t\wedge\rho_{R}}\mathbb{E}_{B}\left(1+\left|X_{r}\right|\right)^{2}dE_{r}+C(R)\int_{0}^{t\wedge\rho_{R}}\mathbb{E}_{B}\left[\left|X_{\tau_{n_{r}}}-X_{r}\right|^{2}\right]dE_{r} \notag \\
    &\quad +C(R)\int_{0}^{t\wedge\rho_{R}}\mathbb{E}_{B}\left[\left|\tilde{X}_{\tau_{n_{r}}}-X_{\tau_{n_{r}}}\right|^{2}\right]dE_{r} \notag \\
    &\leq 6b_{2}E_{T}K_{3}^{2}\delta^{2\eta_{G}}\mathbb{E}_{B}\left[Y_{T}^{(2)}\right]+C(R)\int_{0}^{t\wedge\rho_{R}}\mathbb{E}_{B}\left[\left|X_{\tau_{n_{r}}}-X_{r}\right|^{2}\right]dE_{r}  \notag \\
    &\quad +C(R)\int_{0}^{t\wedge\rho_{R}}\mathbb{E}_{B}\left[\left|\tilde{X}_{\tau_{n_{r}}}-X_{\tau_{n_{r}}}\right|^{2}\right]dE_{r}. \label{I3_1}
  \end{align}
  Subsituting \eqref{I2_3} into \eqref{I3_1} gives
  \begin{equation}
    \mathbb{E}_{B}\left[I_{22}^{2}\right]
    \leq C(R)E_{T}\mathbb{E}_{B}[Y_{T}^{(2h)}]\delta^{\min\{2\eta_{G},\alpha\}}+C(R)\int_{0}^{t\wedge\rho_{R}}\mathbb{E}_{B}\left[Z_{r}^{2}\right]dE_{r}. \label{I22}
  \end{equation}
  Putting together estimates \eqref{I21} and \eqref{I22} gives
  \begin{align}
    \mathbb{E}_{B}\left[Z_{r\wedge\rho_{R}}^{2}\right]
    &\leq C(R)\left(E_{T}^{2}+E_{T}\right)\mathbb{E}_{B}\left[Y_{T}^{(2h)}\right]\delta^{\min\{2\eta_{F},2\eta_{F},\alpha\}}+C(R)\left(1+E_{T}\right)\int_{0}^{t\wedge\rho_{R}}\mathbb{E}_{B}\left[Z_{r}^{2}\right]dE_{r}. \label{I2}
  \end{align}
  Applying the Grönwall-type inequality (Lemma 6.3 in Chapter IX.6a of \cite{Jacod2003}) and taking $\mathbb{E}_{D}$ on both sides with $t\wedge \rho_R=T$ gives
  \begin{equation}
    \mathbb{E}\left[I_{2}^{2}\right]
    \leq C(R)\left(\mathbb{E}\left[E_{T}^{2}\right]+\mathbb{E}\left[E_{T}\right]\right)
    \mathbb{E}\left[Y_{T}^{(2h)}\right]\delta^{\min\{2\eta_{F},2\eta_{G},\alpha\}}
    \cdot \mathbb{E}\left[exp\left(C(R)E_{T}^{2}\right)\right]
    \cdot \mathbb{E}\left[exp\left(C(R)E_{T}\right)\right].\notag
  \end{equation}
  Applying \autoref{Theorem:moment of inverse subordinator E}, \autoref{thm: exact pth moment, probability}, we know that the assumption $\alpha\in(1/2,1)$ allows us to pick $r=2<1/(1-\alpha)$ for \autoref{Theorem: exponential moments of powers of inverse subordinator E}, we have
  \begin{equation}
    \mathbb{E}\left[I_{2}^{2}\right]
    \leq \mathbb{E}\left[V_{1}\right]\cdot\mathbb{E}\left[V_{2}\right], \notag 
  \end{equation}
  where
  \begin{equation}
    \mathbb{E}\left[V_{1}\right]
    =C(R)\left(\frac{2}{\Gamma(2\alpha+1)}T^{2\alpha}+\frac{1}{\Gamma(2\alpha+1)}T^{\alpha}\right)\left(1+2^{h-1}\sum_{k=0}^{\infty}\frac{(2hKt^{\alpha})^{k}}{\Gamma(\alpha k+1)}\mathbb{E}\left[\left|x_{0}\right|^{2h}\right]\right)\delta^{\min\{2\eta_{F},2\eta_{G},\alpha\}}, \notag 
  \end{equation}
  and
  \begin{equation}
    \mathbb{E}\left[V_{2}\right]
    =\sum_{k=0}^{\infty}\frac{C(R)^{k}}{k!}\frac{\Gamma(2k+1)}{\Gamma(2\alpha k+1)}T^{2\alpha k}
    \cdot \sum_{k=0}^{\infty}\frac{(C(R)T^{\alpha})^{k}}{\Gamma(\alpha k+1)}. \notag
  \end{equation}
  Simplifying the constants we have 
  $\mathbb{E}\left[V_{1}\right]\leq C(R,\alpha)\delta^{\min\{2\eta_{F},2\eta_{F},\alpha\}}$ and $\mathbb{E}\left[V_{2}\right]\leq C(R,\alpha)< \infty$, which leads us to $\mathbb{E}\left[I_{2}^{2}\right]\leq C(R,\alpha)\delta^{\min\{2\eta_{F},2\eta_{G},\alpha\}}.$

  For $I_3$, by \autoref{thm: discrete STM, 2nd moment bound, sup} and taking $\mathbb{E}_{D}$, we can easily get 
  \begin{equation}
    \underset{0\leq s\leq t}{\sup}\mathbb{E}\left[\left|\tilde{X}_{s}\right|^{2}\right]	\leq C(t). \notag 
  \end{equation}
  Choosing $q=\delta ^{1/2}$ and $R=1/q$, then by \autoref{thm: exact pth moment, probability} we have
  \begin{equation}
    \frac{q}{2}\mathbb{E}\left[\underset{0\leq s\leq t}{\sup}\left|\tilde{X}_{s}-X_{s}\right|^{2}\right]
    \leq q\mathbb{E}\left[\underset{0\leq s\leq t}{\sup}\left|X_{s}\right|^{2}+\underset{0\leq s\leq t}{\sup}\left|\tilde{X}_{s}\right|^{2}\right]
    \leq C(\alpha,t)\delta^{1/2}.   \label{I3}
  \end{equation}
  For $I_4$, applying \autoref{thm: exact pth moment, probability} we have
  \begin{equation}
    \frac{1}{2q}\mathbb{P}\left(\kappa_{R}<t\right)
    = \frac{1}{2q}\frac{1}{R^{2}}\sum_{k=0}^{\infty}\frac{(2K_1t^{\alpha})^{k}}{\Gamma(\alpha k+1)}\mathbb{E}\left[\left|x_{0}\right|^{2}\right]\leq C(R,\alpha)\delta^{1/2}, \label{I4}
  \end{equation}
  For $I_5$, given an $\varepsilon>0$, by \autoref{thm: FBEM-STM probability} we can choose $\delta$ such that
  \begin{equation}
    \frac{1}{2q}\mathbb{P}\left(\nu_{R}<t\right)
    \leq \frac{1}{2q}\frac{1}{R^{2}}C(R,\alpha)\delta^{\min\{\eta_{F},1/2\}}
    \leq \frac{\varepsilon}{2}\delta^{1/2}. \label{I5}
  \end{equation}
  Finally, putting together \eqref{error}, \eqref{I1}, \eqref{I2}, \eqref{I3}, \eqref{I4} and \eqref{I5} with $t\wedge \rho_R=T$, 
  we obtain the inequality
  \begin{align}
    \mathbb{E}\left[\underset{0\leq t\leq T}{\sup}\left|\tilde{X}_{t}-X_{t}\right|\right]
    &\leq C(R,\alpha)\delta^{\min\{\eta_{F},1/2\}}+C(R,\alpha)\delta^{\min\{\eta_{F},\eta_{G},\alpha/2\}}+C(\alpha, T)\delta^{1/2}+C(R,\alpha)\delta^{1/2}+\frac{\varepsilon}{2}\delta^{1/2} \notag \\
    &\leq C(R,\alpha)\delta^{\min\{\eta_{F}, \eta_{G}, \alpha/2\}}. \notag
  \end{align}
  The proof is now complete.
\end{proof}

\section{Mittag-Leffler stability}
Having established the strong convergence of ST method, we will proceed to the moment stability of the underlying TCSDE \eqref{SDE}. We propose a new stability definition for TCSDE in the sense of Mittag-Leffler and consider Lyapunov's direct method. Suppose $F(0,0)=G(0,0)=0$.
\begin{definition}(Mittag-Leffler Stability)
  The solution of TCSDE \eqref{SDE} is said to be $p$th moment Mittag-Leffler stable if for any initial value $x_0\in \mathbb{R}$ there exist $C,\lambda>0$ and $0<\alpha<1$ such that
  \begin{equation}
    \mathbb{E}\left[\left|X_{t}\right|^{p}\right]\leq C\left|x_{0}\right|^{p}\mathbb{E}\left[e^{-\upsilon E_{t}}\right],
  \end{equation}
  or equivalently,
  \begin{equation}
    \mathbb{E}\left[\left|X_{t}\right|^{p}\right]\leq C\left|x_{0}\right|^{p}E_{\alpha}\left(-\upsilon t^{\alpha}\right),\quad t,p\geq 0,
  \end{equation}
  where 
  \begin{equation}
    E_{\alpha}(z)=\sum_{k=0}^{\infty}\frac{z^{k}}{\Gamma(\alpha k+1)}
  \end{equation}
  is the one parameter Mittag-Leffler function.
\end{definition}

\begin{Theorem}(Mittag-Leffler stability)\label{Mittag-Leffler Stability}
  Let $p,c_{1},c_{2}$ and $c_{3}$ be positive constants. 

  If $V\in C^{1,2}(\mathbb{R}_{+}\times\mathbb{R};\mathbb{R}_{+})$ satisfies:

  1. $c_{1}|x|^{p}\leq V(t,x)\leq c_{2}|x|^{p}$.

  2. $L_{1}V(t,x)+L_{2}V(t,x)\leq-c_3 V(t,x).$

  Here 
  \begin{equation}
      L_{1}V(t,x)=V_{t}(t,x)+V_{x}(t,x)F(t,x)
  \end{equation}
  and
  \begin{equation}
    L_{2}V(t,x)=V_{t}(t,x)+V_{x}(t,x)G(t,x)+\frac{1}{2}V_{xx}(t,x)G^{2}(t,x).
  \end{equation}
  Then for all $x_{0}\in R$,
  \begin{equation}
    \mathbb{E}\left[\left|X_{t}\right|^{p}\right]\leq\frac{c_{2}}{c_{1}}\left|x_{0}\right|^{p}E_{\alpha}\left(-c_{3}t^{\alpha}\right).
  \end{equation}
  That is, the solution of TCSDE \eqref{SDE}  is $p$th moment Mittag-Leffler stable.
\end{Theorem}

\begin{proof}
  Define $\zeta_{m}=m\wedge inf\left\{ t\geq0:\left|\int_{0}^{\kappa_{R}\wedge t}V_{x}(s,X_{s})\sigma(s,X_{s})dB_{E_{s}}\right|\geq m\right\}$.
  Applying the time-changed Itô's formula in \cite{2011Kobayashi} to $E_{\alpha}(c_{3}t^{\alpha})V(t,x)$.
  Then for any $t\geq 0$, 
  \begin{align}
    E_{\alpha}\left(c_{3}\left(t\wedge\kappa_{R}\wedge\zeta_{m}\right)^{\alpha}\right)V\left(t\wedge\kappa_{R}\wedge\zeta_{m},X_{t\wedge\kappa_{R}\wedge\zeta_{m}}\right)
    &=V(0,x_{0})+\int_{0}^{t\wedge\kappa_{R}\wedge\zeta_{m}}e^{c_{3}s}\left[c_{3}V(s,X_{s})+L_{1}V(s,X_{s})\right]ds \notag \\
    &\quad +\int_{0}^{t\wedge\kappa_{R}\wedge\zeta_{m}}e^{c_{3}s}L_{2}V(s,X_{s})dE_{s} \notag \\
    &\quad +\int_{0}^{t\wedge\kappa_{R}\wedge\zeta_{m}}e^{c_{3}s}V_{x}(s,X_{s})\sigma(s,X_{s})dB_{E_{s}}. \notag
  \end{align}
  Since $\int_{0}^{t\wedge\kappa_{R}\wedge\zeta_{m}}e^{c_{3}s}V_{x}(s,X_{s})\sigma(s,X_{s})dB_{E_{s}}$ is a mean zero martingale, taking expectations on both sides and applying conditions (2) and (3) gives
  \begin{align}
    \mathbb{E}\left[E_{\alpha}\left(c_{3}\left(t\wedge\kappa_{R}\wedge\zeta_{m}\right)^{\alpha}\right)V\left(t\wedge\kappa_{R}\wedge\zeta_{m},X_{t\wedge\kappa_{R}\wedge\zeta_{m}}\right)\right]
    &=V(0,x_{0})+\mathbb{E}\left[\int_{0}^{t\wedge\kappa_{R}\wedge\zeta_{m}}e^{c_{3}s}L_{2}V(s,X_{s})dE_{s}\right] \notag \\
    &\quad +\mathbb{E}\left[\int_{0}^{t\wedge\kappa_{R}\wedge\zeta_{m}}e^{c_{3}s}\left[c_{3}V(s,X_{s})+L_{1}V(s,X_{s})\right]ds\right] \notag \\
    &\leq V(0,x_{0}). \notag
  \end{align}
  By condition (1),
  \begin{align}
    c_{1}\mathbb{E}\left[\left|X_{t\wedge\kappa_{R}\wedge\zeta_{m}}\right|^{p}E_{\alpha}\left(c_{3}\left(t\wedge\kappa_{R}\wedge\zeta_{m}\right)^{\alpha}\right)\right]
    &\leq \mathbb{E}\left[E_{\alpha}\left(c_{3}\left(t\wedge\kappa_{R}\wedge\zeta_{m}\right)^{\alpha}\right)V\left(t\wedge\kappa_{R}\wedge\zeta_{m},X_{t\wedge\kappa_{R}\wedge\zeta_{m}}\right)\right] \notag \\
    &\leq V(0,x_{0})\leq c_{2}|x_{0}|^{p}, \notag
  \end{align}
  which implies
  \begin{equation}
    \mathbb{E}\left[\left|X_{t\wedge\kappa_{R}\wedge\zeta_{m}}\right|^{p}E_{\alpha}\left(c_{3}\left(t\wedge\kappa_{R}\wedge\zeta_{m}\right)^{\alpha}\right)\right]
    \leq\frac{c_{2}}{c_{1}}\left|x_{0}\right|^{p}. \notag
  \end{equation}
  Consequently,
  \begin{equation}
  \mathbb{E}\left[I_{[0,\kappa_{R}\wedge\zeta_{m}]}(t)\left|X_{t}\right|^{p}E_{\alpha}(c_{3}t^{\alpha})\right]
  \leq\frac{c_{2}}{c_{1}}\left|x_{0}\right|^{p}. \notag
  \end{equation}
  By the monotone convergence theorem it holds that 
  \begin{align}
    \mathbb{E}\left[\left|X_{t}\right|^{p}E_{\alpha}\left(c_{3}t^{\alpha}\right)\right]
    &=\underset{R,m\rightarrow\infty}{lim}\mathbb{E}\left[I_{[0,\kappa_{R}\wedge\zeta_{m}]}(t)\left|X_{t}\right|^{p}E_{\alpha}\left(c_{3}t^{\alpha}\right)\right] \notag \\
    &\leq\frac{c_{2}}{c_{1}}\left|x_{0}\right|^{p}. \notag
  \end{align}
  Rearranging the equation we obtain the desired result.
\end{proof}

\begin{remark}
  \textnormal{
  (1)  
  The two definitions of Mittag-Leffler stability in Definition 5.1 are equivalent. The latter emphasizes the effect of parameter $\alpha$ on the stability. 
  }

  \textnormal{
  (2) This theorem represents a significant generalization of well-known exponential stability results for classical SDEs. For $\alpha\in(0,1)$, stability exhibits a slower decay rate governed by the Mittag-Leffler function. 
  The smaller the $\alpha$, the larger the jumps in the subordinator, which induces longer waiting periods in the system's evolution. This directly postpones the emergence of stability, as the path to stability is frequently interrupted.
  In the classical limit $\alpha \rightarrow 1$, the Mittag-Leffler stability reduces to exponential stability, thus providing a unified framework that encapsulates both standard exponential stability ($\alpha=1$) and the slower, algebraic-like decay typical of subdiffusive systems ($\alpha<1$).
  }

   \textnormal{
  (3) The duality principle is not just ineffective for analyzing the case of time-dependent coefficients, but particularly ineffective for analyzing Mittag-Leffler stability. Mittag-Leffler stability is an emergent property of the coupled system $\left(X_t, E_t\right)$ evolving in physical time. The duality principle, by attempting to decouple the system and change the time scale, destroys the very context in which this stability manifests.
  }

\end{remark}

\section{Numerical Mittag-Leffler stability}

\begin{definition}\label{numerical stability definition}
  For a given step size $\delta>0$, the ST method is said to be mean-square Mittag-Leffler stable on the TCSDE \eqref{SDE} if for any initial value $x_0\in \mathbb{R}$ there exist $C,\gamma>0$ and $0<\alpha<1$ such that
  \begin{equation}
    \mathbb{E}\left[\left|\tilde{X}_{\tau_{n}}\right|^{2}\right]\leq C\left|x_{0}\right|^{2}E_{\alpha}\left(-\gamma(n\delta)^{\alpha}\right),\quad p,t\geq 0,
  \end{equation}
  where 
  \begin{equation}
    E_{\alpha}(z)=\sum_{k=0}^{\infty}\frac{z^{k}}{\Gamma(\alpha k+1)} \notag
  \end{equation}
  is the one parameter Mittag-Leffler function.
\end{definition}

We need the following assumption.
\begin{Assumption}\label{coercive}
  There exists a positive constant $\lambda$ such that for all $x\in\mathbb{R}^{d}$, 
  \begin{equation}
    \left\langle x,F(t,x)\right\rangle +\frac{1}{2}|G(t,x)|^{2}\leq-\lambda|x|^{2}. 
  \end{equation}
\end{Assumption}

%Next, we give a theorem on the asymptotically stability of the ST method in mean square. 
%This theorem shows that under certain nonlinear conditions, for all step sizes $\delta>0$, the ST method with $\theta\in(\frac{1}{2},1]$ preserves the asymptotic mean square stability of the SDE \eqref{SDE}. Under a stronger assumption, asymptotic mean square stability still holds for the ST method with $\theta\in[0,\frac{1}{2})$, but with restricted step size.

\begin{Theorem}\label{numerical stability}
  Let \autoref{coercive} hold. Assume that there exists a constant $K_5$ such that 
  \begin{equation}
    \left|F(t,x)\right|^{2}\leq K_5|x|^{2}, \label{K5}
  \end{equation}
  1. If $\theta\in[1/2,1]$, then the ST method applied to \eqref{SDE} is mean-square Mittag-Leffler stable for all $\delta>0$.

  2. If $\theta\in[0,1/2)$, then the ST method applied to \eqref{SDE} is mean-square Mittag-Leffler stable for all $0<\delta<\frac{5\lambda}{2K_{5}(1-2\theta)}$.
\end{Theorem}

\begin{proof}
  Using \eqref{discrete STM} we have
  \begin{equation}
  \left|\tilde{X}_{\tau_{n+1}}-\theta\delta F(\tau_{n+1},\tilde{X}_{\tau_{n+1}})\right|^{2}=\left|\tilde{X}_{\tau_{n}}+(1-\theta)\delta F\left(\tau_{n},\tilde{X}_{\tau_{n}}\right)+G\left(\tau_{n},\tilde{X}_{\tau_{n}}\right)\left(B_{(n+1)\delta}-B_{n\delta}\right)\right|^{2}. \notag
  \end{equation}
  Thus
  \begin{align}
    &\left|\tilde{X}_{\tau_{n+1}}\right|^{2}+\left|\theta\delta F\left(\tau_{n+1},\tilde{X}_{\tau_{n+1}}\right)\right|^{2} -2\theta\delta\left\langle \tilde{X}_{\tau_{n+1}},F\left(\tau_{n+1},\tilde{X}_{\tau_{n+1}}\right)\right\rangle \notag \\
    &=\left|\tilde{X}_{\tau_{n}}\right|^{2}+\left|(1-\theta)\delta F\left(\tau_{n},\tilde{X}_{\tau_{n}}\right)\right|^{2}+\left|G\left(\tau_{n},\tilde{X}_{\tau_{n}}\right)\left(B_{(n+1)\delta}-B_{n\delta}\right)\right|^{2} \notag \\
    &\quad +2(1-\theta)\delta\left\langle \tilde{X}_{\tau_{n}},F\left(\tau_{n},\tilde{X}_{\tau_{n}}\right)\right\rangle +2\left\langle \tilde{X}_{\tau_{n}},G\left(\tau_{n},\tilde{X}_{\tau_{n}}\right)\left(B_{(n+1)\delta}-B_{n\delta}\right)\right\rangle \notag \\
    &\quad +2(1-\theta)\delta\left\langle F\left(\tau_{n},\tilde{X}_{\tau_{n}}\right),G\left(\tau_{n},\tilde{X}_{\tau_{n}}\right)\left(B_{(n+1)\delta}-B_{n\delta}\right)\right\rangle. \notag
  \end{align} 
  Taking $\mathbb{E}_{B}$ and $\mathbb{E}_{D}$ consecutively, we have
  \begin{align}
    \mathbb{E}\left[\left|\tilde{X}_{\tau_{n+1}}\right|^{2}\right]
    &=\mathbb{E}\left[\left|\tilde{X}_{\tau_{n}}\right|^{2}\right]+\mathbb{E}\left[\left|(1-\theta)\delta F\left(\tau_{n},\tilde{X}_{\tau_{n}}\right)\right|^{2}\right]+\mathbb{E}\left[\left|G\left(\tau_{n},\tilde{X}_{\tau_{n}}\right)\left(B_{(n+1)\delta}-B_{n\delta}\right)\right|^{2}\right] \notag \\
    &\quad -\mathbb{E}\left[\left|\theta\delta F\left(\tau_{n+1},\tilde{X}_{\tau_{n+1}}\right)\right|^{2}\right]+2\theta\delta\mathbb{E}\left\langle \tilde{X}_{\tau_{n+1}},F\left(\tau_{n+1},\tilde{X}_{\tau_{n+1}}\right)\right\rangle  \notag \\
    &\quad +2(1-\theta)\delta\mathbb{E}\left\langle \tilde{X}_{\tau_{n}},F\left(\tau_{n},\tilde{X}_{\tau_{n}}\right)\right\rangle 
    +2\mathbb{E}\left\langle \tilde{X}_{\tau_{n}},G\left(\tau_{n},\tilde{X}_{\tau_{n}}\right)\left(B_{(n+1)\delta}-B_{n\delta}\right)\right\rangle  \notag \\
    &\quad +2(1-\theta)\delta\mathbb{E}\left\langle F\left(\tau_{n},\tilde{X}_{\tau_{n}}\right),G\left(\tau_{n},\tilde{X}_{\tau_{n}}\right)\left(B_{(n+1)\delta}-B_{n\delta}\right)\right\rangle.
  \end{align}
  Notice that, by definitions, $F\left(\tau_{n},\tilde{X}_{\tau_{n}}\right)$'s are measurable w.r.t. $\mathcal{\mathcal{F}}_{n}:=\sigma\left\{ D_{t},B_{t},0\leq t\leq n\delta\right\}$
  and $\mathcal{\mathcal{F}}_{n}$ is independent of $B_{(n+1)\delta}-B_{n\delta}$. 
  Now, for a given $n\geq 0$, the following inequalities hold
  \begin{equation}
    \mathbb{E}\left[\left|G\left(\tau_{n},\tilde{X}_{\tau_{n}}\right)\left(B_{(n+1)\delta}-B_{n\delta}\right)\right|^{2}\right]=\delta\mathbb{E}\left[\left|G\left(\tau_{n},\tilde{X}_{\tau_{n}}\right)\right|^{2}\right], \label{G}
  \end{equation}
  \begin{equation}
    \mathbb{E}\left\langle \tilde{X}_{\tau_{n}},G\left(\tau_{n},\tilde{X}_{\tau_{n}}\right)\left(B_{(n+1)\delta}-B_{n\delta}\right)\right\rangle =0, \label{X,G}
  \end{equation}
  and
  \begin{equation}
    \mathbb{E}\left\langle F\left(\tau_{n},\tilde{X}_{\tau_{n}}\right),G\left(\tau_{n},\tilde{X}_{\tau_{n}}\right)\left(B_{(n+1)\delta}-B_{n\delta}\right)\right\rangle =0.\label{F,G}
  \end{equation}
  Thus following \autoref{coercive} and \eqref{G}-\eqref{F,G} one finds that
  \begin{align}
    \mathbb{E}\left[\left|\tilde{X}_{\tau_{n+1}}\right|^{2}\right]
    &\leq \mathbb{E}\left[\left|\tilde{X}_{\tau_{n}}\right|^{2}\right]+\mathbb{E}\left[\left|(1-\theta)\delta F\left(\tau_{n},\tilde{X}_{\tau_{n}}\right)\right|^{2}\right]+\delta\mathbb{E}\left[\left|G\left(\tau_{n},\tilde{X}_{\tau_{n}}\right)\right|^{2}\right] \notag \\
    &\quad -\mathbb{E}\left[\left|\theta\delta F\left(\tau_{n+1},\tilde{X}_{\tau_{n+1}}\right)\right|^{2}\right]-2\theta\delta\lambda\mathbb{E}\left[\left|\tilde{X}_{\tau_{n+1}}\right|^{2}\right]
    -2(1-\theta)\delta\lambda\mathbb{E}\left[\left|\tilde{X}_{\tau_{n}}\right|^{2}\right]. \notag
  \end{align}
  Consequently by \eqref{K5},
  \begin{align}
    \mathbb{E}\left[\left|\tilde{X}_{\tau_{n+1}}\right|^{2}\right]
    &\leq\mathbb{E}\left[\left|\tilde{X}_{\tau_{n}}\right|^{2}\right]+(1-\theta)^{2}\delta^{2}K_{5}\mathbb{E}\left[\left|\tilde{X}_{\tau_{n}}\right|^{2}\right]-\frac{1}{2}\delta\lambda\mathbb{E}\left[\left|\tilde{X}_{\tau_{n}}\right|^{2}\right] \notag \\
    &\quad -\theta^{2}\delta^{2}K_{5}\mathbb{E}\left[\left|\tilde{X}_{\tau_{n+1}}\right|^{2}\right]-2\theta\delta\lambda\mathbb{E}\left[\left|\tilde{X}_{\tau_{n+1}}\right|^{2}\right] 
    -2(1-\theta)\delta\lambda\mathbb{E}\left[\left|\tilde{X}_{\tau_{n}}\right|^{2}\right] \notag \\
    &\leq \frac{1+(1-\theta)^{2}\delta^{2}K_{5}-\frac{1}{2}\delta\lambda-2(1-\theta)\delta\lambda}{1+\theta^{2}\delta^{2}K_{5}+2\theta\delta\lambda}\mathbb{E}\left[\left|\tilde{X}_{\tau_{n}}\right|^{2}\right].
  \end{align}
  By iteration, we have
  \begin{equation}
    \mathbb{E}\left[\left|\tilde{X}_{\tau_{n}}\right|^{2}\right]\leq\left(\varphi(\delta)\right)^{n}\mathbb{E}\left[\left|\tilde{X}_{\tau_{0}}\right|^{2}\right], \notag 
  \end{equation}
  where
  \begin{equation}
    \varphi(\delta)=\frac{1+(1-\theta)^{2}\delta^{2}K_{5}-\frac{1}{2}\delta\lambda-2(1-\theta)\delta\lambda}{1+\theta^{2}\delta^{2}K_{5}+2\theta\delta\lambda}. \notag
  \end{equation}
  If $\mathbb{E}\left[\left|\tilde{X}_{\tau_{n}}\right|^{2}\right]$ is contractive, i.e. $\varphi(\delta)<1$,
  then 
  \begin{align}
    \mathbb{E}\left[\left|\tilde{X}_{\tau_{n}}\right|^{2}\right]
    &\leq\mathbb{E}\left[\left|\tilde{X}_{\tau_{0}}\right|^{2}\right]E_{\alpha}\left(-\left(-\frac{ln\varphi(\delta)}{\delta}\right)^{\alpha}(n\delta)^{\alpha}\right) \notag \\
    &=\mathbb{E}\left[\left|\tilde{X}_{\tau_{0}}\right|^{2}\right]E_{\alpha}(-\gamma(n\delta)^{\alpha}), \label{6.10}
  \end{align}
  where $\gamma=\left(-\frac{ln\varphi(\delta)}{\delta}\right)^{\alpha}$, which is the form in Definition 6.1.

  Now we discuss the range of step size. We know $\varphi(\delta)<1$ gives
  \begin{equation}
    (1-\theta)^{2}\delta^{2}K_{5}-\frac{1}{2}\delta\lambda-2(1-\theta)\delta\lambda-\theta^{2}\delta^{2}K_{5}-2\theta\delta\lambda<0. \notag
  \end{equation}
  Simplifying above gives
  \begin{equation}
    \delta K_{5}(1-2\theta)<\frac{5}{2}\lambda. \label{range}
  \end{equation}
  When $\theta\geq\frac{1}{2}$, \eqref{range} holds for any $\delta>0$.  \notag
  When $\theta<\frac{1}{2}$, we have
  \begin{equation}
    \delta<\frac{5\lambda}{2K_{5}(1-2\theta)}. \notag
  \end{equation}
\end{proof}

\begin{remark}
  \textnormal{
  For stability index $\alpha=1$, the generalized Mittag-Leffler stability described by equation \eqref{6.10} reduces precisely to the well-known exponential stability.  
  Since $E_1(z)=e^z$, $E_{\alpha}(-\gamma(n\delta)^{\alpha})$ becomes $E_{1}(-\gamma n\delta)=e^{-\gamma n\delta}$,
  then $\mathbb{E}\left[\left|\tilde{X}_{\tau_{n}}\right|^{2}\right]\leq\mathbb{E}\left[\left|\tilde{X}_{\tau_{0}}\right|^{2}\right]e^{-\gamma n\delta}$ holds, which is the classical definition of mean-square exponential stability, confirming the consistency and generality of the theoretical framework.
  }
\end{remark}

\section{Numerical Simulation}
This section contains several examples illustrating the properties of the TCSDEs and the results we proved earlier. 

For a given step size $\delta$, the exact and numerical solutions of the ST method to \eqref{SDE} are simulated in the following way. 
\begin{itemize}
  \item [1)] One path of the stable subordinator $D_t$ is simulated with step sizes $\delta$. (See for example \cite{jum2015numerical}).
  \item [2)] The corresponding approximated inverse subordinator $\tilde{E_t}$ is found using \eqref{Et}. $\tilde{E_t}$ efficiently approximates $E_t$ by \eqref{approximates}.
  \item [3)] The exact solution is simulated directly by the closed-form solution; if the explicit form of the true solution is hard to obtain, use the numerical solution with step size $\delta_0=10^{-5}$ as the exact solution.
  \item [4)] The numerical solution of the ST method with the step size $\delta$ is computed via \eqref{continuous STM2}.
\end{itemize}

One path of a 0.9-stable subordinator $D_t$ is plotted using $\delta=10^{-4}$ and the corresponding approximated inverse subordinator $\tilde{E_t}$ are drawn in Figure 1. All the remaining simulations build upon a 0.9-stable subordinator $D_t$ as well.

\begin{figure}
  \centering
  \includegraphics[width=0.65\linewidth]{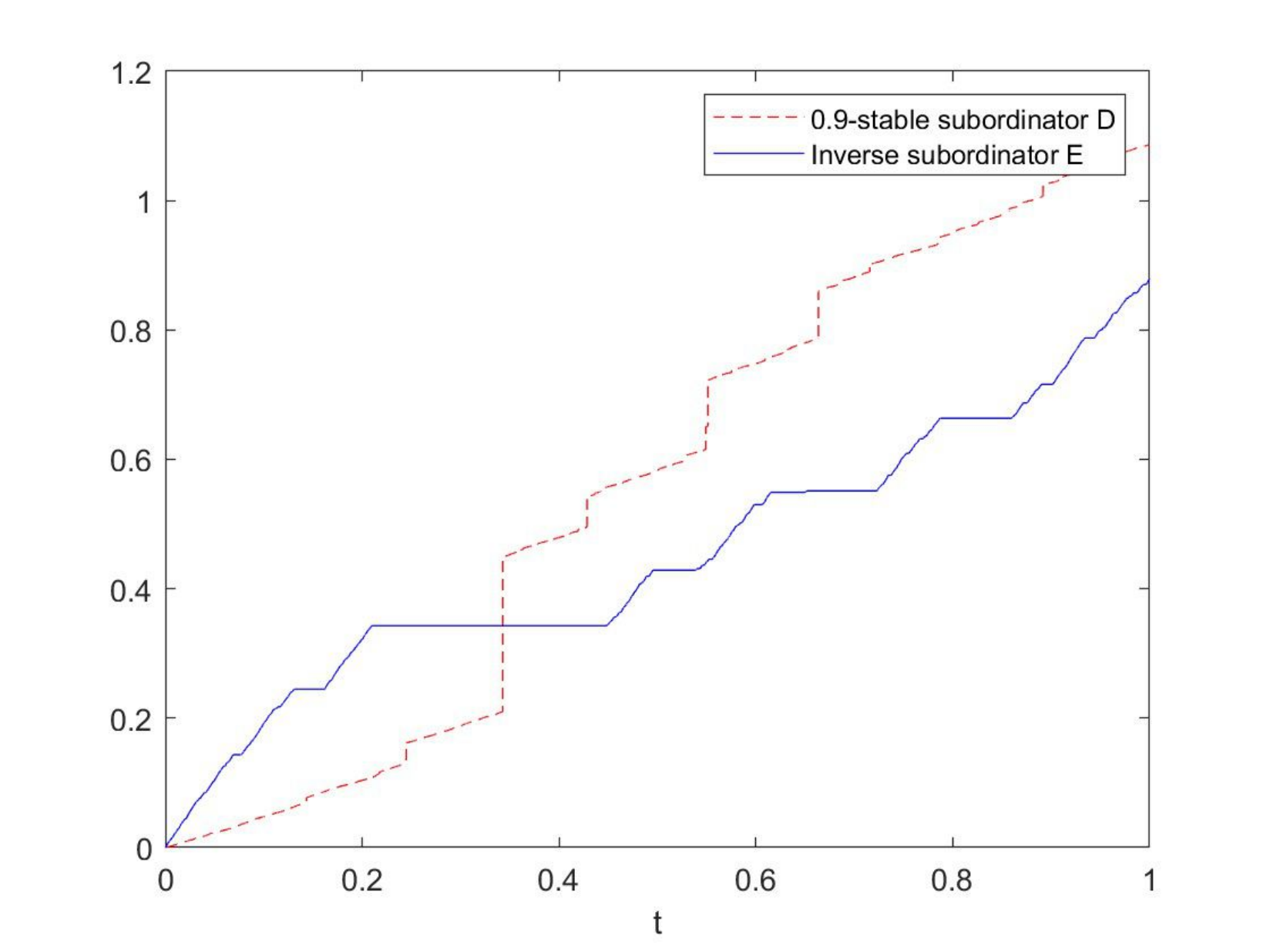}
  \caption{Sample paths of a 0.9-stable subordinator $D$ and the corresponding inverse $E$.}
\end{figure}

\begin{Example}(Time-changed Black-Scholes Model)\label{Example 5.1}
  We first consider the following basic linear TCSDE:
  \begin{equation}
    dX_{t}=\mu X_{t}dE_{t}+\sigma X_{t}dB_{E_{t}},X_{0}=x_{0}, \label{BS}
  \end{equation}
  where $\mu,\sigma$ are real constants and $x_{0}>0,\sigma>0.$
  It can be verified that the conditions in \autoref{thm: convergence} hold.
\end{Example}
This is also an analogue of the so-called Black-Scholes SDE. 
In finance, $\mu$ is the drift rate of annualized stock's return, and $\sigma$ is the standard deviation of the stock's returns.
This classical example has drift coefficient $F(x)=\mu x$ and diffusion coefficient $G(x)=\sigma x$ both satisfying global Lipschitz and linear growth condition, and surely
further satisfying the weaker conditions of \autoref{thm: convergence}.
By Proposition 4.4 of \cite{2011Kobayashi}, the exact solution of \eqref{BS} is given by
\begin{equation}
    X_{t}=x_{0}exp\left\{ \left(\mu-\frac{1}{2}\sigma^{2}\right)E_{t}+\sigma B_{E_{t}}\right\}.  \notag
\end{equation}
Also, the solution of \eqref{BS} has the following asymptotic behavior:
\begin{itemize}
    \item [1)] If $\mu>\sigma^{2}/2$, then $\underset{t\rightarrow\infty}{lim}X_{t}=\infty \ a.s.$
    \item [2)] If $\mu<\sigma^{2}/2$, then $\underset{t\rightarrow\infty}{lim}X_{t}=0\  a.s.$
    \item [3)] If $\mu=\sigma^{2}/2$, then $X_{t}$ asymptotically fluctuates between arbitrarily large and arbitrarity small positive values infinitely often. 
\end{itemize}
%Fig 2
Letting $x_0=1$, consider three different combinations of parameters $\mu$ and $\sigma$, we can see from Figure 2 that 
when $\mu=0.02,\sigma=0.2$, the solution fluctuates between 0.96 and 1.16 infinitely; when $\mu=0.2,\sigma=0.2$, the solution has an upward trend; when $\mu=0.002,\sigma=0.2$, the solution has a doward trend, which verify the statement above.

\begin{figure}
    \centering
    \includegraphics[width=0.7\linewidth]{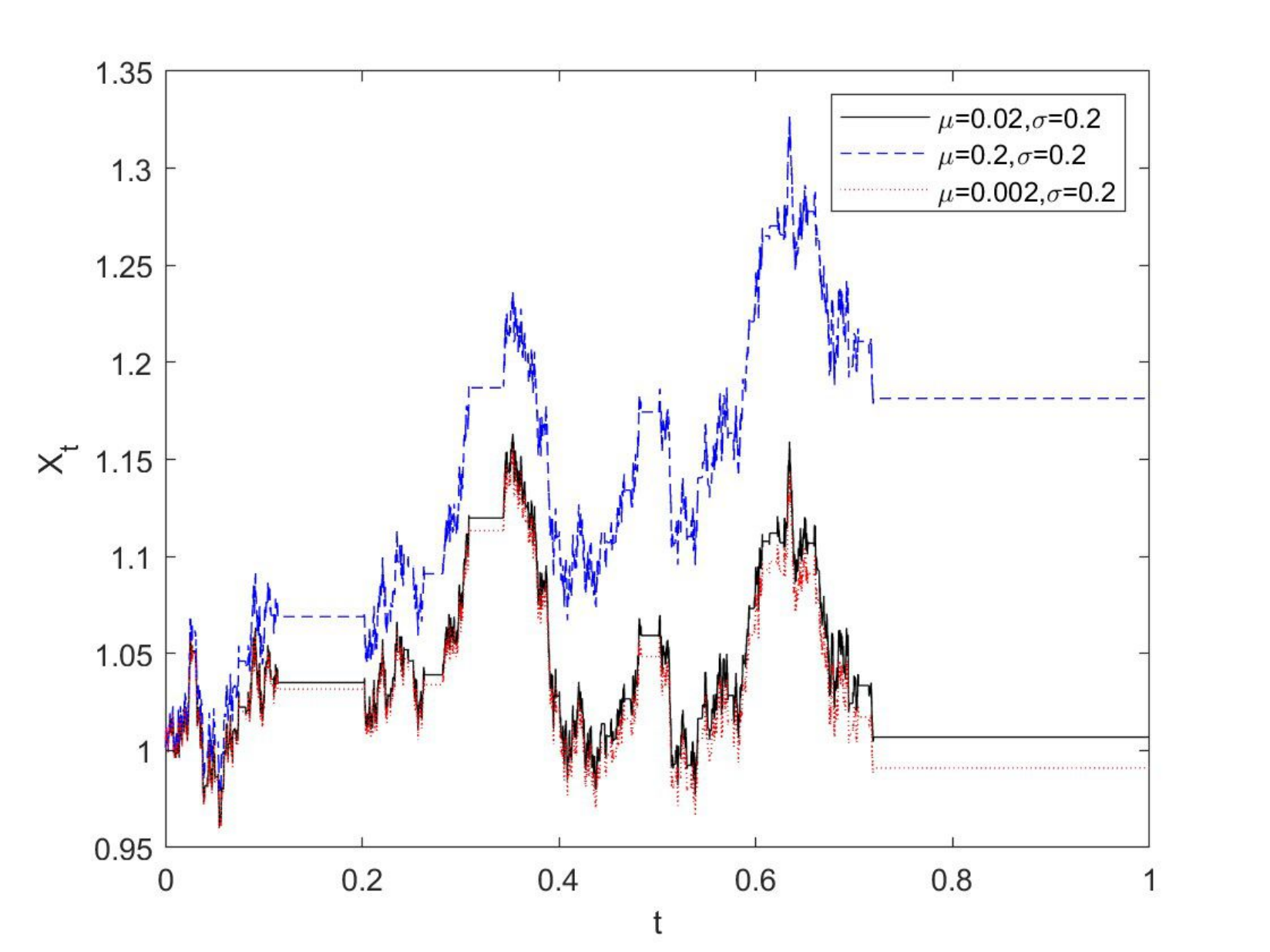}
    \caption{Sample paths of time-changed Black-Scholes SDE with different coefficients.}
\end{figure}

\begin{Example}\label{Example 5.2}
  We consider the following nonlinear TCSDE
  \begin{equation}
    dX_{t}=\left(-X_{t}-\frac{X_{t}^{3}}{1+X_{t}^{2}}\right)dE_{t}+\frac{X_{t}}{\sqrt{1+X_{t}^{2}}}dB_{E_{t}} \notag
  \end{equation}
  with initial data $X_{0}=1$.
\end{Example}
This example is a specific case of a broader class of TCSDEs whose coefficients satisfy global Lipschitz and super-linear growth conditions.
The idea of construction for \autoref{Example 5.2} is that, in drift coefficient, when $|X_{t}|$ is large, $-\frac{X_{t}^{3}}{1+X_{t}^{2}}$ is close to $-X_{t}$, which prevents solution from blowing up; and the boundedness of diffusion coefficient is ensured as when $|X_{t}|\rightarrow\infty,\frac{X_{t}}{\sqrt{1+X_{t}^{2}}}\rightarrow1$. Thus the coefficients are bounded globally.
It is easy to verify that the drift and diffusion coefficients in \autoref{Example 5.2} satisfy Assumptions 3.1-3.4 and 4.1. For \autoref{monotone}, we have
\begin{equation}
  \langle x,F(t,x)\rangle=x\cdot\left(-x-\frac{x^{3}}{1+x^{2}}\right)=-x^{2}\left(1+\frac{x^{2}}{1+x^{2}}\right),\forall x,
\end{equation}
and 
\begin{equation}
  |G(t,x)|^{2}=\left(\frac{x}{\sqrt{1+x^{2}}}\right)^{2}
  =\frac{x^{2}}{1+x^{2}},\forall x.
\end{equation}
Choose $h=1$. 
Then there exists a constant $K_1=1$ such that
\begin{equation}
  \langle x,F(t,x)\rangle+\frac{2h-1}{2}|G(t,x)|^{2}
  =-x^{2}\left(1+\frac{x^{2}}{1+x^{2}}\right)+\frac{1}{2}\frac{x^{2}}{1+x^{2}}
  \leq K_{1}(1+|x|^{2}), \forall x.
\end{equation}
Five sample paths of exact solutions are approximated using the ST method in Figure 3.
\begin{figure}
  \centering
  \includegraphics[width=0.75\linewidth]{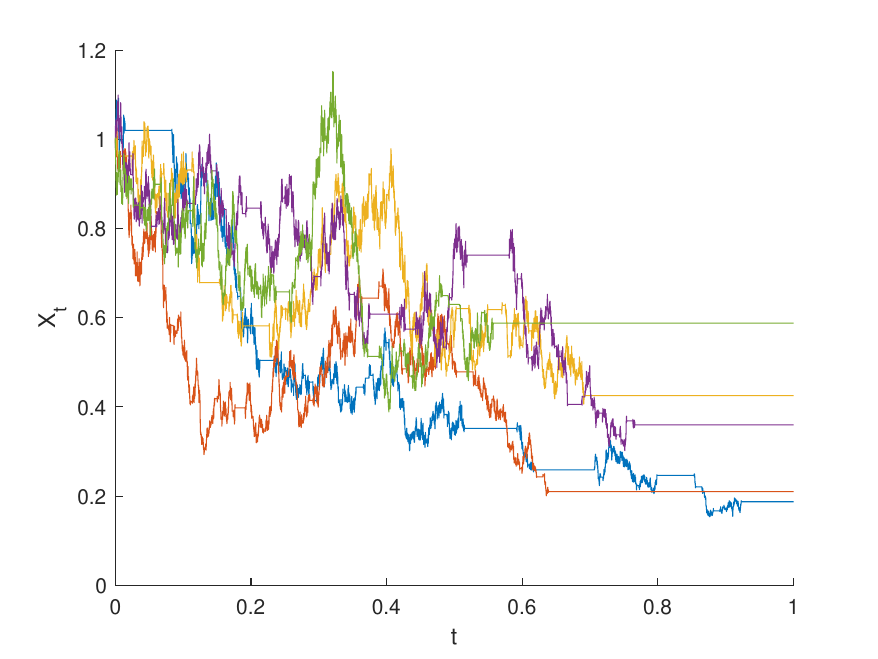}
  \caption{Five sample paths of \autoref{Example 5.2}.}
\end{figure}

\begin{Example}(Time-changed Mean-Reverting Process with state-time-dependent coefficient)\label{Example 7.3}
  We consider the following time-changed mean-reverting process with state-time-dependent coefficient
  \begin{equation}
    dX_{t}=\kappa(\theta_{t}-X_{t})dE_{t}+\sigma_{t}X_t^3 dB_{E_{t}}.
  \end{equation}
  where 
    $\kappa >0$ is a constant representing the speed of mean reversion;
    $|\theta_{t}|\leq M_{\theta}$ is a bounded time-dependent function, representing the mean reversion level. For instance, $\theta_{t}=\theta_{0}+Asin(\omega t)$ models seasonal effects or business cycles; and
    $|\sigma_{t}|\leq M_{\sigma}$ is a bounded time-dependent function representing volatility. 
  \end{Example}
 
Here the time dependence arises from the external physical time $t$ is not synchoronized with the internal stochastic time $E_t$. This makes it impossible to simplify the equation using techniques like the duality principle. 
It is easy to check that the functions $\theta_t$ and $\sigma_t$ can be chosen to satisfy Assumptions 3.1-3.4 and 4.1. 
The drift coefficient $F(t,x)=\kappa(\theta_{t}-X_{t})$ satisfy the global Lipschitz condition, but $G(t,x)=\sigma_t X_t^3$ does not, since its derivative $3\sigma _t X_t^2$ grows without bound as $x\rightarrow \infty$, making it impossible to constrain its growth with a single, finite Lipschitz constant across the entire real line.
For Assumption 3.3, choose $h=1$.
Since $\kappa>0$, and $\theta_{t},\sigma_{t}$ are bounded, then there exists constant $K_{1}=max\left\{ \frac{\kappa}{2}M_{\theta}^{2}+\frac{1}{2}M_{\sigma}^{2},\frac{|\kappa|}{2}\right\}$, such that 
\begin{equation}
  \langle x,F(t,x)\rangle+\frac{2h-1}{2}|G(t,x)|^{2}
  =-\kappa x^{2}+\kappa\theta_{t}x+\frac{1}{2}\sigma_{t}^{2}x^{6}
  \leq K_{1}(1+|x|)^{2}. \notag
\end{equation}

%Fig 3-distribution of sample MSEs
Now we illustrate the strong convergence result.
As per common practice, for a given step size $\delta$, by the law of large numbers, at the final time $T=1$, we focus on the strong mean square error 
\begin{equation}
  \frac{1}{N}\sum_{i=1}^{N}\left|X(E_{T})(\omega_{i})-X(\tilde{E}_{T})(\omega_{i})\right|^{2}
\end{equation}
as a reference baseline.

%Fig 4
We consider \autoref{Example 5.1} with $x_0=1$, $\mu=0.02$ and $\sigma=0.2$.
Three thousands ($N=3000$) sample paths are used to calculate the strong mean square error at time $T$. 
Fix $\theta=1, \delta=10^{-4}$, Figure 4 shows the distribution of errors at time $T=1$ for each of the 3000 sample paths. We can see that the ST method is efficient since the majority of errors lie within the range below 0.001, and the mean square error is 0.0003.
\begin{figure}
  \centering
  \includegraphics[width=0.75\linewidth]{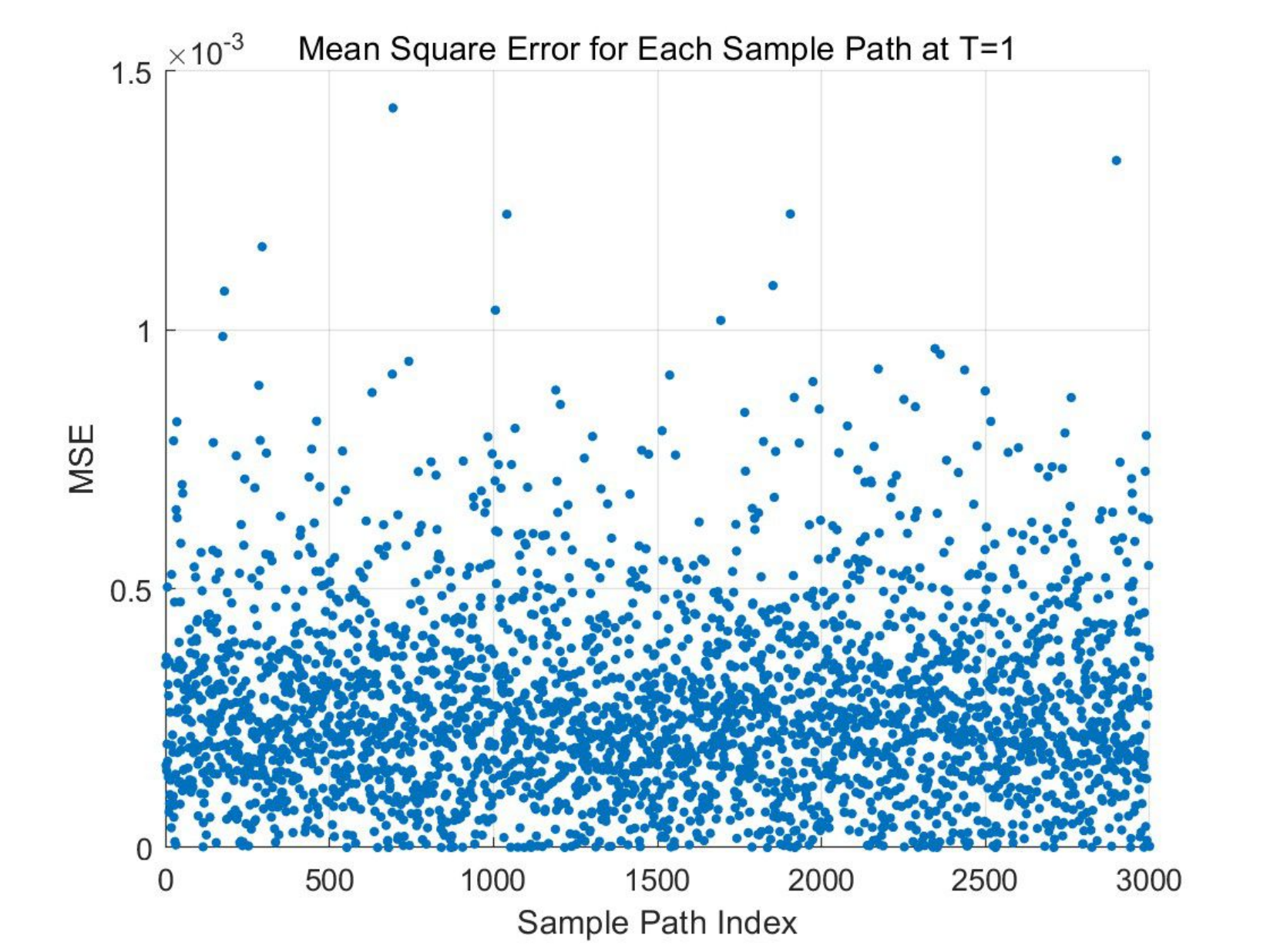}
  \caption{Distribution of sample mean square errors for \autoref{Example 5.1}.} 
\end{figure}

%Fig 5
\begin{figure}
  \centering
  \begin{subfigure}{0.49\linewidth}
    \includegraphics[width=\linewidth]{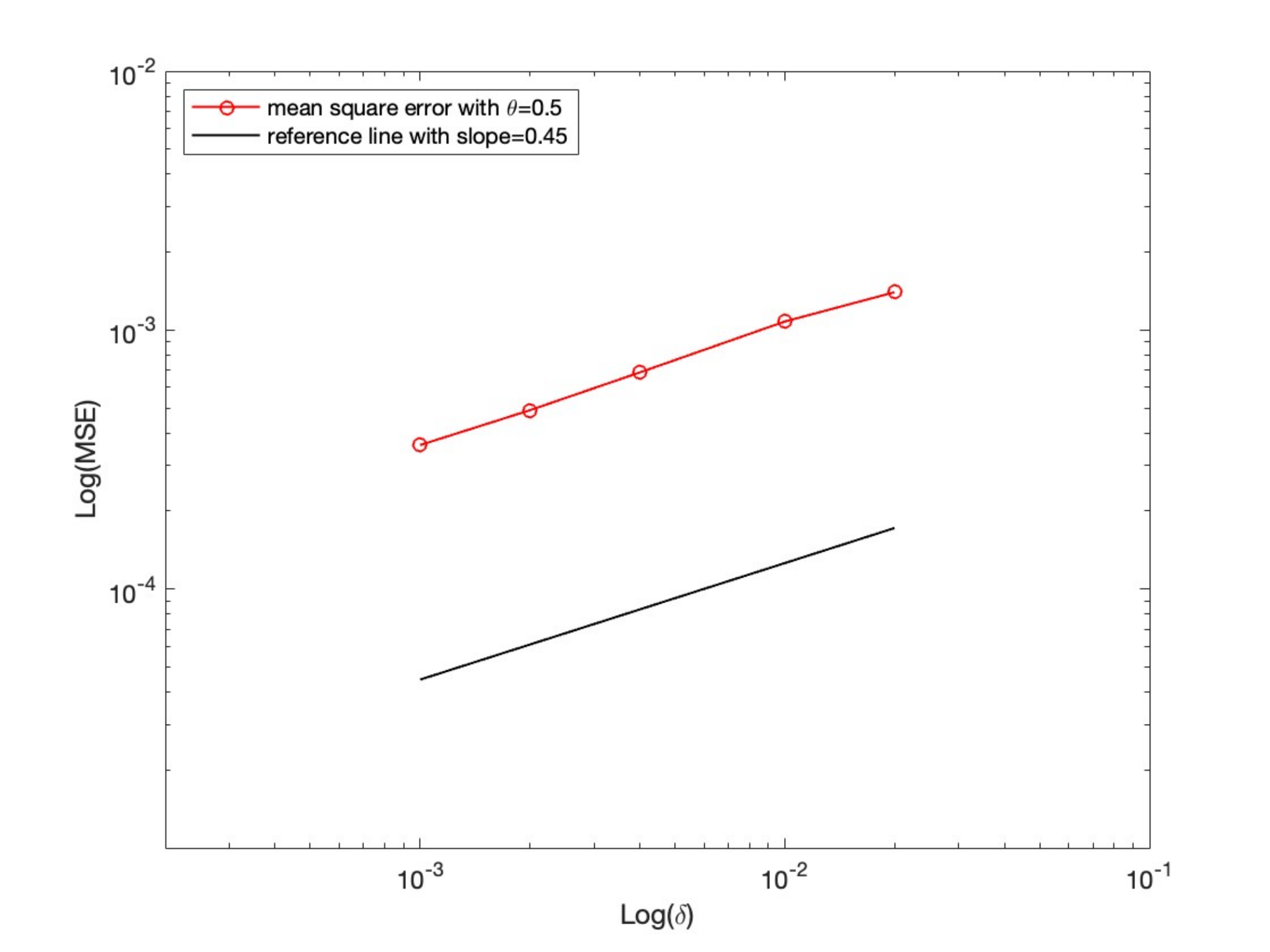}
    \subcaption{$\theta = 0.5$}
  \end{subfigure}
  \begin{subfigure}{0.49\linewidth}
    \includegraphics[width=\linewidth]{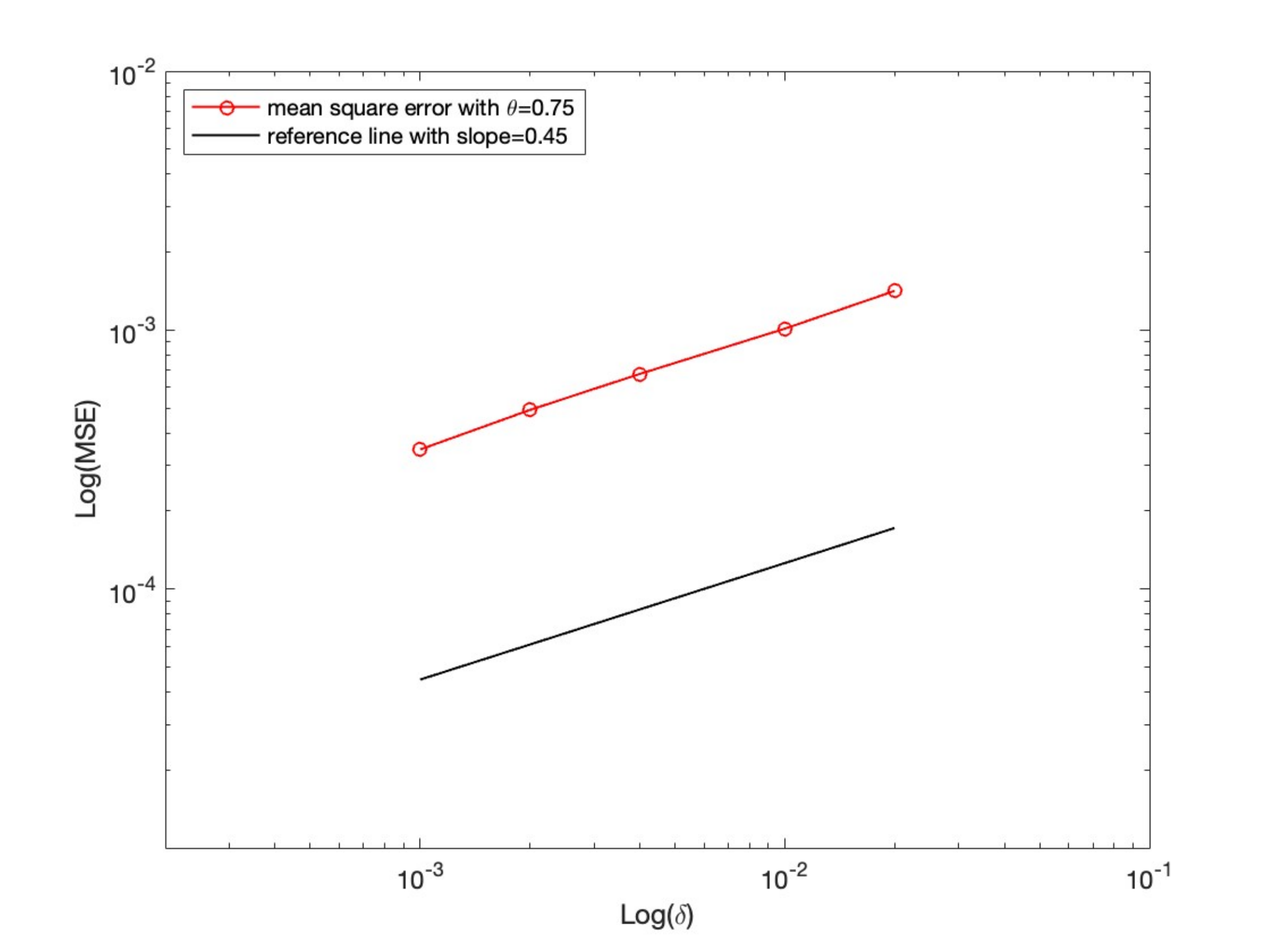} 
    \subcaption{$\theta = 0.75$}
  \end{subfigure}
  %\vspace{0.5cm} % 垂直间距调整
  \begin{subfigure}{0.49\linewidth}
    \includegraphics[width=\linewidth]{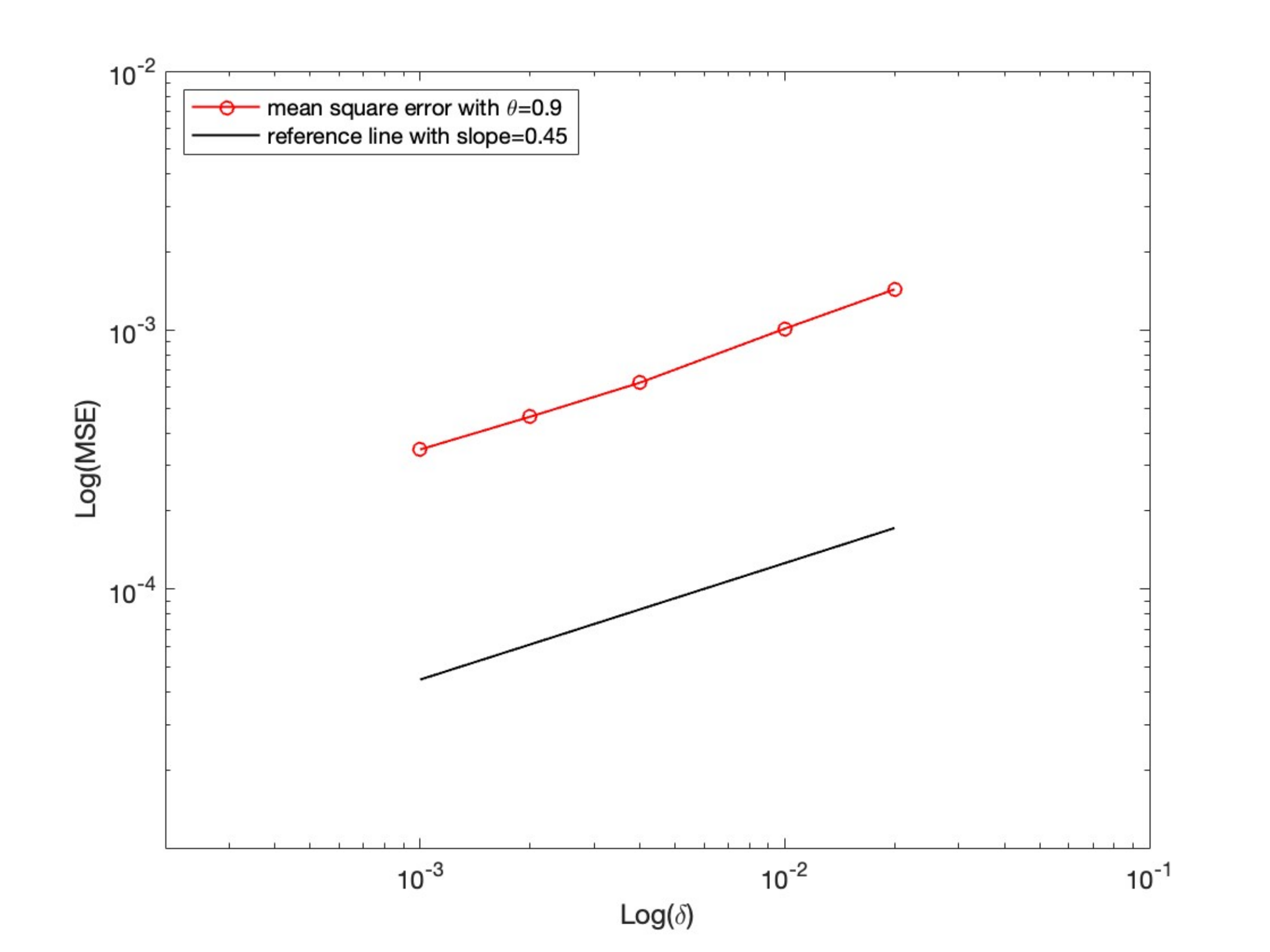} 
    \subcaption{$\theta = 0.9$}
  \end{subfigure}
  \begin{subfigure}{0.49\linewidth}
    \includegraphics[width=\linewidth]{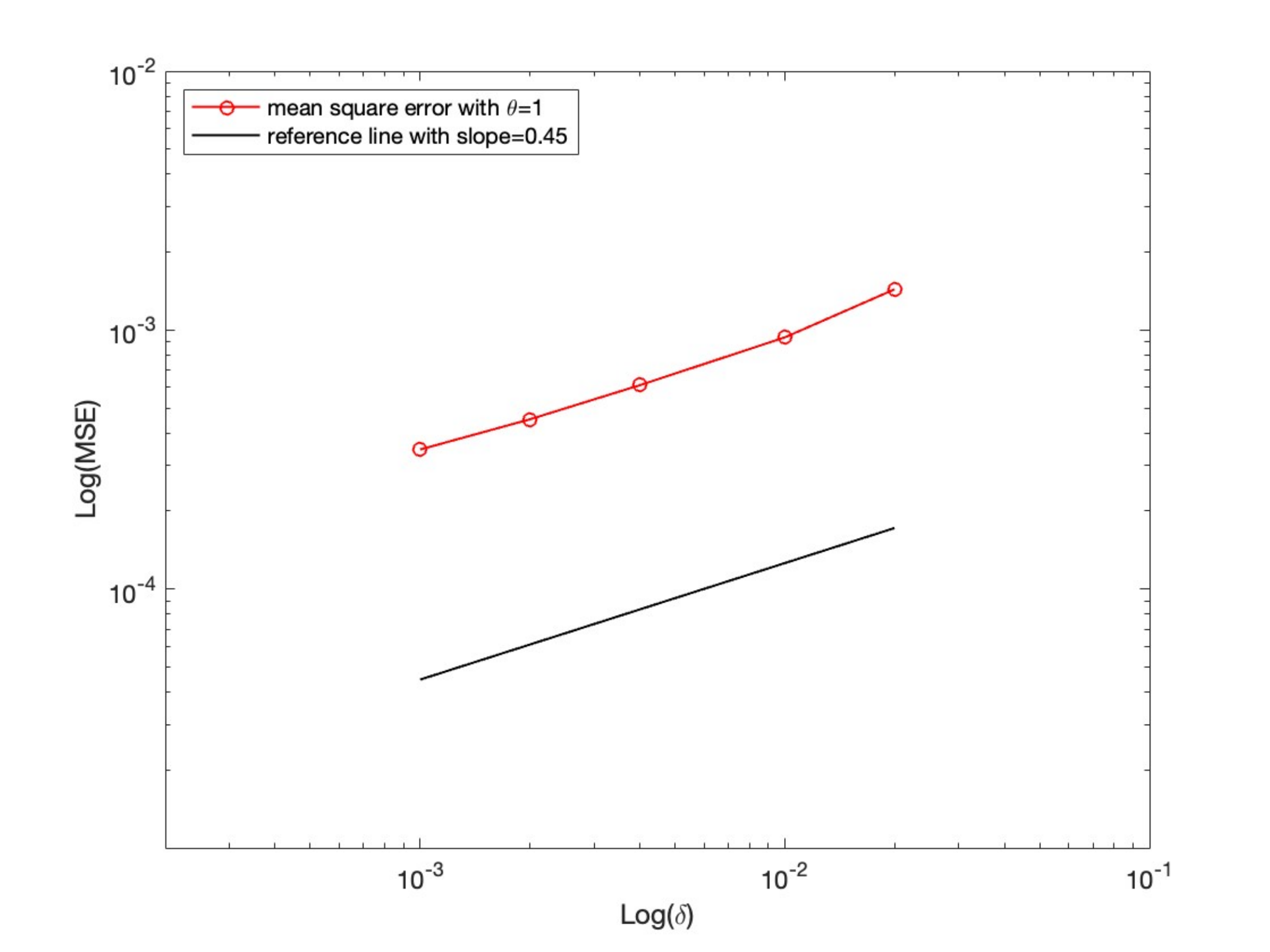} 
    \subcaption{$\theta = 1$}
  \end{subfigure}
  \caption{Convergence order simulation for \autoref{Example 7.3} with $\theta_{t}=0.05+0.03sin(2\pi t)$,
  $\sigma_{t}=0.4\times(1+0.05t)$,
  $\kappa=0.65$}
\end{figure}

%%Fig 4-convergence rate
We show the overall mean square errors across all sample paths for \autoref{Example 7.3} by drawing the Loglog plot of the mean square errors against the step sizes in Figure 5.
Since we choose $\alpha=0.9$, the convergence order suggested by \autoref{thm: convergence} is the minimum value of $\alpha/2=0.45, \eta_F$, and $\eta_G$.
%Example 7.3
We consider 
$\theta_{t}=0.05+0.03sin(2\pi t)$, 
$\sigma_t = 0.4\times(1 + 0.05t)$, and 
$\kappa=0.65$.
The numerical solution of ST method is computed using five different step sizes $\delta =2\times 10^{-2},10^{-2},2\times 10^{-3}, 4\times 10^{-3}$ and $10^{-3}$ for $\theta=0.5, 0.75, 0.9$ and $1$ respectively.
We can see that the ST method is convergent in Figure 5, which verify the theory of \autoref{thm: convergence}.
Graphically, the convergence order is in accordance with the reference convergence order between step sizes $\delta =2\times 10^{-2}$ and $10^{-3}$. 

Next we want to see how convergence order changes as $\alpha$ changes. 
If we choose $\alpha=0.55$, the convergence order suggested by \autoref{thm: convergence} is the minimum value of $\alpha/2=0.275, \eta_F$, and $\eta_G$.
The simulations prove our thereotical findings: the convergence order increases as $\alpha$ increases. 
The approximate convergence orders are 0.448 and 0.18 for \autoref{Example 7.3} in Figure 6(a) and Figure 6(b), respectively.

\begin{figure}
  \centering
  \begin{subfigure}{0.49\linewidth}
    \includegraphics[width=\linewidth]{figs/ex7.3_theta0.9.pdf}
    \subcaption{$\alpha=0.9$}
  \end{subfigure}
  \begin{subfigure}{0.49\linewidth}
    \includegraphics[width=\linewidth]{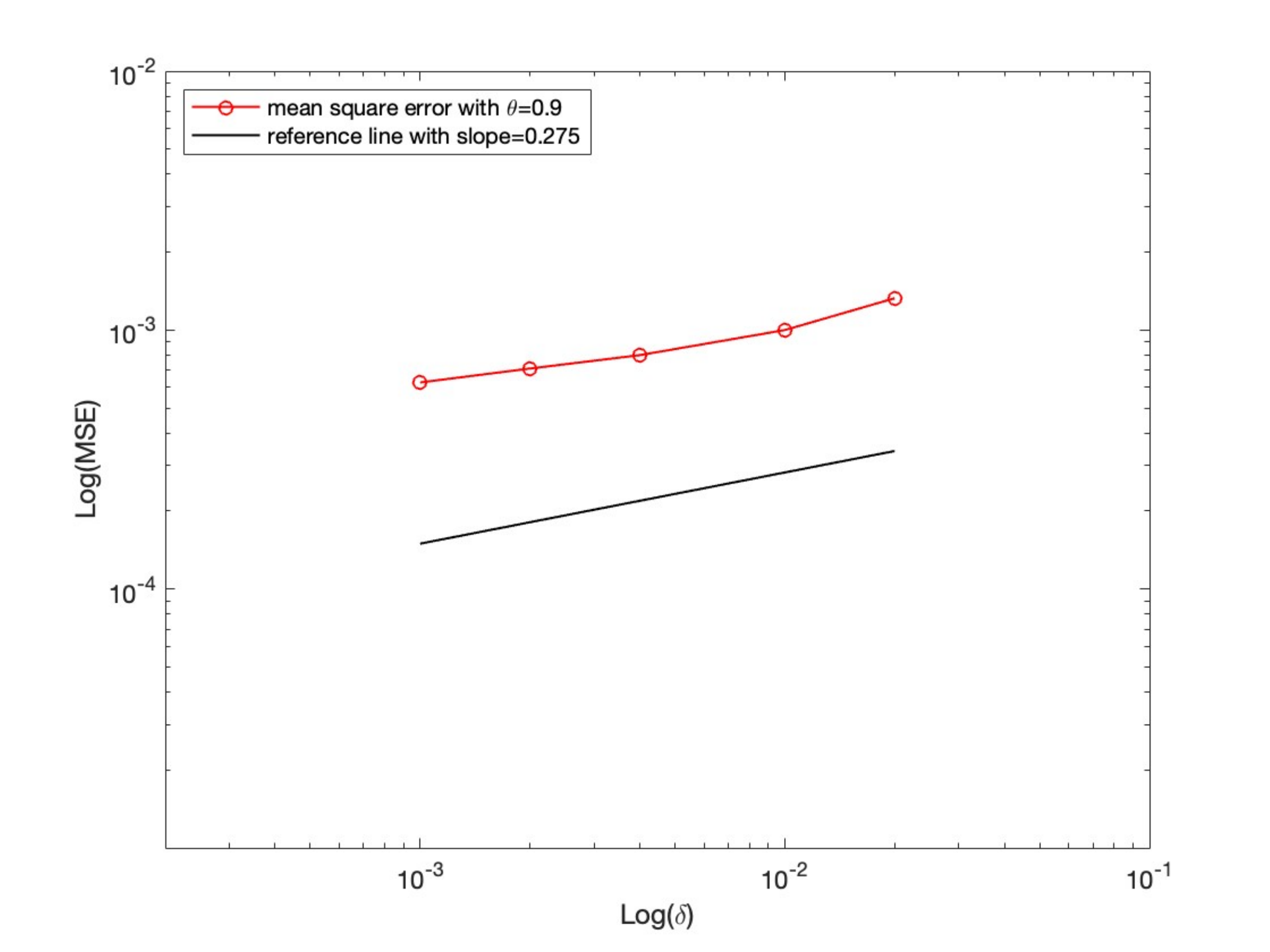} 
    \subcaption{$\alpha=0.55$}
  \end{subfigure}
  \caption{Convergence order vs. $\alpha$  for \autoref{Example 7.3} with $\theta_{t}=0.05+0.03sin(2\pi t)$,
  $\sigma_{t}=0.4\times(1+0.05t)$,
  $\kappa=0.65$}
\end{figure}

%%稳定性分析

We finally verify the Mittag-Leffler stability.
%%Example 7.5
\begin{Example}\label{Example 7.4}
  We consider the following linear TCSDE
  \begin{equation}
    dX_{t}=-2X_{t}dE_{t}+X_{t}dB_{E_{t}}, \notag
  \end{equation}
  with initial data $X_{0}=1$.
\end{Example}
It is easy to check that the drift and diffusion coefficients in \autoref{Example 7.5} satisfy \autoref{coercive} and condition \eqref{K5}, that is, 
\begin{equation}
  \left\langle x,F(t,x)\right\rangle +\frac{1}{2}|G(t,x)|^{2}
  =-2x^{2}-\frac{1}{2}x^{2}
  \leq -\lambda x^2, \forall x. \notag
\end{equation}
and 
\begin{equation}
  \left|F(t,x)\right|^{2}=\left(-2x\right)^{2}\leq K|x|^2, \forall x. \notag
\end{equation}
where $\lambda = 2+1/2$, and $K_5=4$.
This example clearly satisfies the global Lipschitz condition.

In Figure 7, for step size $\delta=2,1,1/2$, we simulate the mean-square stability of the ST method on the \eqref{SDE} with $\theta=0,0.25,0.5$ and $1$, respectively. Three thousand sample paths are generated.
The figure provides a compelling visual argument for the significant benefits of using the ST method in the stability analysis of TCSDEs. The figure systematically demonstrates how the parameter $\theta$, which controls the implicitness of the method, is a critical tool for ensuring and analyzing the numerical stability. We discuss the tunable stability property of the ST method in the following four cases.
\begin{itemize}
  \item For $\theta=0$ (the purely explicit method), the results are highly sensitive to step sizes. The mean-square norm of the numerical solutions fail to decay and instead grows without bound. This is a classic sign of numerical instability: making long-time stability analysis computationally expensive or impossible.
  \item For $\theta=0.25$, the stability improves dramatically. The numerical solution still blows up at $\delta=2$, but is stable in mean square at $\delta=1, 1/2$; 
  \item For $\theta= 0.5$ (the trapezoidal method), the curves for all three step sizes $\delta=2,1,1/2$ now decay to zero. This shows that the method has become stable for the given step size, correctly capturing the system's dissipative nature.
  \item For $\theta$ = 1 (fully implicit method, also known as the Backward Euler-Maruyama (BEM) method), the decay curves are consistently lower and decay more rapidly than those for the $\theta = 0.5$, showing us that methods with higher implicitness is essential for ensuring stability when simulating such complex systems over long time horizons. 
\end{itemize}
 In addition, TCSDEs are driven by time-changed processes like $B_{E_{t}}$, where the inverse subordinator $E_t$ has irregular, Hölder-continuous paths. This "roughness" exacerbates the stability challenges for numerical methods. The ST method, particularly its implicit versions ($\theta\geq 0.5$), is exceptionally well-suited for this, i.e., the implicit treatment of the drift term in the ST method provides a stabilizing force that helps to dampen the numerical errors amplified by the rough paths of the driving process.

\begin{figure}
  \centering
  \begin{subfigure}{0.49\linewidth}
    \includegraphics[width=\linewidth]{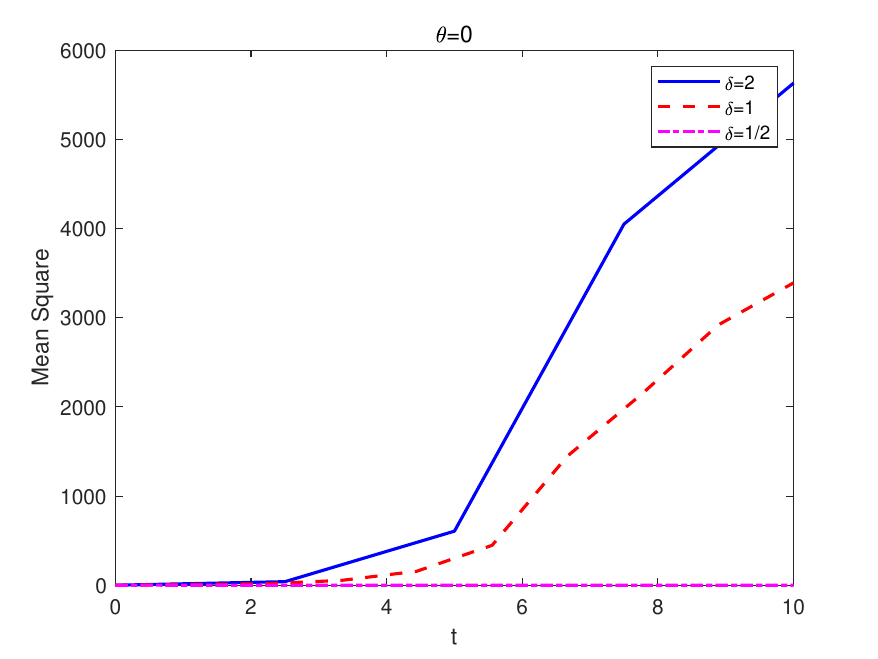}
  \end{subfigure}
  \hfill
  \begin{subfigure}{0.49\linewidth}
    \includegraphics[width=\linewidth]{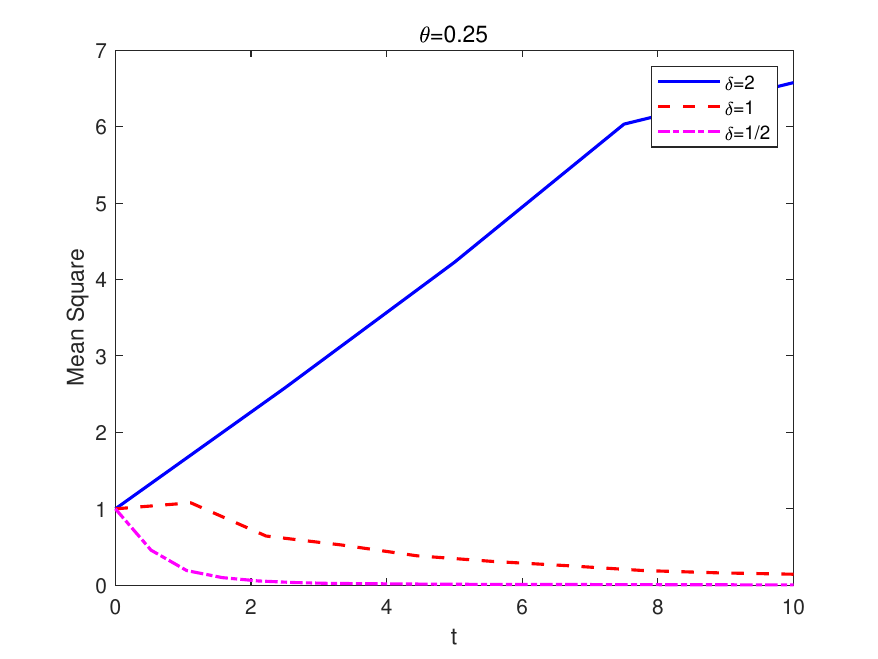}
  \end{subfigure}
  \vspace{0.5cm}
  \begin{subfigure}{0.49\linewidth}
    \includegraphics[width=\linewidth]{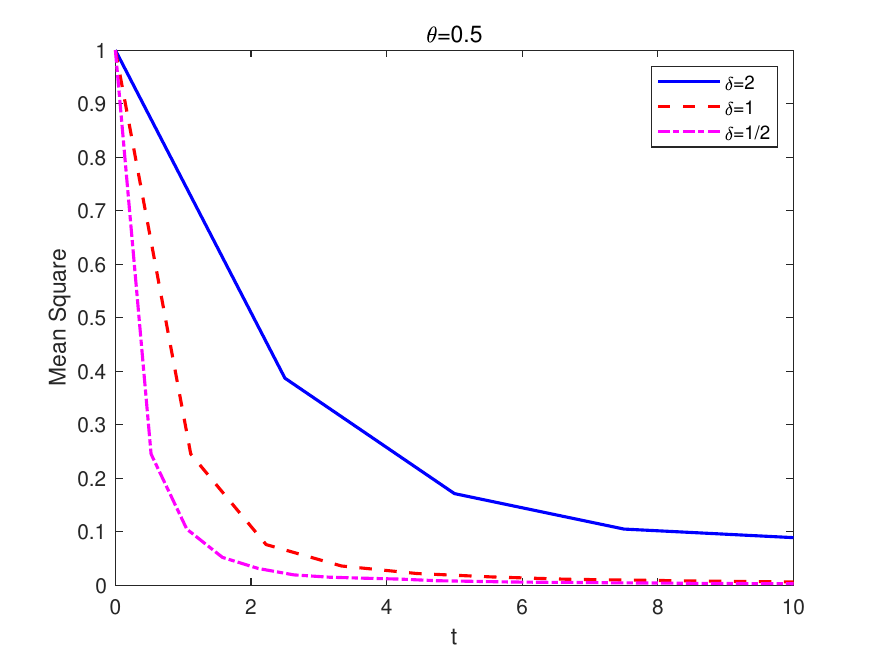}
  \end{subfigure}
  \hfill
  \begin{subfigure}{0.49\linewidth}
    \includegraphics[width=\linewidth]{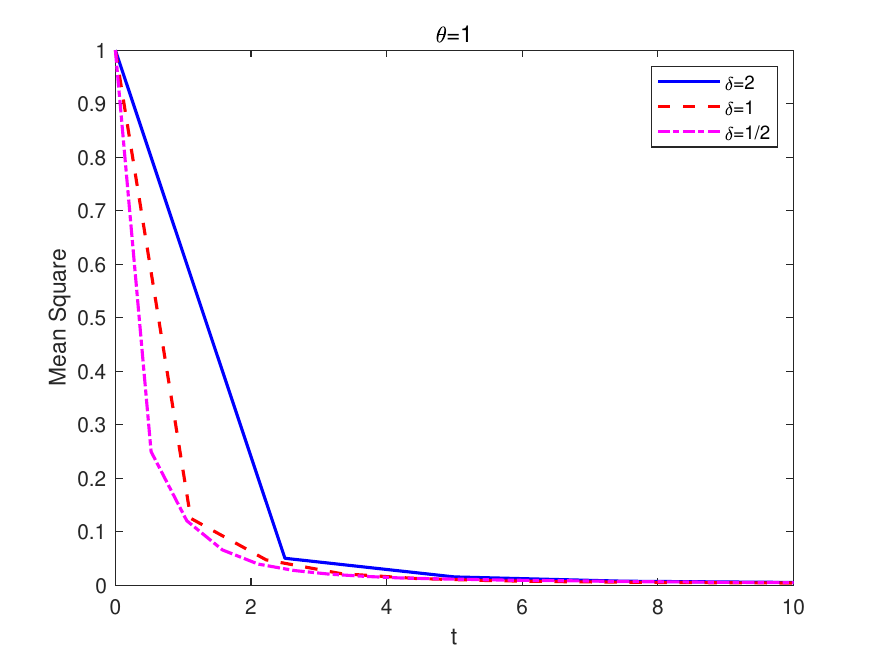}
  \end{subfigure}
  \caption{Mean-square Mittag-Leffler stability of the numerical solutions for \autoref{Example 7.4} with different $\theta$ values.}
\end{figure}

\begin{Example}\label{Example 7.5}
  We consider the following nonlinear TCSDE
  \begin{equation}
    dX_{t}=\left(-2X_{t}-X_{t}^{3}\right)dE_{t}+X_{t}dB_{E_{t}}, \notag
  \end{equation}
  with initial data $X_{0}=1$.
\end{Example}
\begin{Example}\label{Example 7.6}
  We consider the following nonlinear TCSDE
  \begin{equation}
    dX_{t}=\left(-X_{t}-X_{t}^{3}\right)dE_{t}+X_{t}^{2}dB_{E_{t}}, \notag
  \end{equation}
  with initial data $X_{0}=1$.
\end{Example}
\begin{Example}\label{Example 7.7}
  We consider the following nonlinear TCSDE
  \begin{equation}
    dX_{t}=\left((-2t-1)X_{t}-X_{t}^{3}\right)dE_{t}+X_{t}dB_{E_{t}}, \notag
  \end{equation}
  with initial data $X_{0}=1$.
\end{Example}

%Example 7.6 is not GL, is LL
In \autoref{Example 7.5},
the derivative of the drift coefficient, given by $|F'(x)|=|-2-3x^{2}|$, exhibits unbounded growth as $|x|\rightarrow\infty.$
This property precludes the existence of a finite constant $C$ that would satisfy a global Lipschitz condition over the entire real line. 
Nevertheless, when restricted to any bounded domain $[-R, R]$, the derivative remains bounded with $|F'(x)|\leq 2+3R^{2}$. 
Consequently, the function $F(t, x)$ satisfies the local Lipschitz condition, as it is Lipschitz continuous on every compact subset of its domain.
\autoref{Example 7.6} and \autoref{Example 7.7} both satisfy the local Lipschitz conditions as well.

All three cases satisfy \autoref{coercive}. For example, 
in \autoref{Example 7.7},
\begin{equation}
  \left\langle x,F(t,x)\right\rangle +\frac{1}{2}|G(t,x)|^{2}
  =x\cdot[(-2t-1)x-x^{3}]+\frac{1}{2}|x|^{2}
  \leq-\left(-2t-\frac{1}{2}\right)x^{2}-x^{4}
  \leq-\frac{1}{2}x^{2}, \notag
\end{equation}
here $\lambda=1/2$.

\autoref{Example 7.5} and \autoref{Example 7.7} both have a linear diffusion coefficient, while \autoref{Example 7.6} has a highly nonlinear diffusion coefficient. 
\autoref{Example 7.7} has a time-dependent drift coefficient.

Next we discuss how the duality principle fails.
For \autoref{Example 7.7}, the coefficient $(-2t-1)X_t$ is an explicit function of external time $t$, making the dual SDE intractable.
For \autoref{Example 7.5} and \autoref{Example 7.6}, although duality can be applied, 
the stability we concern is of Mittag-Leffler type, and duality principle is not suitable here (see Remark 5.1(3)).

Figure 8 presents numerical simulations of the mean-square stability for the ST method applied to Examples 7.5-7.7, with parameter values $\theta=0.5$ and $\theta=1$, respectively. 
The simulations employ step sizes of $\delta=1,1/2,1/4$, with statistical results obtained from an ensemble of 3,000 sample paths. A key observation is that for $\theta\in[1/2,1]$, the ST method effectively preserves mean-square stability even in the presence of highly nonlinear coefficients. These numerical findings are in full agreement with the theoretical predictions of \autoref{numerical stability}.

All decay curves demonstrate an initial phase of rapid decline followed by an asymptotic approach toward zero. This dynamic is consistent with the decay profile of the Mittag-Leffler function, which governs the rate of decay in the mean-square solution norm. Notably, this Mittag-Leffler-type decay is asymptotically slower than classical exponential decay, the latter corresponding to the special case of the stability index $\alpha=1$.

\begin{figure}
  \centering
  \begin{subfigure}[b]{0.49\linewidth}
    \includegraphics[width=\linewidth]{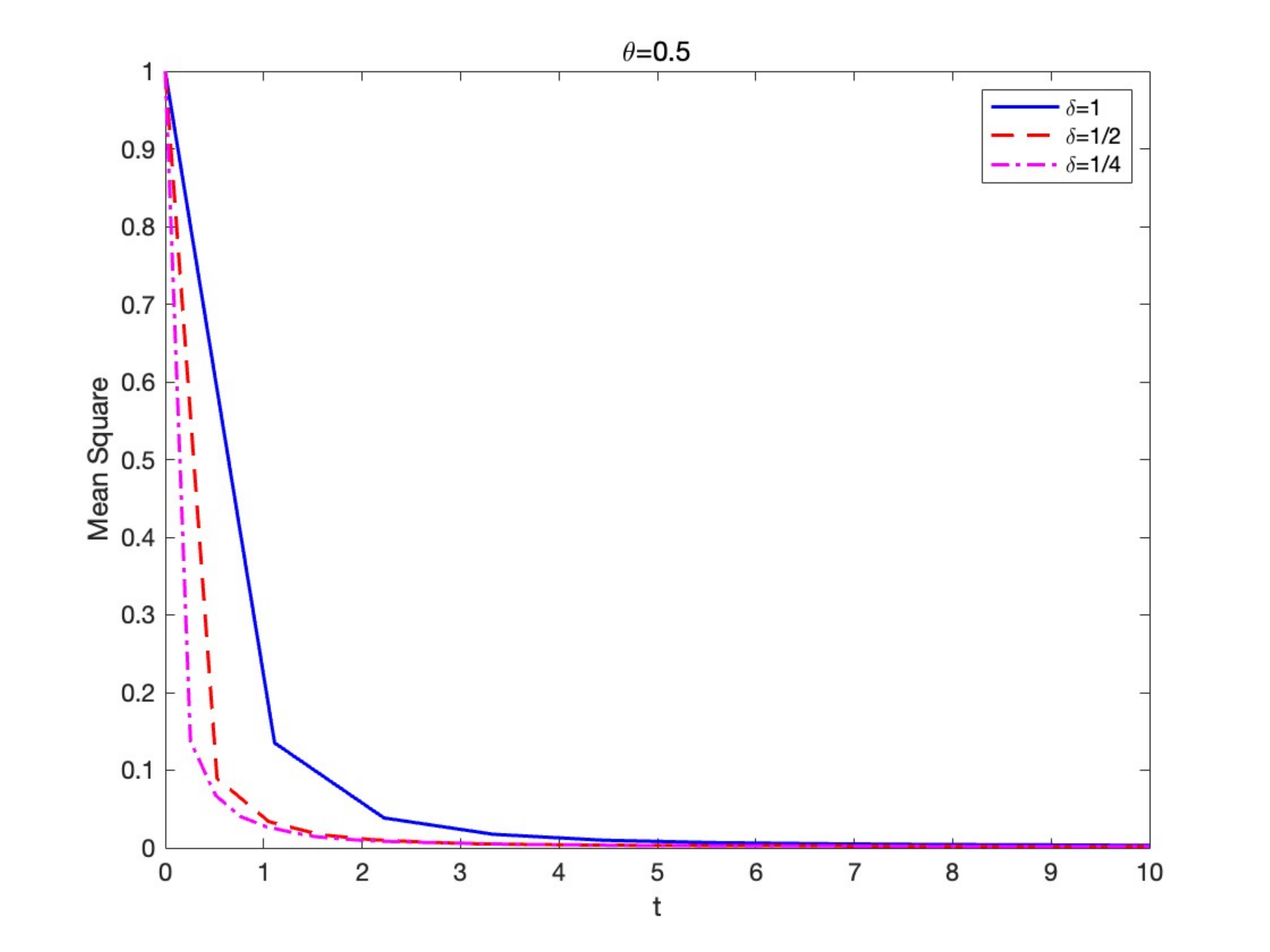}
  \end{subfigure}
  \hfill
    \begin{subfigure}[b]{0.49\linewidth}
    \includegraphics[width=\linewidth]{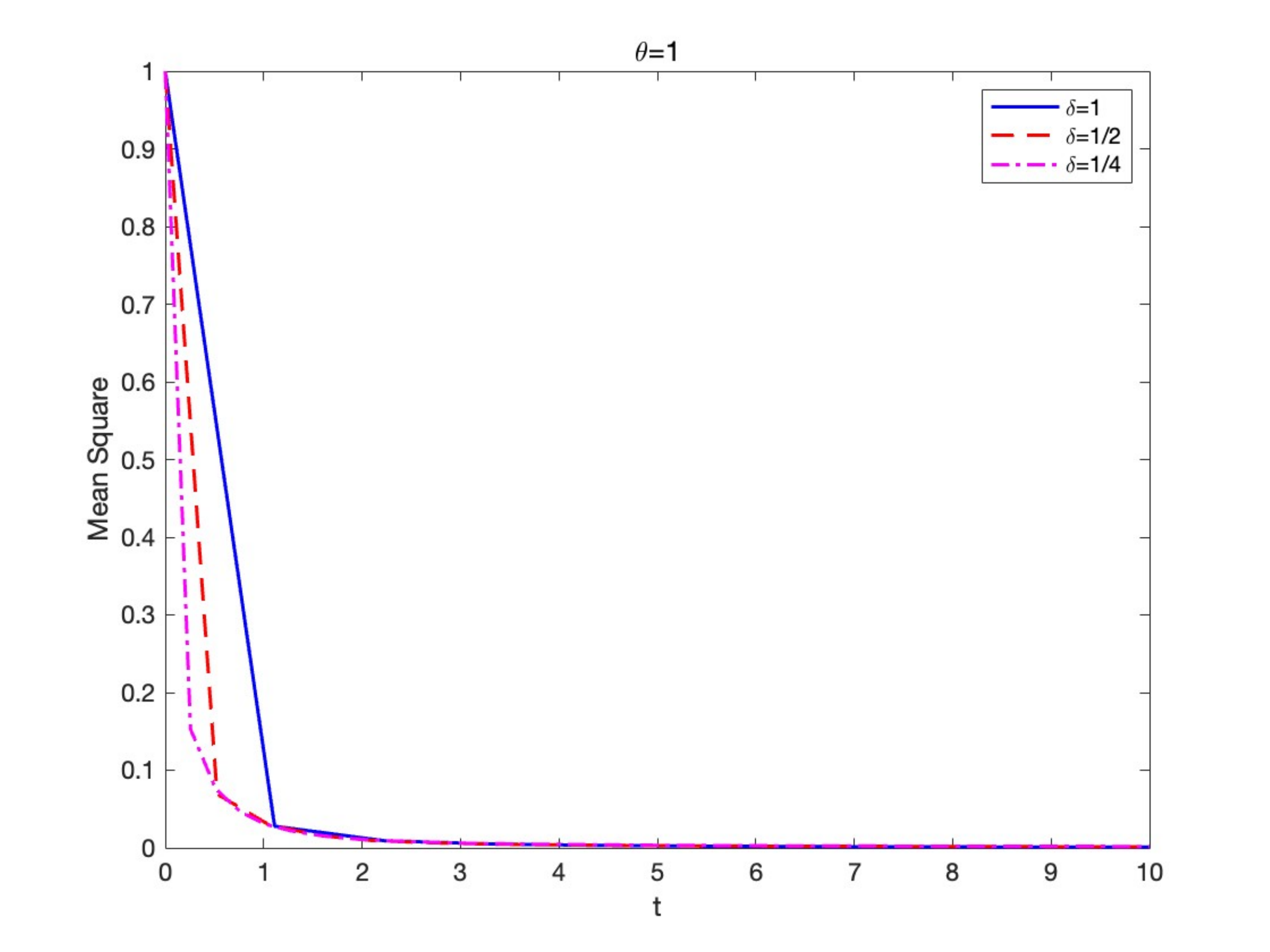}
  \end{subfigure}
    %\vspace{0.5cm}
  \begin{subfigure}[b]{0.49\linewidth}
    \includegraphics[width=\linewidth]{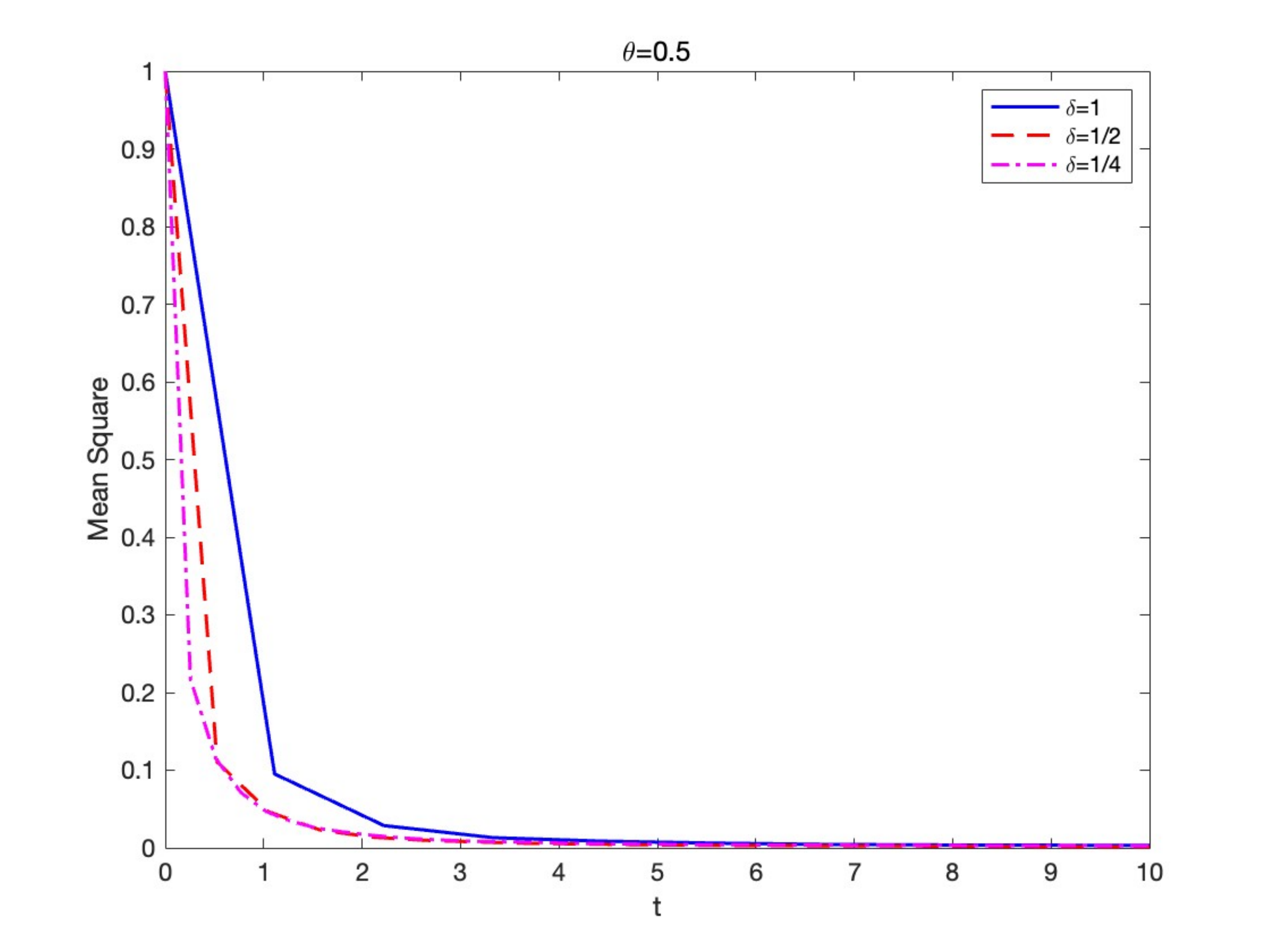}
  \end{subfigure}
  \hfill
  \begin{subfigure}[b]{0.49\linewidth} 
    \includegraphics[width=\linewidth]{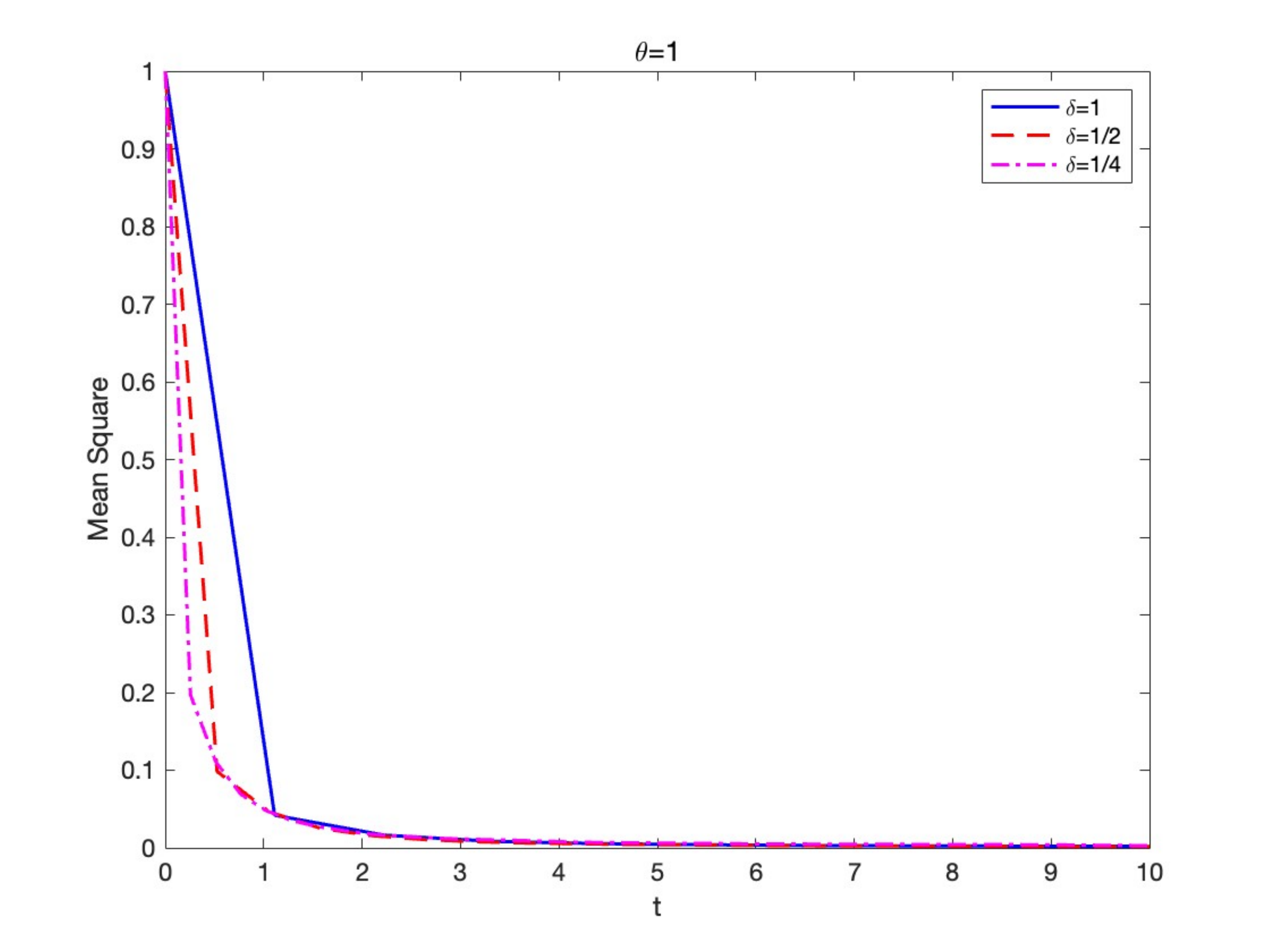}
  \end{subfigure}
      %\vspace{0.5cm}
  \begin{subfigure}[b]{0.49\linewidth}
    \includegraphics[width=\linewidth]{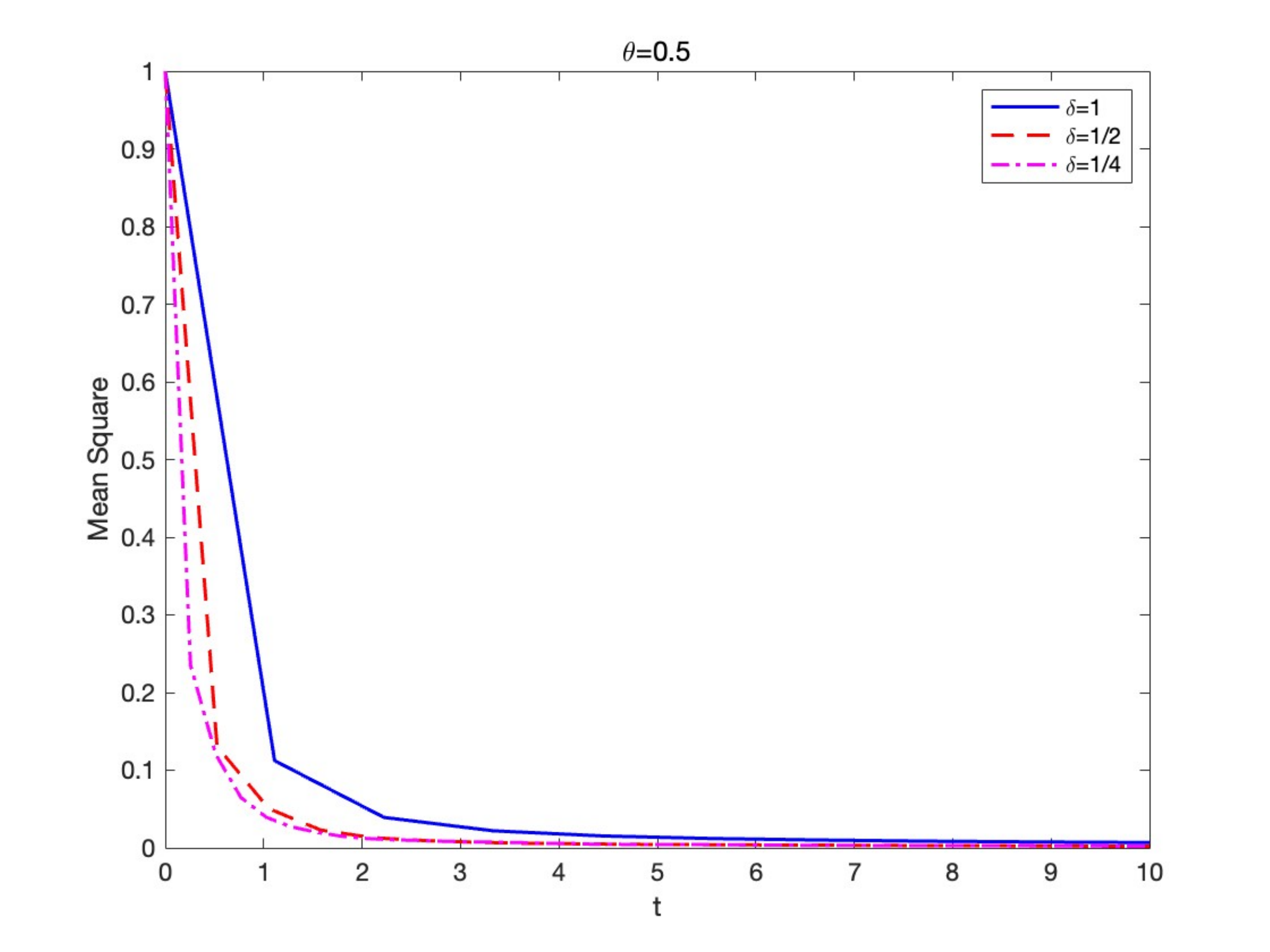}
  \end{subfigure}
  \hfill
  \begin{subfigure}[b]{0.49\linewidth} 
    \includegraphics[width=\linewidth]{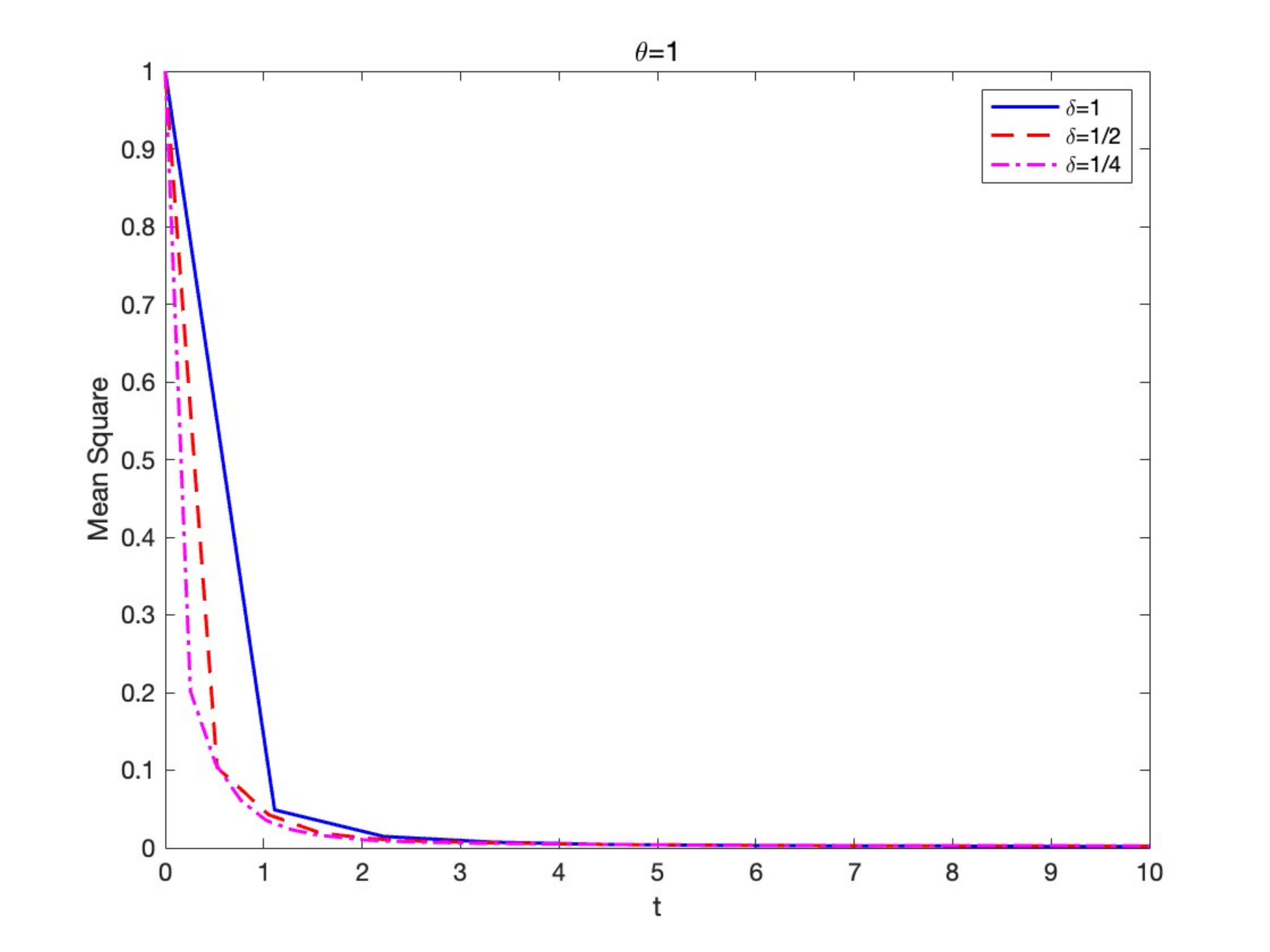}
  \end{subfigure}
  \caption{Mean-square Mittag-Leffler stability of the numerical solutions with different $\theta$ values.(Top: \autoref{Example 7.5}; Middle: \autoref{Example 7.6}; Bottom: \autoref{Example 7.7})}
\end{figure}

\clearpage %put bibliography after figures
\section*{Acknowledgements}
The authors would like to thank the editors and referees for their valuable suggestions. This research did not receive any specific grant from funding agencies in the public, commercial, or not-for-profit sectors.

%\nocite{*}
\bibliographystyle{elsarticle-num}  
\bibliography{main}

\begin{thebibliography}{10}
\expandafter\ifx\csname url\endcsname\relax
  \def\url#1{\texttt{#1}}\fi
\expandafter\ifx\csname urlprefix\endcsname\relax\def\urlprefix{URL }\fi
\expandafter\ifx\csname href\endcsname\relax
  \def\href#1#2{#2} \def\path#1{#1}\fi

\bibitem{magdziarz2011option}
M.~Magdziarz, S.~Orze{\l}, A.~Weron, Option pricing in subdiffusive {Bachelier} model, Journal of Statistical Physics 145 (2011) 187--203.

\bibitem{clark1973subordinated}
P.~K. Clark, A subordinated stochastic process model with finite variance for speculative prices, Econometrica: Journal of the Econometric Society (1973) 135--155.

\bibitem{scalas2006five}
E.~Scalas, Five years of continuous-time random walks in econophysics, in: The Complex Networks of Economic Interactions: Essays in Agent-Based Economics and Econophysics, Springer, 2006, pp. 3--16.

\bibitem{kou2004generalized}
S.~C. Kou, X.~S. Xie, Generalized {Langevin} equation with fractional {Gaussian} noise: Subdiffusion within a single protein molecule, Physical review letters 93~(18) (2004) 180603.

\bibitem{saxton1997single}
M.~J. Saxton, K.~Jacobson, Single-particle tracking: applications to membrane dynamics, Annual Review of Biophysics and Biomolecular Structure 26~(1) (1997) 373--399.

\bibitem{benson2000application}
D.~A. Benson, S.~W. Wheatcraft, M.~M. Meerschaert, Application of a fractional advection-dispersion equation, Water Resources Research 36~(6) (2000) 1403--1412.

\bibitem{umarov2018beyond}
S.~Umarov, M.~Hahn, K.~Kobayashi, Beyond the triangle: \uppercase{B}rownian motion, It{\^o} calculus, and \uppercase{F}okker-\uppercase{P}lanck equation-fractional generalizations, World Scientific, 2018.

\bibitem{hahn2011time}
M.~Hahn, J.~Ryvkina, K.~Kobayashi, S.~Umarov, On time-changed \uppercase{G}aussian processes and their associated \uppercase{F}okker-\uppercase{P}lanck-\uppercase{K}olmogorov equations, Electron. Commun. Probab. 16 (2011) 150--164.

\bibitem{hahn2011fokker}
M.~Hahn, K.~Kobayashi, S.~Umarov, Fokker-\uppercase{P}lanck-\uppercase{K}olmogorov equations associated with time-changed fractional \uppercase{B}rownian motion, Proceedings of the \uppercase{A}merican mathematical Society 139~(2) (2011) 691--705.

\bibitem{meerschaert2009correlated}
M.~M. Meerschaert, E.~Nane, Y.~Xiao, Correlated continuous time random walks, Statistics \& Probability Letters 79~(9) (2009) 1194--1202.

\bibitem{2011Kobayashi}
K.~Kobayashi, Stochastic calculus for a time-changed semimartingale and the associated stochastic differential equations, Journal of Theoretical Probability 24 (2011) 789--820.

\bibitem{hahn2012sdes}
M.~Hahn, K.~Kobayashi, S.~Umarov, \uppercase{SDE}s driven by a time-changed \uppercase{L}{\'e}vy process and their associated time-fractional order pseudo-differential equations, Journal of Theoretical Probability 25 (2012) 262--279.

\bibitem{nane2016stochastic}
E.~Nane, Y.~Ni, Stochastic solution of fractional \uppercase{F}okker-\uppercase{P}lanck equations with space-time-dependent coefficients, Journal of Mathematical Analysis and Applications 442~(1) (2016) 103--116.

\bibitem{magdziarz2016stochastic}
M.~Magdziarz, T.~Zorawik, Stochastic representation of a fractional subdiffusion equation. \uppercase{T}he case of infinitely divisible waiting times, \uppercase{L}{\'e}vy noise and space-time-dependent coefficients, Proceedings of the American Mathematical Society 144~(4) (2016) 1767--1778.

\bibitem{jum2016strong}
E.~Jum, K.~Kobayashi, A strong and weak approximation scheme for stochastic differential equations driven by a time-changed \uppercase{B}rownian motion, Probability and Mathematical Statistics 36 (2016) 201--220.

\bibitem{liu2020truncated}
W.~Liu, X.~Mao, J.~Tang, Y.~Wu, Truncated {Euler}-{Maruyama} method for classical and time-changed non-autonomous stochastic differential equations, Applied Numerical Mathematics 153 (2020) 66--81.

\bibitem{deng2020semi}
C.-S. Deng, W.~Liu, Semi-implicit {Euler}--{Maruyama} method for non-linear time-changed stochastic differential equations, BIT Numerical Mathematics 60 (2020) 1133--1151.

\bibitem{2019JinKobayashi}
S.~Jin, K.~Kobayashi, Strong approximation of stochastic differential equations driven by a time-changed \uppercase{B}rownian motion with time-space-dependent coefficients, Journal of Mathematical Analysis and Applications 476~(2) (2019) 619--636.

\bibitem{jin2021strong}
S.~Jin, K.~Kobayashi, Strong approximation of time-changed stochastic differential equations involving drifts with random and non-random integrators, BIT Numerical Mathematics 61~(3) (2021) 829--857.

\bibitem{wu2016stability}
Q.~Wu, Stability analysis for a class of nonlinear time-changed systems, Cogent Mathematics 3~(1) (2016) 1228273.

\bibitem{zhu2021exponential}
M.~Zhu, J.~Li, D.~Liu, Exponential stability for time-changed stochastic differential equations, Acta Mathematicae Applicatae Sinica, English Series 37~(3) (2021) 617--627.

\bibitem{li2022global}
Z.~Li, L.~Xu, W.~Ma, Global attracting sets and exponential stability of stochastic functional differential equations driven by the time-changed \uppercase{B}rownian motion, Systems \& Control Letters 160 (2022) 105103.

\bibitem{he2024eta}
X.~He, Y.~Zhang, Y.~Wang, Z.~Li, L.~Xu, $\eta$-stability for stochastic functional differential equation driven by time-changed \uppercase{B}rownian motion, Journal of Inequalities and Applications 2024~(1) (2024) 60.

\bibitem{zhang2019razumikhin}
X.~Zhang, C.~Yuan, Razumikhin-type theorem on time-changed stochastic functional differential equations with {Markovian} switching, Open Mathematics 17~(1) (2019) 689--699.

\bibitem{zhang2021asymptotic}
X.~Zhang, Z.~Zhu, C.~Yuan, Asymptotic stability of the time-changed stochastic delay differential equations with {Markovian} switching, Open Mathematics 19~(1) (2021) 614--628.

\bibitem{shen2023class}
G.~Shen, T.~Zhang, J.~Song, J.-L. Wu, On a class of distribution dependent stochastic differential equations driven by time-changed \uppercase{B}rownian motions, Applied Mathematics \& Optimization 88~(2) (2023) 33.

\bibitem{nane2017stability}
E.~Nane, Y.~Ni, Stability of the solution of stochastic differential equation driven by time-changed \uppercase{L}{\'e}vy noise, Proceedings of the American Mathematical Society 145~(7) (2017) 3085--3104.

\bibitem{yin2021stability}
X.~Yin, W.~Xu, G.~Shen, Stability of stochastic differential equations driven by the time-changed \uppercase{L}{\'e}vy process with impulsive effects, International Journal of Systems Science 52~(11) (2021) 2338--2357.

\bibitem{xu2023h}
L.~Xu, Z.~Li, B.~Huang, $ h $-stability for stochastic functional differential equation driven by time-changed \uppercase{L}{\'e}vy process, AIMS Mathematics 8~(10) (2023) 22963--22983.

\bibitem{mao2013strong}
X.~Mao, L.~Szpruch, Strong convergence and stability of implicit numerical methods for stochastic differential equations with non-globally {Lipschitz} continuous coefficients, Journal of Computational and Applied Mathematics 238 (2013) 14--28.

\bibitem{zhang2022convergence}
Y.~Zhang, M.~Song, M.~Liu, Convergence and stability of stochastic theta method for nonlinear stochastic differential equations with piecewise continuous arguments, Journal of Computational and Applied Mathematics 403 (2022) 113849.

\bibitem{meerschaert2004limit}
M.~M. Meerschaert, H.-P. Scheffler, Limit theorems for continuous-time random walks with infinite mean waiting times, Journal of Applied Probability 41~(3) (2004) 623--638.

\bibitem{bingham1971limit}
N.~Bingham, Limit theorems for occupation times of {Markov} processes, Zeitschrift f{\"u}r Wahrscheinlichkeitstheorie und verwandte Gebiete 17 (1971) 1--22.

\bibitem{Protter2004}
P.~Protter, Stochastic Integration and Differential Equations, 2nd Edition, Springer, 2004.

\bibitem{mao1991semimartingales}
X.~Mao, Stability of stochastic differential equations with respect to semimartingales, Longman, 1991.

\bibitem{maobook}
X.~Mao, Stochastic Differential Equations and Applications, 2nd Edition, Elsevier, 2011.

\bibitem{Jacod2003}
J.~Jacod, A.~N. Shiryaev, Limit Theorems for Stochastic Processes, Vol. 288 of Grundlehren der mathematischen Wissenschaften, Springer, Berlin, 2003.

\bibitem{chambers1976method}
J.~M. Chambers, C.~L. Mallows, B.~Stuck, A method for simulating stable random variables, Journal of the American Statistical Association 71~(354) (1976) 340--344.

\bibitem{magdziarz2009stochastic}
M.~Magdziarz, Stochastic representation of subdiffusion processes with time-dependent drift, Stochastic Processes and their Applications 119~(10) (2009) 3238--3252.

\bibitem{zeidler2013nonlinear}
E.~Zeidler, Nonlinear Functional Analysis and its Applications: II/B: Nonlinear Monotone Operators, Springer Science \& Business Media, 2013.

\bibitem{mao2013strongBEM}
X.~Mao, L.~Szpruch, Strong convergence rates for backward \uppercase{E}uler--\uppercase{M}aruyama method for non-linear dissipative-type stochastic differential equations with super-linear diffusion coefficients, Stochastics: An International Journal of Probability and Stochastic Processes 85~(1) (2013) 144--171.

\bibitem{Karatzas1998}
I.~Karatzas, S.~Shreve, Brownian Motion and Stochastic Calculus, 2nd Edition, Springer, New York, 1998.

\bibitem{jum2015numerical}
E.~Jum, Numerical approximation of stochastic differential equations driven by \uppercase{L}{\'e}vy motion with infinitely many jumps, Ph.D. thesis, University of Tennessee - Knoxville (2015).

\end{thebibliography}

\end{document}